\documentclass[11pt,times]{article}

\usepackage[page-ref=none]{acro}

\usepackage{algorithm}
\usepackage{algorithmic}

\usepackage[sumlimits]{amsmath}
\usepackage{amssymb}
\usepackage{amsthm}

\usepackage[backend=bibtex,citestyle=numeric,bibstyle=numeric,maxcitenames=3,maxbibnames=5,giveninits=true]{biblatex}

\DeclareNameAlias{default}{last-first}

\usepackage{color}
\usepackage{enumitem}

\usepackage[mathscr]{eucal}
\usepackage{fancybox}
\usepackage{float}

\usepackage[text={6.5in,9.5in},centering]{geometry}

\usepackage{graphics}
\usepackage{graphicx}

\usepackage[hidelinks]{hyperref}

\usepackage{ifpdf}
\usepackage{ifthen}

\usepackage{listings}
\usepackage{pgfplots}

\usepackage{stackengine}

\usepackage{stmaryrd}
\usepackage{subfig}

\usepackage{tikz}
\usetikzlibrary{arrows}
\usetikzlibrary{fadings}

\usepackage{tikz-cd}

\usepackage{txfonts}    

\usepackage{url}
\usepackage{varioref}
\usepackage{xcolor}

\DeclareAcronym{AMP}{
  short = AMP ,
  long  = asymmetric multicore processor
}

\DeclareAcronym{BVP}{
  short = BVP ,
  long  = boundary value problem
}

\DeclareAcronym{CPPO}{
  short = cppo ,
  long  = complete pointed partial order
}

\DeclareAcronym{DAG}{
  short = DAG ,
  long  = directed acyclic graph
}

\DeclareAcronym{DCPO}{
  short = dcpo ,
  long  = directed-complete partial order
}

\DeclareAcronym{DBF}{
  short = DBF ,
  long  = Digital Beam Former
}

\DeclareAcronym{DSL}{
  short = DSL ,
  long  = domain-specific language
}

\DeclareAcronym{EDSL}{
  short = EDSL ,
  long  = embedded domain-specific language
}

\DeclareAcronym{FBP}{
  short = FBP ,
  long  = free boundary problem
}

\DeclareAcronym{FivePs}{
  short = 5Ps ,
  long  = five Ps
}

\DeclareAcronym{HS}{
  short = HS ,
  long  = hybrid system
}

\DeclareAcronym{ISA}{
  short = ISA ,
  long  = instruction set architecture
}

\DeclareAcronym{IVP}{
  short = IVP ,
  long  = initial value problem
}

\DeclareAcronym{ODE}{
  short = ODE ,
  long  = ordinary differential equation
}

\DeclareAcronym{PC}{
  short = PC ,
  long  = Pulse Compression
}

\DeclareAcronym{PCROP}{
  short = PC-ROP ,
  long  = PDE-constrained rearrangement optimization problem
}

\DeclareAcronym{PDE}{
  short = PDE ,
  long  = partial differential equation
}

\DeclareAcronym{PI}{
  short = PI ,
  long  = principal investigator
}

\DeclareAcronym{POSET}{
  short = poset ,
  long  = partially ordered set
}

\DeclareAcronym{RE}{
  short = r.~e. ,
  long  = recursively enumerable
}

\DeclareAcronym{ROP}{
  short = ROP ,
  long  = rearrangement optimization problem
}

\DeclareAcronym{SAC}{
  short = \textsc{Sac} ,
  long  = Single Assignment C
}

\DeclareAcronym{SBSE}{
  short = SBSE ,
  long  = search-based software engineerign
}

\DeclareAcronym{SOP}{
  short = SOP ,
  long  = shape optimization problem
}

\DeclareAcronym{TTE}{
  short = TTE ,
  long  = Type-II Theory of Effectivity
}

\DeclareAcronym{UIC}{
  short = UIC ,
  long  = uniform ideal continuity
}

\newcommand{\defeq}{\coloneqq} 

\newcommand{\PS}{\mathscr{P}}

\newcommand{\cCPPO}[1]{\mathbf{C}#1} 

\newcommand{\kCPPO}[1]{\mathbf{K}#1} 

\newcommand{\glb}{\ensuremath{\bigsqcap}} 
\newcommand{\lub}{\ensuremath{\bigsqcup}} 

\newcommand{\Po}{\mathbf{Po}}

\newcommand{\id}{\text{id}}

\newcommand{\cpo}[1][s]{\ifthenelse{\equal{#1}{s}}{complete partial
order}{Complete partial order}}

\newcommand{\eg}{e.\,g.}
\newcommand{\ie}{i.\,e.}

\newcommand{\absn}[1]{\lvert\thinspace {#1} \thinspace\rvert}

\newcommand{\D}{\ensuremath{\mathbb{D}}}

\newcommand{\fpoint}[1][s]{\ifthenelse{\equal{#1}{c}}{Floating-point}{floating-point}}

\newcommand{\hrEnclos}[1]{\ensuremath{#1^\Box}} 

\newcommand{\interiorOf}[1]{\ensuremath{{#1}^{\circ}}}

\newcommand{\intvaldom}[1][{\Rinf}]{\ensuremath{{\mathbf{I}#1}}}
\newcommand{\intvalfun}[1][f]{\ensuremath{{\boldsymbol{I}}#1}}

\newcommand{\kc}[1][s]{\ifthenelse{\equal{#1}{s}}{Kolmogorov
complexity}{Kolmogorov complexities}} 

\newcommand{\lep}[1][x]{\ensuremath{{\underline{#1}}}} 
\newcommand{\uep}[1][x]{\ensuremath{{\overline{#1}}}} 

\newcommand{\lift}{_\bot} 

\newcommand{\lowerSet}[1]{\ensuremath{\downarrow \negmedspace #1}}
\newcommand{\upperSet}[1]{\ensuremath{\uparrow \negmedspace #1}}

\newcommand{\ptime}[1][s]{\ifthenelse{\equal{#1}{c}}{Polynomial-time}{polynomial-time}}
\newcommand{\pspace}[1][s]{\ifthenelse{\equal{#1}{c}}{Polynomial-space}{polynomial-space}}

\newcommand{\md}{\thinspace \mathrm{d}} 
\newcommand{\me}{\mathrm{e}} 

\newcommand{\N}{\ensuremath{\mathbb{N}}}

\newcommand{\NFL}[1][c]{\ifthenelse{\equal{#1}{s}}{no free lunch}{No
Free Lunch}}
\newcommand{\norm}[1]{\ensuremath{{\lVert\thinspace {#1}
\thinspace\rVert}}} 
\newcommand{\ntuple}[3]{\ensuremath{\left( {#1}_{#2}, \ldots, {#1}_{#3} \right)}}

\newcommand{\Q}{\ensuremath{\mathbb{Q}}}

\newcommand{\R}{\ensuremath{\mathbb{R}}}

\newcommand{\RId}[1]{\ensuremath{\mathrm{RId}(#1)}} 

\newcommand{\restrictTo}[2]{\ensuremath{#1 \negmedspace \restriction_{#2}}} 

\newcommand{\Scott}[1]{\ensuremath{\Sigma ({#1})}} 

\newcommand{\sequence}[5][i]{\ifthenelse{\equal{#1}{i}}{\ensuremath{\left< {#2}_{#3}, {#2}_{#4},
\ldots,{#2}_{#5}, \ldots \right>}}{\ensuremath{\left< {#2}_{#3}, {#2}_{#4},
\ldots,{#2}_{#5} \right>}}}

\newcommand{\set}[1]{\ensuremath{\left\{{#1}\right\}}}
\newcommand{\setbarTall}[2]{\ensuremath{\left\{{#1}\;\vrule\;
{#2}\right\}}} 
\newcommand{\setbarNormal}[2]{\ensuremath{\{{#1} \mid
{#2}\}}}

\newcommand{\symExpand}{\ensuremath{\oplus}} 

\newcommand{\wayAboves}[1]{\ensuremath{\Uparrow \negmedspace #1}}
\newcommand{\wayBelows}[1]{\ensuremath{\Downarrow \negmedspace #1}}

\theoremstyle{plain} 

\newtheorem{definition}{Definition}[section]
\newtheorem{assumption}[definition]{Assumption}
\newtheorem{corollary}[definition]{Corollary}
\newtheorem{lemma}[definition]{Lemma}
\newtheorem{notation}[definition]{Notation}
\newtheorem{proposition}[definition]{Proposition}
\newtheorem{theorem}[definition]{Theorem}

\theoremstyle{definition}

\newtheorem{example}[definition]{Example}
\newtheorem{remark}[definition]{Remark}

\newcounter{defenumalph}

\newcounter{defenum}

\floatstyle{ruled} 
\restylefloat{figure}
\restylefloat{table}

\newcounter{saveeqn}

\newcommand{\arrayoptions}[2]{\setlength{\arraycolsep}{#1}\renewcommand{\arraystretch}{#2}}

\DeclareMathOperator{\atan}{atan}

\DeclareMathOperator{\fix}{fix} 

\DeclareMathOperator{\precAbsBas}{\vartriangleleft}

\DeclareMathOperator{\widthOf}{w} 

\title{\textbf{Recursive Solution of Initial Value Problems with Temporal Discretization}}

\author{Abbas Edalat\thanks{Department of Computing, Imperial College,
    180 Queen’s Gate, London, SW7 2BZ, United Kingdom, Email:
    \href{mailto:A.Edalat@ic.ac.uk} {\texttt{A.Edalat@ic.ac.uk}}} \and
  Amin Farjudian\thanks{\textbf{(Corresponding author)} School of Mathematics, University of Birmingham, Edgbaston, Birmingham, B15 2TT, United Kingdom, Email:
    \href{mailto:A.Farjudian@bham.ac.uk}
    {\texttt{A.Farjudian@bham.ac.uk}}} \and Yiran Li\thanks{School
    of Computer Science, University of Nottingham Ningbo China, 199
    Taikang East Road, Ningbo, 315100, Zhejiang, China, Email:
    \href{mailto:Yiran.Li@nottingham.edu.cn}
    {\texttt{Yiran.Li@nottingham.edu.cn}}} }

\date{}

\begin{document}

\maketitle

\begin{abstract}

  We construct a continuous domain, as a model of interval analysis,
  for temporal discretization of differential equations. By using this
  domain, and the domain of Lipschitz maps, we formulate a
  generalization of the Euler operator, which exhibits second-order
  convergence. We prove computability of the operator within the
  framework of effectively given domains. The operator only requires
  the vector field of the differential equation to be Lipschitz
  continuous, in contrast to the related operators in the literature
  which require the vector field to be at least continuously
  differentiable. Within the same framework, we also analyze temporal
  discretization and computability of another variant of the Euler
  operator formulated according to Runge-Kutta theory. We prove that,
  compared with this variant, the second-order operator that we
  formulate directly, not only imposes weaker assumptions on the
  vector field, but also exhibits superior convergence rate. We
  implement the first-order, second-order, and Runge-Kutta Euler
  operators using arbitrary-precision interval arithmetic, and report
  on some experiments. The experiments confirm our theoretical
  results. In particular, we observe the superior convergence rate of
  our second-order operator compared with the Runge-Kutta Euler and
  the common (first-order) Euler operators.
  
\end{abstract}


\noindent
\textbf{MSC classes:}
\textbf{06B35}, 
\textbf{65G20}, 
68Q55, 
65L05, 
65L20 

\vspace{1em}
\noindent
\textbf{Declarations of interest:} none

\tableofcontents

\section{Introduction}
\label{sec:intro}

Many concepts and phenomena in modern science have been formulated
using differential equations. To analyze systems which have been
formulated using differential equations, apart from rare cases with
closed-form solutions, numerical methods must be used. The methods of
classical numerical analysis are adequate when infrequent errors are
acceptable. In safety-critical systems, however, even isolated errors
can lead to disproportionate damages, or more importantly, loss of
life. Therefore, in analyses of such systems, computations must be
\emph{validated}.

A powerful approach to validated computation is provided by interval
analysis~\parencite{Moore:2009:IIA,Alefeld:Interval_analysis:2000}. The
foundations of modern interval analysis were laid in 1950s by
\textcite{Warmus-Intval_Arith-1956},
\textcite{Sunaga:interval_algebra:1958}, and
\textcite{Moore:Automatic_Error_Analysis:1959}, although, the most
influential work has been the well-known book by
\textcite{Moore:1966:IA}. Ever since, interval methods have evolved
considerably, and have been applied to various tasks, {\eg}, interval
Newton method~\parencite{Moore_Jones:Safe_Starting_Regions:1977},
optimization~\parencite{Hansen-Global2004}, solution of systems of
linear
equations~\parencite{Minamihata_et_al:interval_linear_systems:2020},
eigenvalue problems~\parencite{Hladik_et_al:interval_eigenvalue:2011},
solution of
\acp{PDE}~\parencite{Schwandt:Almost_Globally_Convergent:1986,Hoffmann_Marciniak_Szyszka:Interval_Central_Difference:2013,Hoffmann_Marciniak:Interval_generalized_Poisson:2016},
various topics in functional
analysis~\parencite{Moore_Cloud:Computational_Functional_Analysis:2007},
and security analysis of neural
networks~\parencite{Wang_et_al:symbolic_intervals:2018}, to name just
a few. There is also a rich literature on the solution of
\acp{IVP}---which is the subject of the current article---using
interval methods, {\eg},
\parencite{Berz_Makino:VerifiedIO:1998,Nedialkov_et_al:High_Order_Interval_ODE:2001,RauhValencia-SCAN.2006.47,Marciniak:Selected_Interval_Methods:2009,Moore:2009:IIA,Tucker:Validated_Num:Book:2011,AlexandreDitSandretto:Validated_Runge_Kutta:2016,Marciniak_Szyszka:Interval_Runge-Kutta:2019}.

Interval analysis has also had a productive interplay with domain
theory~\parencite{AbramskyJung94-DT,Gierz-ContinuousLattices-2003,Goubault-Larrecq:Non_Hausdorff_topology:2013},
and each subject has enriched the other. In Dana Scott's seminal paper
on domain theory~\parencite{Scott:Outline:1970}, the interval lattice
appears as one of the only two examples of concrete data types that
are presented.\footnote{With the lattice of natural numbers being the
  other example.}  In the
monograph~\parencite{Scott:Outline:Oxford_Lab_Version:1970}---which is
a revised and expanded version
of~\parencite{Scott:Outline:1970}---Moore's book on interval
analysis~\parencite{Moore:1966:IA} is mentioned as one of the few
references related to the mathematical theory that is
developed. Subsequent results in casting mathematical analysis in a
domain-theoretic framework also owe much of their inspiration and
conceptual foundation to the earlier work on interval analysis. For
instance, interval methods for integration and differential equation
solving appear in~\parencite{Moore:1966:IA}, well before they were
recast in a domain-theoretic
framework~\parencite{Edalat_Escardo:Integ_realPCF:2000,EdalatPattinson2006-Euler-PARA,Edalat_Pattinson2007-LMS_Picard}. Domain
theory, in turn, provides a refined framework for interval analysis,
which enables systematic analyses of convergence, computability, and
complexity of the methods involved. To be more concrete, domain theory
provides a model of interval analysis whereby a notion of
\emph{approximation} of an interval $I \defeq [\lep[I],\uep[I]]$ by
another interval $J \defeq [\lep[J],\uep[J]]$ is introduced, which
requires $J$ to contain $I$ in its interior, {\ie},
$[\lep[I],\uep[I]] \subseteq (\lep[J],\uep[J])$.\footnote{Also,
  see~\eqref{eq:way_below_compact_sets} on
  page~\pageref{eq:way_below_compact_sets}.} This notion, together
with the notion of convergence of a nested set of intervals to its
intersection, allows for an analysis of computability of (say) the
solution of an \ac{IVP}, bringing the study of \acp{ODE} and theories
of computability together.

In the current article, we focus on \acp{ODE}. Specifically, we
consider the \ac{IVP}:

\begin{equation}
  \label{eq:main_ivp}
    \left\{
      \begin{array}{r@{\hspace{0.5ex}=\hspace{0.5ex}}l}
        y'(t) & f(y(t)),\\
        y(0) & (\underbrace{0, \ldots, 0}_{n}),\\
      \end{array}
    \right.
  \end{equation}
  in which $f: [-K,K]^n \to [-M,M]^n$ is a continuous vector field,
  for some natural number $n \geq 1$, and positive rational numbers
  $K, M \in \Q_+$. 
  
  By merely assuming the vector field $f$ to be continuous, Peano's
  theorem implies that the \ac{IVP}~\eqref{eq:main_ivp} has a
  differentiable solution
  $y : [ -\frac{K}{M}, \frac{K}{M}] \to [ -K,
  K]^n$~\cite[Theorem~2.19]{Teschl:2012:ODE}. Furthermore, if $f$ is
  Lipschitz continuous, by Picard-Lindel\"{o}f theorem, the \ac{IVP}
  must have a unique solution~\cite[Theorem~2.2]{Teschl:2012:ODE}. The
  Lipschitz condition is sufficient but not necessary for uniqueness
  of solutions~\parencite{MillerMichel-ODE-1982}. Hiroshi Okamura
  provided a necessary and sufficient condition for uniqueness of
  solutions, which may be found
  in~\parencite{YoshizawaHayashi-Uniqueness-ODE-1950}.
  
  Lipschitz continuity is, nevertheless, the common assumption that is
  imposed on the vector field in developing computational methods. This is
  partly because convergence analysis is easier using Lipschitz
  properties of the maps involved. We also focus on Lipschitz
  continuous functions due to their desirable computational
  properties~\parencite{Edalat_Lieutier:Domain_Calculus_One_Var:MSCS:2004,EdalatLieutierPattinson:2013-MultiVar-Journal,Edalat:Lipschitz:JACM:2022}.

  Methods of classical numerical analysis for solving differential
  equations are commonly implemented using floating-point
  numbers. These implementations have inherent inaccuracies, mainly
  due to round-off and truncation errors. A full error analysis may be
  carried out in some
  cases~\parencite{Higham:2002:ASN,Boldo_Melquiond:Computer_Arithmetic:2017},
  but in general, error analysis of floating-point computations is
  very hard, and the computations may be deceptively unstable, even
  when the operations involved are rather
  simple~\parencite{Muller:Arithmetique_des_ordinateurs:Book:1989,Menissier2005,Rump:recurrences:2020,Rump:addendum_recurrences:2020}. An
  effective alternative is developing algorithms which are correct
  \emph{by construction}, {\ie}, do not require a separate error
  analysis. Validated numerics offer a sound alternative of this kind,
  where results are provided with absolute guarantees of correctness.

By assuming the vector field to be analytic, a validated solution of
the \ac{IVP}~(\ref{eq:main_ivp}) can be obtained in polynomial
time~\parencite{Mueller-Poly_IVP:1993} by using a Taylor model. Without the
analyticity assumption, \ac{IVP} solving becomes PSPACE-complete, even
with $C^\infty$ assumption on the vector
field~\parencite{Kawamura-Lipschitz_IVP-2009}. Despite such complexity
constraints, for sufficiently regular vector fields, there are
feasible validated methods available in the literature, including
high-order Taylor
methods~\parencite{Berz_Makino:VerifiedIO:1998,Nedialkov_et_al:High_Order_Interval_ODE:2001,Tucker:Validated_Num:Book:2011}
and Runge-Kutta
methods~\parencite{Marciniak:Selected_Interval_Methods:2009,AlexandreDitSandretto:Validated_Runge_Kutta:2016,Marciniak_Szyszka:Interval_Runge-Kutta:2019}.

By combining powerful tools from order theory and topology, domain
theory provides a rich semantic model for validated numerics. In
contrast to Taylor methods, a characteristic feature of a
domain-theoretic approach is that the vector field is interpreted as
the limit (under Scott topology) of a sequence of
finitely-representable approximations. An added advantage of
developing differential equation solvers in a domain-theoretic
framework is that they can be subsequently incorporated into
domain-theoretic models for other applications, such as reachability
analysis of hybrid
systems~\parencite{Edalat_Pattinson:Hybrid:2007,Moggi_Farjudian_Duracz_Taha:Reachability_Hybrid:2018,Moggi_Farjudian_Taha:System_Analysis_and_Robustness:2019,Moggi_Farjudian_Taha:System_Analysis_and_Robustness:ICTCS:2019}.

Domain theoretic methods for solving \acp{IVP} have indeed been
studied for Picard
method~\parencite{Edalat_Lieutier:Domain_Calculus_One_Var:MSCS:2004,Edalat_Pattinson2007-LMS_Picard},
(first-order) Euler method~\parencite{EdalatPattinson2006-Euler-PARA},
and a second-order variant of Euler
method~\parencite{Edalat_Farjudian_Mohammadian_Pattinson:2nd_Order_Euler:2020:Conf}.\footnote{A
  numerical scheme for solving \acp{IVP} is said to be of order $p$
  if, for every step size $h$, the local error incurred is of order
  $O(h^{p+1})$~\cite[Page 8]{Iserles:2009}. See
  Section~\ref{subsec:convergence_analysis} for more details.} These
methods are all based on the interval domain
model~\parencite{Escardo96-tcs}. The Picard operators
of~\parencite{Edalat_Lieutier:Domain_Calculus_One_Var:MSCS:2004,Edalat_Pattinson2007-LMS_Picard}
have a functional style of definition. As the underlying interval
domain models may be enriched with an effective structure, it is
possible to study computability of the Picard operator using the
classical interval domain. There are several advantages in adopting a
functional---as opposed to imperative---style of programming, such as
correctness of program construction, equational reasoning, modularity,
and automatic optimization, to name a
few~\parencite{Backus:Von_Neumann_Style:1978,Hughes:Why_FP_Matters:1989,Hu_Hughes_Wang:How_FP_mattered:2015}. In
general, functional programs are significantly more amenable to
automatic reasoning and verification, compared with imperative
programs. These are particularly relevant when programs become larger
and more complex, as in the context of \ac{IVP} solving.

The Euler operators
of~\parencite{EdalatPattinson2006-Euler-PARA,Edalat_Farjudian_Mohammadian_Pattinson:2nd_Order_Euler:2020:Conf},
however, have been defined in an imperative style. As it turns out,
this is not a problem that can be rectified directly in the classical
interval domain models of
(say)~\parencite{Escardo96-tcs,Edalat_Lieutier:Domain_Calculus_One_Var:MSCS:2004,Edalat_Pattinson2007-LMS_Picard}.\footnote{These
  claims will be justified formally in
  Section~\ref{sec:domain_temporal_discretization}, after introducing
  the necessary technical background in
  Section~\ref{sec:preliminaries}.} In essence, the main difference
between Euler and Picard methods is that Euler methods proceed
according to a \emph{temporal discretization}. This temporal
discretization is the main obstacle which renders classical interval
domain models ineffective. Runge-Kutta methods also proceed according
to a temporal discretization. As a result, the classical interval
domain is not adequate for functional definition of Runge-Kutta
operators either.

\subsection{Contributions}
\label{sec:contributions}

At a fundamental level, the main contribution of this article is the
construction of an effectively-given domain for temporal
discretization of differential equations. The construction of this
domain model requires a novel approach because the classical interval
domain models are too restrictive. For further clarification, we
expand on this claim informally here, while the technical details and
precise formulations will be presented in
Section~\ref{sec:domain_temporal_discretization}.

We let $\R$ denote the set of real numbers, and let
$(\intvaldom[\R], \sqsubseteq)$ denote the interval domain, {\ie}, the
\ac{POSET} $\setbarNormal{[a,b]}{a,b \in \R, a \leq b}$ of non-empty
compact intervals ordered under reverse inclusion:
\begin{equation*}
  [a,b] \sqsubseteq [c,d] \iff [c,d] \subseteq [a,b].
\end{equation*}
\noindent
For every interval function $f: \R \to \intvaldom[\R]$ , we let
$\lep[f]$ and $\uep[f]$ denote the lower and upper bounds of $f$,
{\ie}, $\forall x \in \R: f(x) = [\lep[f](x), \uep[f](x)]$. It is
well-known that $f$ is continuous with respect to the Euclidean
topology over $\R$ and Scott topology over $\intvaldom[\R]$ if and
only if $\lep[f]: \R \to \R$ is lower semi-continuous and
$\uep[f]: \R \to \R$ is upper
semi-continuous~\parencite{Edalat_Lieutier:Domain_Calculus_One_Var:MSCS:2004}. For
temporal discretization of differential equations, these requirements
turn out to be too restrictive. A typical method that proceeds by
temporal discretization---{\eg}, Euler or Runge-Kutta method---may
only require semi-continuity \emph{from the
  left}~(Definition~\ref{def:left_semi_continuity}). The problem is
that, the \ac{POSET} of interval functions with left semi-continuous
bounds---which we denote by ${\mathcal{D}}_{\Q}$ and define formally
in~\eqref{eq:def_D_Q}---is not even
continuous~(Corollary~\ref{cor:D_Q_not_cont}).\footnote{The subscript
  $\Q$ appears in ${\mathcal{D}}_{\Q}$ because we take the partition
  points to be rational numbers.} This should be contrasted with the
\ac{POSET} of interval functions with semi-continuous bounds, which is
indeed a continuous domain, and has been used extensively, {\eg}, in
domain-theoretic analysis of the Picard
method~\parencite{Edalat_Lieutier:Domain_Calculus_One_Var:MSCS:2004,Edalat_Pattinson2007-LMS_Picard}.

To address this problem, we present a method for constructing
continuous domains for function spaces $[X \to D]$, in which $X$ is an
arbitrary topological space---{\ie}, not necessarily
core-compact---and $D$ is an arbitrary bounded-complete continuous
domain. Then, following this general construction, and using a
suitable abstract basis, we obtain the fundamental $\omega$-continuous
domain of the paper as a special case, which we denote by
${\mathcal{W}}_{\mathcal{D}}$ and define formally in
Definition~\ref{def:W_Q}. Although this domain is not isomorphic to
${\mathcal{D}}_{\Q}$, they are related to each other via a Galois
connection~(which is a special case of the general Galois connection
of Theorem~\ref{thm:Galois_connection_X_to_D}). As such,
${\mathcal{W}}_{\mathcal{D}}$ may be regarded as an `optimal'
substitute for the non-continuous \ac{POSET}~${\mathcal{D}}_{\Q}$.

Suitability of ${\mathcal{W}}_{\mathcal{D}}$ for temporal discretization will be
demonstrated by developing domain-theoretic formulation of two
operators for validated solution of \acp{IVP}:

\begin{enumerate}[label=(\arabic*)]
\item A second-order Euler operator $E^2$, which has its foundation in the
  results
  of~\parencite{Edalat_Farjudian_Mohammadian_Pattinson:2nd_Order_Euler:2020:Conf}.
\item An Euler operator $E^{\mathrm{R}}$, which is formulated
  according to Runge-Kutta theory.
\end{enumerate}

We will show that, the operator $E^2$ has superior theoretical and
practical properties compared with the operator $E^{\mathrm{R}}$. In
particular, it exhibits superior convergence properties under weaker
differentiability assumptions on the vector field of the differential
equation.

In the formulation of Picard and Euler operators as presented
in~\parencite{Edalat_Lieutier:Domain_Calculus_One_Var:MSCS:2004,Edalat_Pattinson2007-LMS_Picard}
and~\parencite{EdalatPattinson2006-Euler-PARA}, respectively, the
vector field is required to be Lipschitz continuous to guarantee uniqueness
of the solutions. The local Lipschitz properties of the vector field,
however, are not used in the algorithms. The second-order Euler method
introduced
in~\parencite{Edalat_Farjudian_Mohammadian_Pattinson:2nd_Order_Euler:2020:Conf}
uses the local Lipschitz properties of the vector field to speed up the
convergence. In~\parencite{Edalat_Farjudian_Mohammadian_Pattinson:2nd_Order_Euler:2020:Conf},
several fundamental properties of the method have been proven, such as
soundness and completeness, together with some basic convergence and
algebraic complexity results. The convergence analysis
of~\parencite{Edalat_Farjudian_Mohammadian_Pattinson:2nd_Order_Euler:2020:Conf},
however, is not fine enough to reflect the second-order nature of the
operator. As such, another important contribution of the current
article is a detailed convergence analysis of the Euler operator
$E^2$, which proves that it is truly second-order.

To summarize, the main contributions of the current article are as
follows:

\begin{enumerate}[label=(\roman*)]

\item A general construction that provides continuous domains for
  function spaces $[X \to D]$, for arbitrary topological spaces $X$
  and arbitrary bounded-complete continuous domains $D$.
  
\item Introduction of a continuous domain for temporal discretization
  and solution of \acp{IVP}, as a special case of the general
  construction.

  \item Domain-theoretic formulations of the Euler operator according
    to two approaches: a direct approach, and one which conforms to
    the Runge-Kutta theory.

  \item Computable analysis of the Euler operators within the
    framework of effectively given domains.

  \item A detailed convergence analysis of the Euler operator $E^2$,
    which shows that the operator is indeed second-order.

  \item Experiments which confirm the theoretical results.
    
  \end{enumerate}
  
  \begin{remark}[Initial values]
    For simplicity, we take the initial values in the
    \ac{IVP}~(\ref{eq:main_ivp}) to be all zeros. The results can be
    generalized, in a straightforward manner, to initial conditions of
    the form $y( t_0) = (q_1, \ldots, q_n)$, in which
    $t_0, q_1, \ldots, q_n \in \Q$. Indeed, in our experiments, we
    consider \acp{IVP} that have rational initial values. Generalizing
    further to initial values which are irrational points, or
    non-degenerate intervals, is not so straightforward, and must be
    studied separately. We do not follow this line in the current
    work, but point out that handling uncertain initial values has
    been studied for the Picard
    operator~\parencite{Konecny_Duracz_Farjudian_Taha:Picard_Uncertain:2014}.
  \end{remark}

\begin{remark}[Autonomous versus non-autonomous \acp{IVP}]
  We consider only autonomous \acp{IVP} of the
  form~(\ref{eq:main_ivp}). Given a non-autonomous equation
  $y'(t)=f(t,y(t))$ with $y(t) = (y_1(t), \ldots, y_n(t))$, the
  function defined by $Y(t) \defeq (t, y(t))$ satisfies the autonomous
  equation $Y'(t)=g(Y(t))$ with
  $g( t, \vec{\theta}) \defeq ( 1, f( t, \vec{\theta}))$. Thus, any
  non-autonomous \ac{IVP} can be converted to an autonomous one by
  adding an extra variable. The uniqueness of solutions for the
  non-autonomous \ac{IVP} is guaranteed if the vector field $f$ is
  Lipschitz continuous in its second argument, in contrast to the
  autonomous \ac{IVP}~(\ref{eq:main_ivp}), for which Lipschitz
  continuity of $f$ is required in all arguments. As such, we lose a
  slight bit of generality.
\end{remark}


  \subsection{Structure of the Paper}
\label{sec:structure}

The rest of this article is structured as follows:

\begin{itemize}
\item The necessary technical background and notation is introduced in
  Section~\ref{sec:preliminaries}.

\item The general construction of the main continuous domain is
  presented in Section~\ref{sec:domain_for_fun_spaces}.

\item In Section~\ref{sec:domain_temporal_discretization}, we
  investigate the continuous domain constructed in
  Section~\ref{sec:domain_for_fun_spaces} for the specific case of
  temporal discretization.

\item Section~\ref{susec:2nd_Order_Euler_Operator} contains the
  domain-theoretic analyses of the second-order Euler operator,
  including the formulation of the operator, computability, and
  convergence analyses.

\item The domain-theoretic analyses of the Runge-Kutta Euler operator
  are presented in Section~\ref{sec:Runge_Kutta_Euler_Operator}.

\item The results of some of our experiments will be presented in
  Section~\ref{sec:Experiments}, where different variants of the Euler
  operator will be compared in terms of their performance.

\item The article will be concluded by some remarks in Section~\ref{sec:concluding_remarks}.

\end{itemize}

\section{Preliminaries}
\label{sec:preliminaries}

In this section, we present the technical background needed for the
rest of the article.

\subsection{Domain Theory}
\label{subsec:domain_theory}

Domain theory has its roots in topological
algebra~\parencite{Compendium:Book:1980,Keimel:Domain_Ramifications_Interactions:2017}, and
it has enriched computer science with powerful methods from order
theory, algebra, and topology. Domains gained prominence when they
were introduced as a mathematical model for lambda calculus
by~\textcite{Scott:Outline:1970}. We present a brief reminder of the
concepts and notations that will be needed later. The interested
reader may refer
to~\parencite{AbramskyJung94-DT,Gierz-ContinuousLattices-2003,Goubault-Larrecq:Non_Hausdorff_topology:2013} for
more on domains in general, and refer
to~\parencite{Escardo96-tcs,Edalat:Domains_Physics:1997} for the
interval domain, in particular. A succinct and informative survey of
the theory may be found
in~\parencite{Keimel:Domain_Ramifications_Interactions:2017}.

For any subset $X$ of a \acf{POSET} $(D, \sqsubseteq)$, we write
$\bigsqcup X$ and $\bigsqcap X$ to denote the least upper bound, and
the greatest lower bound, of $X$, respectively, whenever they
exist. We also define the lower and upper sets of $X$ by
$\lowerSet{X} \defeq \setbarNormal{y \in D}{\exists x \in X: y
  \sqsubseteq x}$ and
$\upperSet{X} \defeq \setbarNormal{y \in D}{\exists x \in X: x
  \sqsubseteq y}$, respectively. When $X$ is a singleton $\set{x}$, we
may simply write $\lowerSet{x}$ and $\upperSet{x}$, respectively.

In our discussion, $x \sqsubseteq y$ may be interpreted as `$y$
contains more information than $x$'.  A subset $X \subseteq D$ is said
to be \emph{directed} if $X\neq \emptyset$ and every two elements of
$X$ have an upper bound in $X$, {\ie},
$\forall x, y \in X: \exists z \in X: x \sqsubseteq z \wedge y
\sqsubseteq z$. A \ac{POSET} $(D, \sqsubseteq)$ is said to be
\emph{directed complete} if every directed subset $X \subseteq D$ has
a least upper bound in $D$. The \ac{POSET} $(D, \sqsubseteq)$ is said
to be \emph{pointed} if it has a bottom element, which we usually
denote by $\bot_D$, or if $D$ is clear from the context, simply by
$\bot$. In the current article, we assume that every \ac{DCPO} is
pointed.

Directed completeness is a basic requirement for almost all the
\acp{POSET} in our discussion. A notable exception is the \ac{POSET}
of partitions, which is essential in temporal discretization. In what
follows, we assume that $a > 0$ is a rational number. Although most of
the abstract theory may be presented by just assuming $a$ to be real,
for implementation and computable analysis of the algorithms (in
\ac{TTE}~\parencite{Weihrauch2000:book}) it is more convenient to take
$a \in \Q$.

\begin{definition}[Partitions]
  A partition of $[ 0, a]$ is a finite sequence $( q_0, \ldots, q_k)$
  of real numbers such that $k \geq 1$ and $0 = q_0 < \cdots <
  q_k=a$. Furthermore:

  \begin{enumerate}[label=(\roman*)]
  \item We let ${\mathcal{P}}$ denote the set of all partitions of an
    interval of interest, {\eg}, $[ 0, a]$.

  \item We let ${\mathcal{P}}_{\Q}$ denote the set of rational partitions,
    {\ie}, those that satisfy
    $\forall i \in \set{0, \ldots, k}: q_i \in \Q$.
    
  \item The norm $|Q|$ of a partition $Q = ( q_0, \ldots, q_k)$ is
    given by
    $|Q| \defeq \max \setbarNormal{q_i - q_{i-1}}{1 \leq i \leq k} $.

  \item We define
    ${\mathcal{P}}_1 \defeq \setbarNormal{ Q \in {\mathcal{P}}}{|Q| \leq 1}$ and
    ${\mathcal{P}}_{1,\Q} \defeq \setbarNormal{ Q \in {\mathcal{P}}_{\Q}}{|Q|
      \leq 1} = {\mathcal{P}}_1 \cap {\mathcal{P}}_{\Q}$.

  \item The minimal width of a partition $Q = ( q_0, \ldots, q_k)$ is
    given by
    $m(Q) \defeq \min \setbarNormal{q_i - q_{i-1}}{1 \leq i \leq k} $.
  \item We let $r_Q \defeq \frac{|Q|}{m(Q)}$. A partition
    $Q = ( q_0, \ldots, q_k)$ is said to be equidistant if and only if $r_Q = 1$,
    {\ie},
    $\forall i \in \set{1, \ldots, k}: q_i - q_{i-1} = \frac{a}{k}$.
  \item A partition $Q = ( q_0, \ldots, q_k)$ is said to refine
    another partition $P = ( p_0, \ldots, p_{\ell})$---denoted by
    $P \sqsubseteq Q$--- if and only if
    $\set{p_0, \ldots, p_{\ell}} \subseteq \set{q_0, \ldots, q_k}$,
    with $p_0 = q_0 = 0$ and $p_{\ell} = q_k = a$.
  \end{enumerate}
\end{definition}

As each partition must have only finitely many points, none of the
\acp{POSET} $({\mathcal{P}}, \sqsubseteq)$,
$({\mathcal{P}}_1, \sqsubseteq)$, $({\mathcal{P}}_{\Q}, \sqsubseteq)$,
or $({\mathcal{P}}_{1,\Q}, \sqsubseteq)$ may be directed
complete. They are, however, all pointed, and have lattice or
semilattice structures:

\begin{proposition}
  For any interval $[0,a]$ with $a \in \Q$, the sets ${\mathcal{P}}$ and
  ${\mathcal{P}}_{\Q}$ form lattices under the refinement order
  $\sqsubseteq$, while ${\mathcal{P}}_1$ and ${\mathcal{P}}_{1,\Q}$ form join
  semilattices.
\end{proposition}

\begin{proof}
  The fact that they form \acp{POSET} is immediate. Next, assume that
  $P \equiv \set{p_0, \ldots, p_{\ell}} \in {\mathcal{P}}$ and
  $Q \equiv \set{q_0, \ldots, q_k} \in {\mathcal{P}}$. Then, the partitions
  $P \sqcup Q$ and $P \sqcap Q$ are formed from the partition points
  $\set{p_0, \ldots, p_{\ell}} \cup \set{q_0, \ldots, q_k}$ and
  $\set{p_0, \ldots, p_{\ell}} \cap \set{q_0, \ldots, q_k}$,
  respectively. 
\end{proof}

We let $\kCPPO{\R^n\lift}$ denote the \ac{POSET} with the carrier set
$\setbarNormal{C \subseteq \R^n}{C \text{ is non-empty and compact}}
\cup \set{\R^n}$, ordered by reverse inclusion, {\ie},
$\forall X, Y \in \kCPPO{\R^n\lift}: X \sqsubseteq Y \iff Y \subseteq
X$. By further requiring the subsets to be convex, we obtain the
sub-poset $\cCPPO{\R^n\lift}$. Finally, we let $\intvaldom[\R^n\lift]$
denote the \ac{POSET} of hyper-rectangles of $\R^n$---{\ie}, subsets of the
form $\prod_{i=1}^n [a_i,b_i]$---ordered under reverse inclusion, with
$\R^n$ added as the bottom element. The three \acp{POSET} are \acp{DCPO}
and $\lub X = \bigcap X$, for any directed subset $X$. 

At the heart of domain theory lies the concept of \emph{way-below}
relation. Assume that $(D, \sqsubseteq)$ is a \ac{DCPO} and let $x,y
\in D$. The element $x$ is said to be \emph{way-below} $y$---written
as $x \ll y$---if for every directed subset $X$ of $D$, if $y
\sqsubseteq \lub X$, then there exists an element $d \in X$ such that
$x \sqsubseteq d$. Intuitively, the relation $x \ll y$ may be phrased
as `$x$ is a finitary approximation of $y$', or even as `$x$ is a lot
simpler than $y$'~~\cite[Section~2.2.1]{AbramskyJung94-DT}. An element
$x \in D$ is said to be \emph{finite} if $x \ll x$. The way-below
relation is, in general, stronger than the information order
$\sqsubseteq$. For instance, over $\kCPPO{\R^n\lift}$ (hence, also
over $\intvaldom[\R^n\lift]$ and $\cCPPO{\R^n\lift}$) we have the
following characterization:
\begin{equation}
  \label{eq:way_below_compact_sets}
  \forall K_1, K_2 \in \kCPPO{\R^n\lift}: K_1 \ll K_2 \iff K_2 \subseteq \interiorOf{K_1},
\end{equation}
in which $\interiorOf{K_1}$ denotes the interior of $K_1$. In
particular, $\kCPPO{\R^n\lift}$ (hence, also $\intvaldom[\R^n\lift]$
and $\cCPPO{\R^n\lift}$) has no finite element except $\bot$, which
is always finite.

For every element $x$ of a \ac{DCPO} $(D, \sqsubseteq)$, let
$\wayBelows{x} \defeq \setbarNormal{a \in D}{a \ll x}$. In
domain-theoretic terms, the elements of $\wayBelows{x}$ are the true
approximants of $x$. In fact, the way-below relation is also known as
the order of approximation~\cite[Section~2.2.1]{AbramskyJung94-DT}. A
subset $B$ of a \ac{DCPO} $(D, \sqsubseteq)$ is said to be a
\emph{basis} for $D$, if for every element $x \in D$, the set
$B_x \defeq \wayBelows{x} \cap B$ is a directed subset with supremum
$x$, {\ie}, $x = \lub B_x$. A \ac{DCPO} $(D, \sqsubseteq)$ is said to
be:
  \begin{enumerate}[label=(\roman*)]
  \item continuous if it has a basis;
  \item $\omega$-continuous if it has a countable basis;
  \item algebraic if it has a basis consisting entirely of finite
    elements;
  \item $\omega$-algebraic if it has a countable basis consisting
    entirely of finite elements.
  \end{enumerate}
  In this article, we call $(D, \sqsubseteq)$ a domain if it
    is a continuous \ac{DCPO}.\footnote{Recall that, in this article,
  we assume every \ac{DCPO} to be pointed.}

  Assume that $\Sigma$ is a finite alphabet. Let $\Sigma^*$ denote the
  set of finite strings over $\Sigma$, and let $\Sigma^\omega$ denote
  the set of countably infinite sequences over $\Sigma$. A typical
  example of an $\omega$-algebraic domain is given by the set
  $\Sigma^\infty \defeq \Sigma^* \cup \Sigma^\omega$, ordered under
  prefix relation:
  \begin{equation*}
    \forall s,t \in \Sigma^\infty: s \sqsubseteq t \iff \left( s \in
    \Sigma^\omega \wedge s = t \right) \vee \left(s \in \Sigma^* \wedge \exists z \in
    \Sigma^\infty: sz = t \right),
  \end{equation*}
  in which $sz$ denotes the concatenation of $s$ and $z$. The set
  $\Sigma^*$ is a basis of finite elements for
  $\Sigma^\infty$. Although algebraic domains have been used in real
  number
  computation~\parencite{Gianantonio:Real_Domain:IC:1996,Farjudian:Shrad:2007},
  non-algebraic domains have proven more suitable for computation over
  continuous spaces, {\eg}, dynamical
  systems~\parencite{Edalat95:DT-fractals}, exact real number
  computation~\parencite{Escardo96-tcs,Edalat:Domains_Physics:1997},
  differential equation
  solving~\parencite{Edalat_Pattinson2007-LMS_Picard}, reachability
  analysis of hybrid
  systems~\parencite{Edalat_Pattinson:Hybrid:2007,Moggi_Farjudian_Duracz_Taha:Reachability_Hybrid:2018},
  and robustness analysis of neural
  networks~\parencite{Zhou_Shaikh_Li_Farjudian:Robust_NN:MSCS:2023},
  to name a few. Hence, in this article, we will also be mainly
  working in the framework of non-algebraic $\omega$-continuous
  domains.
  
  The three \acp{DCPO} $\intvaldom[\R^n\lift]$, $\kCPPO{\R^n\lift}$,
  and $\cCPPO{\R^n\lift}$, are all $\omega$-continuous domains,
  because:

\begin{itemize}
\item
  $B_{\intvaldom[\R^n\lift]} \defeq \set{\R^n} \cup \setbarNormal{C
    \in \intvaldom[\R^n\lift]}{C \text{ is a hyper-rectangle with
      rational coordinates}}$ is a basis for $\intvaldom[\R^n\lift]$;

\item $B_{\kCPPO{\R^n\lift}} \defeq \set{\R^n} \cup \setbarNormal{C \in
  \kCPPO{\R^n\lift}}{C \text{ is a finite union of hyper-rectangles with
    rational coordinates}}$ is a basis for $\kCPPO{\R^n\lift}$;

\item and
  $B_{\cCPPO{\R^n\lift}} \defeq \set{\R^n} \cup \setbarNormal{C \in
    \cCPPO{\R^n\lift}}{C \text{ is a convex polytope with rational
      coordinates}}$ is a basis for $\cCPPO{\R^n\lift}$.

\end{itemize}
They are all non-algebraic because their only finite element is the
bottom element $\bot$.

The following \emph{interpolation property} is one of the most
important features of the way-below relation over domains:

\begin{lemma}[{\cite[Lemma~2.2.15]{AbramskyJung94-DT}}]
  \label{lemma:interpolation_property_cont_domains}
  Assume that $(D, \sqsubseteq)$ is a continuous domain. Let $z \in D$
  and let $X \subseteq D$ be a finite set satisfying
  $\forall x \in X: x \ll z$. Then, there exists an element $y \in D$
  interpolating between $X$ and $z$, {\ie}, satisfying
  $\forall x \in X: x \ll y \ll z$, which we write as $X \ll y \ll
  z$. Furthermore, if $B$ is a basis for $D$, then $y$ can be chosen
  from $B$.
\end{lemma}

\begin{remark}
  In what follows, we use `enclosure' as a generic term referring to
  intervals, hyper-rectangles, polytopes, etc., or functions that take
  such set values.
\end{remark}

\begin{remark}
  \Acp{DCPO} which are not continuous seldom appear naturally in
  applications, and examples of such posets are usually manufactured
  for providing insight (see,
  {\eg},~\cite[Section~2.2.3]{AbramskyJung94-DT}). One such instance,
  however, appears naturally in our discussion. The \ac{POSET}
  ${\mathcal{D}}_{\Q}$ of left semi-continuous enclosures, which we
  will study in Section~\ref{subsec:left_semi_cont_enclosures}, is a
  non-continuous \ac{DCPO}.
\end{remark}

Apart from order-theoretic structure, domains also have a topological
structure. Assume that $(D, \sqsubseteq)$ is a \ac{POSET}. A subset
$O \subseteq D$ is said to be \emph{Scott open} if it has the
following properties:
  \begin{enumerate}[label=(\arabic*)]
  \item It is an upper set, {\ie}, $\forall x \in O, \forall y \in D:
    x \sqsubseteq y \implies y \in O$.
  \item For every directed set $X \subseteq D$ for which $\lub X$
    exists, if $\lub X \in O$ then $X \cap O \neq \emptyset$.
  \end{enumerate}

  The collection of all Scott open subsets of a \ac{POSET} forms a
  $T_0$ topology, refered to as the Scott topology. A function
  $f: (D_1, \sqsubseteq_1) \to (D_2, \sqsubseteq_2)$ is said to be
  Scott continuous if it is continuous with respect to the Scott
  topologies on $D_1$ and $D_2$. Scott continuity can be stated purely
  in order-theoretic terms, {\ie}, a map
  $f: (D_1, \sqsubseteq_1) \to (D_2, \sqsubseteq_2)$ between two
  posets is Scott continuous if and only if it is monotonic and
  preserves the suprema of directed sets, {\ie}, for every directed
  set $X \subseteq D_1$ for which $\lub X$ exists, we have
  $f(\lub X) = \lub
  f(X)$~\parencite[Proposition~4.3.5]{Goubault-Larrecq:Non_Hausdorff_topology:2013}.
  
  For every element $x$ of a \ac{DCPO} $(D, \sqsubseteq)$, let
  $\wayAboves{x} \defeq \setbarNormal{a \in D}{x \ll a}$.
  \begin{proposition}[{\parencite[Proposition~2.3.6]{AbramskyJung94-DT}}]
    \label{prop:Scott_wayaboves}
    Let $D$ be a domain with a basis $B$. Then, for each $x \in D$,
    the set $\wayAboves{x}$ is open, and the collection
    ${\mathcal{O}} \defeq \setbarTall{\wayAboves{x}}{x \in B}$ forms a
    base for the Scott topology.
  \end{proposition}

  The maximal elements of $\intvaldom[\R^n\lift]$,
  $\kCPPO{\R^n\lift}$, and $\cCPPO{\R^n\lift}$ are singletons, and the
  sets of maximal elements may be identified with $\R^n$. As a
  corollary of Proposition~\ref{prop:Scott_wayaboves}, if
  ${\mathcal{O}}_S$ is the Scott topology on $\intvaldom[\R^n\lift]$,
  $\kCPPO{\R^n\lift}$, or $\cCPPO{\R^n\lift}$, then the restriction of
  ${\mathcal{O}}_S$ over $\R^n$ is the Euclidean topology. Thus, the
  sets of maximal elements are indeed homeomorphic to $\R^n$. For
  simplicity, we write $x$ to denote a maximal element $\set{x}$. For
  any $K \geq 0$, by restricting to $[-K,K]^n$, we obtain the
  $\omega$-continuous domains $\intvaldom[{[-K,K]^n}]$,
  $\kCPPO{[-K,K]^n}$, and $\cCPPO{[-K,K]^n}$, respectively.

  In general, if $D$ and $E$ are two domains, the space of
  Scott-continuous functions from $D$ to $E$, under pointwise
  ordering, may not be a domain. The study of Cartesian closed
  categories of domains has a rich literature, and the interested
  reader may refer to,
  {\eg},~\parencite{AbramskyJung94-DT,Goubault-Larrecq:Non_Hausdorff_topology:2013},
  and the references therein for more details. We say that a
  \ac{POSET} $(D,\sqsubseteq)$ is \emph{bounded-complete} if each
  bounded pair $x,y \in D$ has a
  supremum. By~\cite[Corollary~4.1.6]{AbramskyJung94-DT}, the category
  of bounded-complete domains is Cartesian closed. The domains
  $\intvaldom[{[-K,K]^n}]$, $\kCPPO{[-K,K]^n}$, and $\cCPPO{[-K,K]^n}$
  are all bounded-complete. In fact, all the domains that we use in
  the framework developed in this article are bounded-complete.
  
  Assume that $(X, \Omega(X))$ is a topological space, and let
  $(D, \sqsubseteq_D)$ be a bounded-complete continuous domain. We let
  $(D, \Scott{D})$ denote the topological space with carrier set $D$
  under the Scott topology $\Scott{D}$. The space $[X \to D]$ of
  functions $f: X \to D$ which are $(\Omega(X), \Scott{D})$ continuous
  can be ordered pointwise by defining:
\begin{equation*}
  \forall f, g \in [X \to D]: \quad f \sqsubseteq g \iff \forall x \in X:
  f(x) \sqsubseteq_D g(x).
\end{equation*}
It is straightfoward to verify that the \ac{POSET}
$([X \to D], \sqsubseteq)$ is directed-complete and
$\forall x \in X: (\lub_{i \in I} f_i)(x) = \lub
\setbarNormal{f_i(x)}{i \in I}$, for any directed subset
$\setbarNormal{f_i}{i \in I}$ of $[X \to D]$. A central question in
our discussion is whether this \ac{POSET} is continuous or not.

Consider the \ac{POSET} $(\Omega(X), \subseteq)$ of open subsets of
$X$ ordered under subset relation. For any topological space $X$, this
\ac{POSET} is a complete lattice, with $\emptyset$ as the bottom
element, and $X$ as the top element. Furthermore, we have:
\begin{equation*}
  \forall A \subseteq \Omega(X): \quad \lub A = \bigcup A \text{ and
  } \glb A = \interiorOf{(\bigcap A)}.
\end{equation*}
A topological space $(X, \Omega(X))$ is said to be \emph{core-compact}
if the lattice $(\Omega(X), \subseteq)$ is continuous. It is
well-known that:
\begin{theorem}
  \label{thm:core_compact}
  For any topological space $(X, \Omega(X))$ and non-singleton
  bounded-complete continuous domain $(D, \sqsubseteq_D)$, the
  function space $([X \to D], \sqsubseteq)$ is a bounded-complete
  continuous domain $\iff (X, \Omega(X))$ is core-compact.
\end{theorem}

\begin{proof}
  For the ($\Leftarrow$) direction,
  see~\parencite[Proposition~2]{Erker_et_al:way_below:1998}. A proof
  of the ($\Rightarrow$) direction can also be found
  on~\parencite[pages 62 and 63]{Erker_et_al:way_below:1998}.
\end{proof}

\begin{notation}[$X \Rightarrow D$]
  \label{notation:X_Rightarrow_D}
  Whenever $(X, \Omega(X))$ is a core-compact topological space, and
  $(D, \sqsubseteq_D)$ is a bounded-complete continuous domain, we use
  the notation $X \Rightarrow D$ to denote the bounded-complete
  continuous domain $([X \to D], \sqsubseteq)$.
  \end{notation}

  Every locally compact space is
  core-compact~\parencite[Theorem~5.2.9]{Goubault-Larrecq:Non_Hausdorff_topology:2013}. It
  is well-known that Euclidean spaces are locally
  compact. Furthermore, every domain is locally compact in its Scott
  topology~\parencite[Corollary~5.1.36]{Goubault-Larrecq:Non_Hausdorff_topology:2013}. As
  such, Euclidean spaces and continuous domains form two important
  cases of core-compact spaces that are relevant to our discussion.

  \begin{definition}[$D_m^{(0)}(X)$]
    \label{def:D_0_m}
    For any set $X \subseteq \R^n$ which is locally compact under the
    restriction of the Euclidean topology , we let $D_m^{(0)}(X)$
    denote the function space $X \Rightarrow \intvaldom[\R^m\lift]$,
    with $X$ under the restriction of the Euclidean topology and
    $\intvaldom[\R^m\lift]$ under the Scott topology.
\end{definition}

If $(X, \Omega(X))$ is any topological space, then $f: X \to \R$ is
said to be:

\begin{itemize}
\item upper semi-continuous at $x_0 \in X$ $\iff$ for every
  $y > f(x_0)$, there exists a neighborhood $U \in \Omega(X)$ of $x_0$
  such that $\forall x \in U: f(x) < y$.

\item lower semi-continuous at $x_0 \in X$ $\iff$ for every
  $y < f(x_0)$, there exists a neighborhood $U \in \Omega(X)$ of $x_0$
  such that $\forall x \in U: f(x) > y$.

\item upper (respectively, lower) semi-continuous $\iff$ it is upper
  (respectively, lower) semi-continuous at every $x_0 \in X$.

\end{itemize}
The following is a well-known fact (see,
{\eg},~\parencite{Edalat_Lieutier:Domain_Calculus_One_Var:MSCS:2004})
and follows from the definition of semi-continuity:
\begin{proposition}
  \label{prop:D_0_n:Upper_Lower_SC}
  A function
  $f \equiv (f_1, \ldots, f_n ): [0,a] \to \intvaldom[\R^n\lift]$ is
  in $D_n^{(0)}([0,a])$ if and only if for every $j \in \set{1, \ldots, n}$,
  \uep[f_j] is upper semi-continuous and \lep[f_j] is lower
  semi-continuous.
\end{proposition}

As $\intvaldom[\R^n\lift]$ is a sub-poset of $\kCPPO{\R^n\lift}$, for
any compact set $K \in \kCPPO{\R^n\lift}$, we may define the
hyper-rectangular closure as
$\hrEnclos{K} \defeq \lub{\setbarNormal{R \in \intvaldom[\R^n\lift]}{K
    \subseteq R}}$, {\ie}, the smallest axes-aligned hyper-rectangle
containing $K$.

\begin{proposition}[{\parencite[Corollary~2.17]{Zhou_Shaikh_Li_Farjudian:Robust_NN:MSCS:2023}}]
      \label{prop:box_map_scott_cont}
      The map
    $\hrEnclos{(\cdot)} : \kCPPO{\R^n\lift} \to \intvaldom[\R^n\lift]$
    is Scott-continuous.
\end{proposition}

\begin{definition}[Extension, Canonical Interval Extension {$\intvalfun[f]$}, Approximation]{\ }
  \label{def:ext_canonical_ext}
  
  \begin{enumerate}[label=(\roman*)]

  \item \label{item:extension} A map
    $u: \kCPPO{\R^n\lift} \to \kCPPO{\R^m\lift}$ is said to be an
    extension of $f: \R^n \to \R^m$ iff
    $\forall x \in \R^n: u(\set{x}) = \set{f(x)}$.

  \item A map $u: \intvaldom[\R^n\lift] \to \intvaldom[\R^m\lift]$ is
    said to be an interval extension of
    $f: \R^n \to \kCPPO{\R^m\lift}$ iff
    $\forall x \in \R^n: u(\set{x}) = \hrEnclos{f(x)}$.

  \item \label{item:canonical_interval_extension} For any
    $f : \R^n \to \kCPPO{\R^m\lift}$, we define the canonical interval
    extension
    ${\intvalfun[f]} : \intvaldom[\R^n\lift] \to
    \intvaldom[\R^m\lift]$ by:
    \begin{equation}
      \label{eq:canonical_ext}
      \forall \alpha \in \intvaldom[\R^n\lift]: \quad
      {\intvalfun[f]}(\alpha) 
      \defeq \bigsqcap_{x \in \alpha} \hrEnclos{f(x)}.
    \end{equation}

  \item A map $u: \intvaldom[\R^n\lift] \to \intvaldom[\R^m\lift]$ is
    said to be an interval approximation of
    $f: \R^n \to \kCPPO{\R^m\lift}$ if $u \sqsubseteq \intvalfun[f]$.   
  \end{enumerate}
\end{definition}

\begin{proposition}
\label{prop:canonical_interval_extension}
For every Euclidean-Scott-continuous $f : \R^n \to \kCPPO{\R^m\lift}$,
the canonical interval extension ${\intvalfun[f]}$ defined
in~\eqref{eq:canonical_ext} is the maximal extension of $f$ among all
the interval extensions in the domain
$\intvaldom[\R^n\lift] \Rightarrow \intvaldom[\R^m\lift]$. In
particular, ${\intvalfun[f]}$ is Scott-continuous.
\end{proposition}

\begin{proof}
  Given a map $f : \R^n \to \kCPPO{\R^m\lift}$, we define
  $\hrEnclos{f}: \R^n \to \intvaldom[\R^m\lift]$ by
  $\hrEnclos{f} \defeq \hrEnclos{(\cdot)} \circ f$, {\ie},
  $\forall x \in \R^n: \hrEnclos{f}(x) = \hrEnclos{f(x)}$. If $f$ is
  Euclidean-Scott-continuous, then, by
  Proposition~\ref{prop:box_map_scott_cont}, so is $\hrEnclos{f}$. It
  is straightforward to verify that a map
  $u: \intvaldom[\R^n\lift] \Rightarrow \intvaldom[\R^m\lift]$ is an
  interval approximation of $f$ if and only if it is an interval
  approximation of~$\hrEnclos{f}$.

  Thus, it suffices to prove the proposition for the special case of
  $f: \R^n \to \intvaldom[\R^m\lift]$. This has already been
  proven in~\parencite[Lemma~3.4]{Edalat_Escardo:Integ_realPCF:2000}
  for the case of $n = m = 1$. But, the proof given
  in~\parencite{Edalat_Escardo:Integ_realPCF:2000} is independent of
  the values of $n$ and $m$, because the crucial property that is
  needed is that $\intvaldom[\R^m\lift]$ is a continuous
  $\glb$-semilattice, for any $m \in \N$.
\end{proof}

As the restriction of the Scott topology of $\kCPPO{\R^m\lift}$ over
$\R^m$ is the Euclidean topology, we may consider any continuous map
$f: \R^n \to \R^m$ also as a function of type
$\R^n \to \kCPPO{\R^m\lift}$, and construct its canonical interval extension
accordingly, which will be Scott-continuous.

In lambda calculus, a fixpoint combinator is used for recursive
definitions. One common fixpoint term is the so-called $Y$ combinator
$Y \defeq \lambda f. (\lambda x. f( xx)) (\lambda x. f( xx))$. As
previously mentioned, domains were introduced to construct
mathematical models of lambda
calculus~\parencite{Scott:Outline:1970}. The denotation of a fixpoint
combinator may be provided by the (domain-theoretic) fixpoint
operator, as specified in the following theorem:

  \begin{theorem}[fixpoint operator: $\fix$] 
    Assume that $D$ is a \ac{DCPO} with bottom element $\bot$. Then:

    \begin{enumerate}[label=(\roman*)]
    \item Every Scott-continuous $f: D \to D$ has a least fixpoint
      given by $\bigsqcup_{n \in \N} f^n(\bot)$.

    \item The fixpoint operator $\fix: [D \to D] \to D$ defined by
      $\fix f \defeq \bigsqcup_{n \in \N} f^n(\bot)$ is Scott
      continuous.
    \end{enumerate}
  \end{theorem}

  \begin{proof}
    See~\cite[Theorem~2.1.19]{AbramskyJung94-DT}.
  \end{proof}

  Domains provide a natural setting for the concept of
  approximation. In particular, objects which are not
  finitely-representable may be constructed via their finite
  approximations. A common approach in this context is through the use
  of the fixpoint operator. One of the main contributions of the
  current paper is the formulation of Euler and Runge-Kutta operators
  using the fixpoint operator (Theorem~\ref{thm:E_2_MFPS_Ideal_Compl}
  and Definition~\ref{def:E_R_Ideal_Compl}).

\subsection{Domain-Theoretic Derivative}
\label{subsec:domain_theoretic_derivative}

We recall the concept of a domain-theoretic derivative for a function
$f: U \to \R$ defined on an open set $U \subseteq \R^n$. Assume that
$(X, \Omega(X))$ is a topological space, and $(D, \sqsubseteq)$ is a
\ac{DCPO}, with bottom element $\bot$. Then for any open set
$O \in \Omega(X)$, and any element $b \in D$, we define the single-step
function $b \chi_O: X \to D$ as follows:
\begin{equation}
  \label{eq:single_step_fun}
  b \chi_O (x) \defeq
  \left\{
    \begin{array}{ll}
      b, &  \text{if } x \in O,\\
      \bot, & \text{if } x \in X \setminus O,
    \end{array}
  \right.
\end{equation}  

\begin{definition}[$L$-derivative]
  \label{def:L_derivative}

  Assume that $O \subseteq U \subseteq \R^n$, both $O$ and $U$ are
  open, and $b \in \cCPPO{\R^n\lift}$:

  \begin{enumerate}[label=(\roman*)]

  \item The single-step tie $\delta(O, b)$ is the set of all
    functions $f : U \to \R$ which
    satisfy:
    \begin{equation*}
      \forall x, y \in O:\quad b(x-y) \sqsubseteq f(x) - f(y).
    \end{equation*}
    The set $b(x-y)$ is obtained by taking the inner product of
    every element of $b$ with the vector $x-y$, and $\sqsubseteq$ is
    the reverse inclusion order on $\mathbf{C}\R^1\lift$.

  \item The $L$-derivative of any function $f: U \to \R$ is defined
    as:
    \begin{equation*}
     L(f) \defeq \bigsqcup \setbarTall{b \chi_O}{f \in \delta(O, b)}. 
    \end{equation*}

\end{enumerate}
  
\end{definition}

The $L$-derivative is a Scott-continuous function. When $f$ is
classically differentiable at $x \in U$, the $L$-derivative and the
classical derivative
coincide~\parencite{Edalat:2008:Continuous_Derivative}. Many of the
fundamental properties of the classical derivative can be generalized
to the domain theoretic one, {\eg}, additivity and the chain
rule~\parencite{Edalat_Pattinson:2005:Inverse_Implicit}. A
generalization of the mean value theorem~\cite[Theorem
5.4]{EdalatLieutierPattinson:2013-MultiVar-Journal} is essential in
the proof of Lemma~\ref{lemma:Taylor_Lipschitz_OneVar}, which
underlies the soundness of the Euler and Runge-Kutta operators. This
generalization follows from the corresponding result for the
Clarke-gradient~\cite[Theorem
2.3.7]{Clarke:Opt_Non_Smooth_Analysis-Book:1990}, and the fact that
the Clarke-gradient coincides with the domain theoretic derivative,
which was first proven for finite dimensional Banach spaces
by~\citeauthor{Edalat:2008:Continuous_Derivative}~\parencite{Edalat:2008:Continuous_Derivative},
and later generalized to infinite dimensional Banach spaces by
\citeauthor{Hertling:Clarke_Edalat:2017}~\parencite{Hertling:Clarke_Edalat:2017}.

In this paper, instead of working with the general convex sets in
$\cCPPO{\R^n\lift}$, we work with the simpler hyper-rectangular ones
in $\intvaldom[\R^n\lift]$:
\begin{definition}[$\overline{L}$-derivative]
  \label{def:Lhat_derivative}
  Let $U \subseteq \R^n$ be an open set. In the definition of
  $L$-derivative (Definition~\ref{def:L_derivative}), if
  $\cCPPO{\R^n\lift}$ is replaced with $\intvaldom[\R^n\lift]$, then
  we obtain the concept of $\overline{L}$-derivative for functions of type
  $U \to \R$.
\end{definition}
Clearly, the $\overline{L}$-derivative is, in general, coarser than the
$L$-derivative:

\begin{example}
  Let $D_n$ denote the closed unit disc in $\R^n$, and define
  $f: \R^n \to \R$ by $f(x) = \norm{x}_2$, in which $\norm{\cdot}_2$
  is the Euclidean norm. Then, $L(f)(0) = D_n$, while
  $\overline{L}(f)(0) = [-1,1]^n$.
\end{example}

\begin{definition}
  For every open set $U \subseteq \R^n$ and vector valued function
  $f \equiv ( f_1, \ldots, f_m): U \to \R^m$, we define:
  \begin{equation*}
    \overline{L}(f) \defeq  \left( \overline{L}(f_1),  \ldots, \overline{L}(f_m)
    \right)^\intercal,
  \end{equation*}
  in which, $(\cdot)^\intercal$ denotes the transpose of a matrix. In
  other words, for each $1 \leq i \leq m$ and $x \in U$, let
  $\overline{L}(f_i)(x) \equiv ( \alpha_{i,1}, \ldots, \alpha_{i,n}) \in \intvaldom[\R^n\lift]$. Then, for
  each $x \in U$, $\overline{L}(f)(x)$ is the $m \times n$ interval matrix
  $[\alpha_{i,j}]_{1 \leq i \leq m, 1 \leq j \leq n}$.
\end{definition}

The solution of the \ac{IVP}~(\ref{eq:main_ivp}) is a function of type
$y \equiv (y_1, \dots, y_n): [0, a] \to [-K,K]^n$. For each component
$y_j$, with $1 \leq j \leq n$, if $M$ is a local Lipschitz constant
for $y_j$ around $x$, then
$[-M, M] \sqsubseteq L(y_j)(x) = \overline{L}(y_j)(x)$. We also have:
  \begin{equation*}
    \overline{L}(y) = (L(y_1), \dots, L(y_n))^\intercal : [0, a] \to \intvaldom[\R^n\lift].
  \end{equation*}

  \noindent
  In general, $\overline{L}(y)(x)$ contains the generalized Jacobian of $y$
  at $x$.

\begin{definition}[${\mathcal{V}}^1$]
  \label{def:V1}
  Consider the domain:
  \begin{equation*}
    {\mathcal{V}}^* \defeq \left( \intvaldom[{[-K,K]^n}] \Rightarrow
      \intvaldom[{[-M,M]^n}] \right) \times \left(
      \intvaldom[{[-K,K]^n}] \Rightarrow \intvaldom[{[-M_1,M_1]^{n^2}}]
    \right).  
  \end{equation*}
  We say that a pair $(u,u') \in {\mathcal{V}}^*$ is consistent if
  there exists an $h: [-K,K]^n \to [-M,M]^n$ satisfying
  $ u \sqsubseteq \intvalfun[h]$ and
  $u' \sqsubseteq \intvalfun[\overline{L}(h)]$. We define the domain
  ${\mathcal{V}}^1$ to be the sub-domain of ${\mathcal{V}}^*$
  consisting of consistent pairs $( u, u') \in {\mathcal{V}}^*$.
\end{definition}

We use consistent pairs to approximate functions and their
derivatives. Specifically, in a consistent pair $(u,u')$, the
component $u$ is meant to approximate a function while $u'$
approximates the derivative.\footnote{The presence of $K$, $M$, and
  $M_1$ in Definition~\ref{def:V1} reflects the fact that the pair
  $(u,u')$ is meant to approximate the vector field $f$
  of~\eqref{eq:main_ivp} and its derivative.} For more on consistency
and strong consistency, the reader may refer
to~\parencite{Edalat_Lieutier:Domain_Calculus_One_Var:MSCS:2004,EdalatLieutierPattinson:2013-MultiVar-Journal},
where it has also been shown that ${\mathcal{V}}^1$ is indeed a
domain.

The domain ${\mathcal{V}}^1$ can be equipped with an effective
structure~\parencite{EdalatLieutierPattinson:2013-MultiVar-Journal}. This
is crucial when dealing with imprecisely given input data as we need
an effective method for verifying the consistency of a given pair
$(u,u') \in {\mathcal{V}}^*$. The effective structure will also be
central in computable analysis, as will be seen when we study
computability of the Euler operator in
Section~\ref{subsec:computability}.

\begin{remark}
  Although the $\overline{L}$-derivative is coarser than the
  $L$-derivative, we will still obtain completeness and second-order
  convergence for our \ac{IVP} solver
  (Theorem~\ref{thm:Second_Order_Convergence_int_Lip}). Implementing
  hyper-rectangles is easier and more efficient compared with compact
  convex sets. Furthermore, the consistency relation is decidable over
  rational hyper-rectangles, but it remains open whether it is
  decidable over rational convex polytopes in
  $\cCPPO{\R^n\lift}$~\parencite{EdalatLieutierPattinson:2013-MultiVar-Journal}.
\end{remark}

\begin{definition}[${\mathcal{V}}_{\negmedspace f}^1$]
  \label{def:V1_f}
  We define the continuous lattice ${\mathcal{V}}_{\negmedspace f}^1$
  to be the sub-domain of ${\mathcal{V}}^1$ with the carrier set:
  \begin{equation*}
    {\mathcal{V}}_{\negmedspace f}^1 \defeq \setbarNormal{(u,u') \in  {\mathcal{V}}^1}{u \sqsubseteq \intvalfun[f]
  \text{ and }u' \sqsubseteq \intvalfun[\overline{L}(f)]},
  \end{equation*}
  in which $f: [-K,K]^n \to [-M,M]^n$ is a Lipschitz continuous vector
  field with $M_1$ as a Lipschitz constant.\footnote{In the current
    article, we will be mainly concerned with the continuous lattice
    ${\mathcal{V}}_{\negmedspace f}^1$ where $f$ is the vector field
    of the \ac{IVP}~\eqref{eq:main_ivp}.}
\end{definition}

\begin{notation}
  If $f \equiv (f_1, \dots, f_n): [a, b] \to \R^n$ is $k$-times
  differentiable and $0 \leq i \leq k$, we write
  $f^{(i)}: [a, b] \to \R^n$ for the $i$-th classical derivative of
  $f$. In particular, $f^{(0)} = f$, $f^{(1)} = f'$, and
  $f^{(2)}= f''$.
\end{notation}

The following is a generalization of the domain of scalar $C^1$
functions introduced
in~\parencite{Edalat_Lieutier:Domain_Calculus_One_Var:MSCS:2004}. The
constants $M, M_1, \ldots, M_p$ will be fixed depending on the
application:

\begin{definition}[$\hat{D}^{(p)}, \hat{D}^{(p)}_f$]
  \label{def:domain_of_Cp_funs}
  
  For every $p \in \N$, we define the domain:
  \begin{equation*}
  D^{(p)}_* \defeq \left( \intvaldom[{[-K,K]}] \Rightarrow
    \intvaldom[{[-M,M]}] \right) \times \left( \intvaldom[{[-K,K]}] \Rightarrow
    \intvaldom[{[-M_1,M_1]}] \right) \times \ldots \times \left( \intvaldom[{[-K,K]}] \Rightarrow
    \intvaldom[{[-M_p,M_p]}] \right).
  \end{equation*}
  We let
  $\hat{D}^{(0)} \defeq D^{(0)}_*$, and for every $p \geq 1$, we let
  $\hat{D}^{(p)}$ denote the sub-domain of consistent tuples of
  $D^{(p)}_*$, {\ie}:
  \begin{equation*}
    \hat{D}^{(p)} \defeq \setbarTall{(u_0, \ldots, u_p) \in D^{(p)}_*}{
      \exists g \in C^{p-1}([-K,K]): \forall i \in \set{0, \ldots, p-1}: u_i
      \sqsubseteq \intvalfun[g^{(i)}] \wedge u_p \sqsubseteq \intvalfun[\overline{L}(g^{(p-1)})]}.
  \end{equation*}
  For a fixed $f \in C^{p-1}([-K,K])$, we define:
  \begin{equation*}
    \hat{D}^{(0)}_f \defeq \setbarTall{u \in \hat{D}^{(0)}}{u
      \sqsubseteq \intvalfun[f]},
  \end{equation*}
  and for $p \geq 1$, we let $\hat{D}^{(p)}_f$ denote the sub-domain
  of $\hat{D}^{(p)}$ consisting of those tuples that approximate $f$
  and its derivatives, {\ie}:
  \begin{equation}
    \label{eq:hat_D_p_f}
    \hat{D}^{(p)}_f \defeq \setbarTall{(u_0, \ldots, u_p) \in
      \hat{D}^{(p)}}{\forall i \in \set{0, \ldots, p-1}: u_i
      \sqsubseteq \intvalfun[f^{(i)}] \wedge u_p \sqsubseteq \intvalfun[\overline{L}(f^{(p-1)})]}.
  \end{equation}

\end{definition}

Note that in~\eqref{eq:hat_D_p_f}, for $\hat{D}^{(p)}_f$ to be
non-empty, the following must hold:

\begin{enumerate}[label=(\roman*)]
\item $\forall x \in [-K,K]: f(x) \in [-M,M]$.
\item
  $\forall x \in [-K,K], \forall i \in \set{1, \ldots, p-1}:
  f^{(i)}(x) \in [-M_i,M_i]$.
\item The function $f^{(p-1)}$ must be Lipschitz continuous with $M_p$
  as a Lipschitz constant.
\end{enumerate}

\subsection{Interval Analysis}
\label{subsec:interval_analysis}

We will also use some concepts from interval analysis, especially for
convergence analysis in Section~\ref{subsec:convergence_analysis}.

  \begin{definition}[Width]{\ }
    \begin{enumerate}[label=(\roman*)]
    \item For any interval
      $\alpha = [{\lep[\alpha]}, {\uep[\alpha]}]$, the \emph{width} of
      $\alpha$ is defined by ${\widthOf(\alpha)} \defeq {\uep[{\alpha}]} - {\lep[{\alpha}]}.$

    \item For an $n$-dimensional hyper-rectangle ({\ie}, interval
      vector) $\vec{\alpha} = {\ntuple{\alpha}{1}{n}}$, the width is
      defined by
      ${\widthOf(\vec{\alpha})} \defeq \max_{1 \leq i \leq n}
      {\widthOf({\alpha}_{i})} $.

    \item For an interval matrix $A=[A_{ij}]_{m \times n}$, we define
      $\widthOf(A) \defeq \max \setbarTall{\widthOf(A_{ij})}{1 \leq
        i \leq m, 1 \leq j \leq n}$.
      
    \item If width is given over a set $B$, then it can be lifted to
      functions from any set $A$ to $B$ by:
      \begin{equation*}
        \forall f \colon A \to B \colon \quad {\widthOf(f)} \coloneqq \sup \{ {\widthOf(f(x))} \mid x \in A \}.        
      \end{equation*}
      Note that it is possible to have $\widthOf(f) = \infty$.
    \end{enumerate}
  \end{definition}

  \begin{definition}[$\norm{I}, \norm{A}_1, \norm{A}_\infty$]

  \label{def:interval_norms}

  For every interval $I = [a,b]$, we define
  $\norm{I} \defeq \max(\absn{a}, \absn{b})$. For an interval matrix $A = [A_{ij}]_{m \times n}$, the norms
    $\norm{A}_1$ and $\norm{A}_\infty$ are defined as follows:
    \begin{equation}
      \label{eq:intval_matrix_norms}
      \left\{
        \arrayoptions{0.3ex}{1.3}
      \begin{array}{lcl}
      \norm{A}_1 & \defeq & \max_{1 \leq j \leq n} \sum_{i=1}^{m} \norm{A_{ij}},\\        
      \norm{A}_\infty & \defeq & \max_{1 \leq i \leq m} \sum_{j=1}^{n} \norm{A_{ij}}.
      \end{array}
    \right.
    \end{equation}
  \end{definition}

  \begin{remark}

    Assume that $A = [A_{ij}]_{m \times n}$ is a matrix of real
    numbers, viewed as a linear operator $A : \R^n \to \R^m$. For each
    $p \in [1,\infty]$, if we use the $p$-norm for vectors on $\R^n$
    and $\R^m$, then the corresponding operator norm for $A$ is:
    \begin{equation*}
      \norm{A}_p = \sup_{x \in \R^n \setminus \set{0}} \frac{\norm{Ax}_p}{\norm{x}_p}.
    \end{equation*}
    It can be shown that for the specific cases of $p \in \set{1,
      \infty}$, we have:

    \begin{equation}
      \label{eq:real_matrix_norms}
      \left\{
        \arrayoptions{0.3ex}{1.3}
      \begin{array}{lcl}
      \norm{A}_1 & = & \max_{1 \leq j \leq n} \sum_{i=1}^{m} \absn{A_{ij}},\\        
      \norm{A}_\infty & = & \max_{1 \leq i \leq m} \sum_{j=1}^{n} \absn{A_{ij}}.
      \end{array}
    \right.
  \end{equation}
    As such, the interval matrix norms
    of~\eqref{eq:intval_matrix_norms} are generalizations of the real
    matrix norms of~\eqref{eq:real_matrix_norms}.    
  \end{remark}
  
  For error analysis of interval computations, a metric structure is
  required. \textcite[page~52]{Moore:2009:IIA} used the Hausdorff
  distance for interval analysis. The following extends their
  definition to tuples and functions:
  \begin{definition}[Interval Distance]{\ }
  \label{def:interval_distance}

  \begin{enumerate}[label=(\roman*)]
  \item \label{item:Moore_Intval_Dist}  For any pair of intervals $x = [\lep[x], \uep[x]]$ and
  $y = [\lep[y], \uep[y]]$, let $d(x,y) \defeq \max\left( \absn{\lep[x] - \lep[y]} ,  \absn{\uep[x]
        - \uep[y]}\right)$.
  
  \item \label{item:Intval_Dist_IRn} For $\alpha = ( \alpha_1, \ldots, \alpha_n)$ and $\beta = (
    \beta_1, \ldots, \beta_n)$ in $\intvaldom[\R^n]$, we let $d( \alpha, \beta) \defeq \max \setbarNormal{d( \alpha_i,
        \beta_i)}{ 1 \leq i \leq n}$.

  \item Let $X$ be an arbitrary set, and let
    $f, g: X \to \intvaldom[\R^n]$. Then, we define $d( f, g) \defeq \sup \setbarNormal{d( f(x), g(x))}{ x \in X}$.

  \end{enumerate}
\end{definition}
\begin{definition}[Symmetric Expansion]
  \label{def:symmetric_expansion}
  For any
  $\alpha = ( [ \lep[a_1], \uep[a_1]], \ldots, [\lep[a_n],
  \uep[a_n]]) \in \intvaldom[\R^n]$ and $r \geq 0$, we define the
  symmetric expansion of the interval vector $\alpha$ with the real
  constant $r$ as follows:

  \begin{equation*}
    \alpha \symExpand r \defeq ( [\lep[a_1] - r, \uep[a_1] + r], \ldots, [\lep[a_n] - r, \uep[a_n] + r]).
  \end{equation*}
\end{definition}
%
%
\noindent
Note that  $\alpha \symExpand r$ is the Minkowski sum of
  $\prod_{1 \leq i \leq n} [ \lep[a_i], \uep[a_i]]$ and $[-r,r]^n$, as
  subsets of $\R^n$.

We say that a function
$u: \intvaldom[\R^n\lift] \to \intvaldom[\R^m\lift]$, with
$n, m \in \N$, is interval
Lipschitz~\cite[Definition~6.1]{Moore:2009:IIA} with constant
$L \geq 0$ if and only if:

  \begin{equation*}
    \forall \alpha \in \intvaldom[\R^n]: \quad \widthOf( u(
    \alpha)) \leq L \widthOf( \alpha).
  \end{equation*}
\begin{remark}
  \label{rem:intval_Lipschitz_total}
  If $u$ is interval Lipschitz, then it is \emph{total}, {\ie}, it
  maps every maximal element of $\intvaldom[\R^n\lift]$ to a maximal
  element of $\intvaldom[\R^m\lift]$, because
  $\forall x \in \R^n: \widthOf(u(\set{x})) \leq L \widthOf(\set{x}) =
  0$.
\end{remark}

\section{A Domain for Function Spaces}

\label{sec:domain_for_fun_spaces}

The interval $[0,a]$ under the Euclidean topology is
core-compact. Thus, the function space $D_n^{(0)}([0,a])$ of
Definition~\ref{def:D_0_m} is a continuous domain. The domain
$D_n^{(0)}([0,a])$ is suitable for \ac{IVP} solving using Picard
method~\parencite{Edalat_Lieutier:Domain_Calculus_One_Var:MSCS:2004,Edalat_Pattinson2007-LMS_Picard}. For
methods such as Euler and Runge-Kutta, which proceed according to a
temporal discretization, the Euclidean topology is not
suitable. Instead, we must work with (a variant of) the so-called
\emph{upper limit} topology on $[0,a]$, which is not core-compact
(Proposition~\ref{prop:0_a_not_core_compact}). This necessitates the
construction of a continuous domain for dealing with function spaces
$[X \to D]$ when $X$ is not core-compact.

A function $f: X \to D$ is said to be a step function if it is the
supremum of a finite set of single-step functions, {\ie},
$f = \lub_{i \in I} b_i \chi_{O_i}$, in which $I$ is finite and
$b_i \chi_{O_i}$ is as defined in~\eqref{eq:single_step_fun}. Assume
that $B(X)$ is a base of the topology $\Omega(X)$, and $B(D)$ is a
basis---in the domain-theoretic sense---of the continuous domain
$D$. We let $\mathcal{B}$ denote the set of step functions defined
using elements of $B(X)$ and $B(D)$, {\ie}
\begin{equation}
  \label{eq:def_B}
  \mathcal{B} \defeq \setbarTall{f: X \to D}{f = \lub_{i \in I} b_i
    \chi_{O_i}, I \text{ is finite}, \forall i \in I: O_i \in B(X) \wedge
    b_i \in B(D) }.
\end{equation}
\noindent
In~\eqref{eq:def_B}, we have tacitly assumed that each $f$ is
well-defined. More precisely, the supremum
$\lub_{i \in I} b_i \chi_{O_i}$ exists if and only if
$\setbarTall{b_i \chi_{O_i}}{i \in I}$ satisfies the following
consistency condition:
\begin{equation*}
  \forall J \subseteq I: \quad \bigcap_{j \in J} O_j \neq \emptyset
  \implies \exists b_J \in B(D): \forall j \in J: b_j \sqsubseteq b_J.
\end{equation*}

When $X$ is core-compact, the set $\mathcal{B}$ provides a basis for
the continuous domain $[X \to
D]$~\parencite{Erker_et_al:way_below:1998}. If $X$ is not core-compact,
then by Theorem~\ref{thm:core_compact}, the function space $[X \to D]$
is not a continuous domain. In applications of domain theory, normally
a continuous domain is constructed first, and then one of its bases is
identified for further analysis. At times, however, it is useful to
take the opposite approach, {\ie}, start with a given structure as a
basis, and then construct a continuous domain over that basis. The
essential properties required of such a structure to enable the
construction of a continuous domain are quite minimal. This is
captured by the concept of an \emph{abstract basis}:

\begin{definition}[Abstract basis]
  \label{def:abstract_basis}
  A pair $(B, \prec)$ consisting of a set $B$ and a binary relation
  $\prec\ \subseteq B \times B$ is said to be an abstract basis if the
  relation $\prec$ is transitive and satisfies the following
  interpolation property:
  \begin{itemize}
  \item For every finite subset $A \subseteq B$ and element
    $x \in B: A \prec x \implies \exists y \in B: A \prec y \prec x$.
  \end{itemize}
  Here, by $A \prec x$ we mean $\forall a \in A: a \prec x$. 
\end{definition}

Continuous domains may be constructed over abstract bases by a
completion process which is similar to how the set $\R$ of real
numbers is constructed by completion of the set $\Q$ of rational
numbers.\footnote{To be more specific, the completion process is
  similar to the construction of real numbers via Dedekind cuts of
  rational
  numbers~\cite[Exercise~III-4.18]{Gierz-ContinuousLattices-2003}.}
The process is refered to as \emph{ideal completion}, which requires a
certain type of ideals, called \emph{rounded ideals}:

\begin{definition}[Rounded ideal]
  \label{def:rounded_ideal}
  A subset $I$ of an abstract basis $(B, \prec)$ is said to be a
  rounded ideal if it satisfies the following conditions:

  \begin{enumerate}[label=(\roman*)]
  \item \label{item:def_rounded_ideal_non_empty} $I$ is non-empty.
  \item $I$ is a lower set, {\ie}, $\forall y \in I, x \in B: x \prec
    y \implies x \in I$.

  \item \label{item:def_rounded_ideal_directed} $I$ is directed, {\ie}, for every finite set $A \subseteq I$,
    there exists an element $z \in I$ such that $A \prec z$.
  \end{enumerate}
\end{definition}

\begin{remark}
  Condition~\ref{item:def_rounded_ideal_non_empty} in
  Definition~\ref{def:rounded_ideal} is redundant as it follows from
  condition~\ref{item:def_rounded_ideal_directed} by taking
  $A = \emptyset$. Nonetheless, we keep non-emptiness in the statement
  as it makes it explicit and also makes the proof of some later
  statements ({\eg},
  Proposition~\ref{prop:f_star_rounded_ideal_X_to_D}) more intuitive.
\end{remark}

A continuous domain can be constructed over an abstract basis
$(B, \prec )$ by taking the set of rounded ideals of $B$ under the
subset relation. In other words, if we denote the set of the rounded
ideals of $(B, \prec)$ by $\RId{B, \prec}$, then $(\RId{B, \prec}, \subseteq)$ is a
continuous domain. For more on construction of domains using abstract
bases, the reader may refer to~\cite[Section~2.2.6]{AbramskyJung94-DT}
or~\cite[Section~III-4]{Gierz-ContinuousLattices-2003}.

\begin{remark}
  The symbol `$\prec$' is commonly used to denote the order over
  abstract bases. In what follows, we will use the symbol
  `$\precAbsBas$' instead to distinguish the order over the abstract
  bases of step functions from the general case of
  Definition~\ref{def:abstract_basis}. Furthermore, this helps us
  avoid confusion with the common `strictly less than' relation
  symbol, which will also be used over the elements of some abstract
  bases of real valued function that we will need later on (especially
  in Section~\ref{sec:domain_temporal_discretization}).
\end{remark}

Given a step function $\lub_{i \in I} b_i \chi_{O_i}: X \to D$, for
every $x \in X$ we define $I_x \defeq \setbarTall{i \in I}{x \in
  O_i}$. We observe that:
\begin{equation}
  \label{eq:step_fun_eval}
  \forall x \in X: \quad (\lub_{i \in I} b_i \chi_{O_i})(x) = \lub_{x
    \in O_i} b_i = \lub_{i
    \in I_x} b_i.
\end{equation}

\begin{proposition}
  For any pair of step functions $\lub_{i \in I} b_i \chi_{O_i}$ and
  $\lub_{j \in J} b'_j \chi_{O'_j}$, we have:

  \begin{equation}
    \label{eq:step_fun_sqsubset_formulation}
  \lub_{i \in I} b_i \chi_{O_i} \sqsubseteq  \lub_{j \in J} b'_j
  \chi_{O'_j} \iff \forall i_0 \in I: O_{i_0} \subseteq (\lub_{j \in J} b'_j
  \chi_{O'_j})^{-1}(\upperSet{b_{i_0}}). 
\end{equation}
\end{proposition}

\begin{proof}
  To prove the $(\Rightarrow)$ direction, for any given $i_0 \in I$
  and $x \in O_{i_0}$, we derive:
  \begin{eqnarray*}
    b_{i_0} & \sqsubseteq & \lub_{i \in I_x} b_i \\
    (\text{by equation~\eqref{eq:step_fun_eval}}) & = & (\lub_{i \in
                                                        I} b_i
                                                        \chi_{O_i})(x)\\
    (\text{by assumption} \lub_{i \in I} b_i \chi_{O_i} \sqsubseteq  \lub_{j \in J} b'_j
  \chi_{O'_j}) & \sqsubseteq & (\lub_{j \in J} b'_j
  \chi_{O'_j})(x),
  \end{eqnarray*}
  which implies that $O_{i_0} \subseteq (\lub_{j \in J} b'_j
  \chi_{O'_j})^{-1}(\upperSet{b_{i_0}})$.

  To prove the $(\Leftarrow)$ direction, we fix an $i_0 \in I$. From
  the assumption
  $O_{i_0} \subseteq (\lub_{j \in J} b'_j
  \chi_{O'_j})^{-1}(\upperSet{b_{i_0}})$ we deduce that
  $\forall x \in O_{i_0}: b_{i_0} \chi_{O_{i_0}}(x) = b_{i_0}
  \sqsubseteq (\lub_{j \in J} b'_j \chi_{O'_j})(x)$. This, combined
  with the fact that
  $\forall x \in X \setminus O_{i_0}: b_{i_0} \chi_{O_{i_0}}(x) =
  \bot$, implies that
  $b_{i_0} \chi_{O_{i_0}} \sqsubseteq \lub_{j \in J} b'_j
  \chi_{O'_j}$. As $i_0$ was chosen arbitrarily, by taking the
  supremum we obtain
  $\lub_{i \in I} b_i \chi_{O_i} \sqsubseteq \lub_{j \in J} b'_j
  \chi_{O'_j}$.
\end{proof}

\begin{definition}[Stable space]
  \label{def:stable}
A core-compact space $(X, \Omega(X))$ is called \emph{stable} if
$U \ll V$ and $U \ll V'$ imply $U \ll V \cap V'$, for any
$U, V, V' \in \Omega(X)$.   
\end{definition}

When $X$ is stable, based
on~\parencite[Lemma~1 and Proposition~5]{Erker_et_al:way_below:1998},
we know that:
  \begin{equation}
    \label{eq:step_fun_way_below_formulation}
  \lub_{i \in I} b_i \chi_{O_i} \ll  \lub_{j \in J} b'_j
  \chi_{O'_j} \iff \forall i \in I: O_i \ll (\lub_{j \in J} b'_j
  \chi_{O'_j})^{-1}(\wayAboves{b_i}). 
\end{equation}
We point out that, for the $(\Leftarrow)$ implication to hold, it
suffices for $X$ to be
core-compact~\parencite[Lemma~1]{Erker_et_al:way_below:1998}. Our aim
is, however, to develop a framework for any arbitrary topological
space $X$. In fact, our focus will be on spaces $X$ which are not
core-compact. Hence, with~\eqref{eq:step_fun_sqsubset_formulation}
and~\eqref{eq:step_fun_way_below_formulation} in mind, we define an
approximation order $\precAbsBas$ on the step functions in
$\mathcal{B}$ (defined in~\eqref{eq:def_B}) as follows:
\begin{equation}
  \label{eq:def_precAbsBas}
  \lub_{i \in I} b_i \chi_{O_i} \precAbsBas  \lub_{j \in J} b'_j
  \chi_{O'_j} \iff \forall i \in I: O_i \subseteq (\lub_{j \in J} b'_j
  \chi_{O'_j})^{-1}(\wayAboves{b_i}). 
\end{equation}
For the relation $\precAbsBas$ to be well-defined, we must show
that~\eqref{eq:def_precAbsBas} is independent of how step functions
are presented. The definition is clearly independent of how
$\lub_{j \in J} b'_j \chi_{O'_j}$ is presented. Hence, we must show
that, if
$\lub_{i \in I_1} b_i \chi_{O_i} = \lub_{i \in I_2} c_i \chi_{U_i}$,
then
$\lub_{i \in I_1} b_i \chi_{O_i} \precAbsBas \lub_{j \in J} b'_j
\chi_{O'_j} \iff \lub_{i \in I_2} c_i \chi_{U_i} \precAbsBas \lub_{j
  \in J} b'_j \chi_{O'_j}$. This follows from
item~\ref{item:prop:sqsubseteq_prec_implies_prec} of the following
Proposition:

\begin{proposition}
  \label{prop:prec_step_funs_basic_properties}
  For all $\theta, \theta', \theta'' \in \mathcal{B}$, we have:
  \begin{enumerate}[label=(\roman*)]
  \item \label{item:prop:sqsubseteq_prec_implies_prec} $\theta \sqsubseteq \theta' \precAbsBas \theta'' \implies
    \theta \precAbsBas \theta''$.
    
  \item \label{item:prop:prec_sqsubseteq_implies_prec}
    $\theta \precAbsBas \theta' \sqsubseteq \theta'' \implies
    \theta \precAbsBas \theta''$.

  \item \label{item:prop:prec_implies-sqsubseteq}
    $\theta \precAbsBas \theta' \implies \theta \sqsubseteq \theta'$.
    
  \item \label{item:prop:step_fun_sup_prec}
    $(\theta \precAbsBas \theta'') \wedge (\theta' \precAbsBas
    \theta'') \implies \lub \set{\theta, \theta'} \precAbsBas
    \theta''$.
  \end{enumerate}
\end{proposition}

\begin{proof}
  Let us assume that 
    $\theta = \lub_{i \in I} b_i \chi_{O_i}$, $\theta' = \lub_{j \in J} b'_j
    \chi_{O'_j}$, and $\theta'' =  \lub_{k \in K} b''_k \chi_{O''_k}$.

    \begin{enumerate}[label=(\roman*)]

    \item From the assumption $\theta \sqsubseteq \theta'$, we deduce
      that
      $\forall i \in I: b_i \chi_{O_i} \sqsubseteq \lub_{j \in J} b'_j
      \chi_{O'_j}$. From~\eqref{eq:step_fun_eval}
      and~\eqref{eq:step_fun_sqsubset_formulation}, for any $i \in I$
      and $x \in O_i$, by setting
      $J_x \defeq \setbarTall{j \in J}{x \in O'_j}$, we have:
      \begin{equation}
        \label{eq:theta_double_prime_inverse_j_Jx}
        b_i \sqsubseteq \lub_{j \in J_x} b'_j
        \implies \bigcap_{j \in J_x} \wayAboves{b'_j} \subseteq
        \wayAboves{b_i} \implies (\theta'')^{-1}(\bigcap_{j \in J_x}
        \wayAboves{b'_j}) \subseteq (\theta'')^{-1}(\wayAboves{b_i}).
      \end{equation}
      On the other hand, from the assumption $\theta' \precAbsBas
      \theta''$, we obtain $\forall j \in J_x: O'_j \subseteq
      (\theta'')^{-1} (\wayAboves{b'_j})$. Hence:
      \begin{eqnarray*}
        \bigcap_{j \in J_x} O'_j & \subseteq & \bigcap_{j \in J_x}
                                               (\theta'')^{-1}(\wayAboves{b'_j})
        \\
        \left(\text{In general, } f^{-1}( \bigcap_{s \in S}X_s) =
        \bigcap_{s \in S}f^{-1}(X_s) \right) & = &
                                                   (\theta'')^{-1}(\bigcap_{j \in J_x} \wayAboves{b'_j}) \\
        (\text{By~\eqref{eq:theta_double_prime_inverse_j_Jx}}) & \subseteq & (\theta'')^{-1}(\wayAboves{b_i}).
      \end{eqnarray*}
      In other words
      $\forall x \in O_i: x \in (\theta'')^{-1}(\wayAboves{b_i})$,
      which implies that
      $O_i \subseteq (\theta'')^{-1}(\wayAboves{b_i})$. As $i$ was
      chosen arbitrarily, by~\eqref{eq:def_precAbsBas} we infer that
      $\theta \precAbsBas \theta''$.

    \item From $\theta \precAbsBas \theta'$
      and~\eqref{eq:def_precAbsBas} we deduce
      $\forall i \in I: O_i \subseteq
      (\theta')^{-1}(\wayAboves{b_i})$. On the other hand, from
      $\theta' \sqsubseteq \theta''$ we obtain
      $\forall i \in I: (\theta')^{-1}(\wayAboves{b_i}) \subseteq
      (\theta'')^{-1}(\wayAboves{b_i})$. Thus, we have
      $\forall i \in I: O_i \subseteq
      (\theta'')^{-1}(\wayAboves{b_i})$, which implies that
      $\theta \precAbsBas \theta''$.
      
    \item If $\theta \precAbsBas \theta'$, then
      by~\eqref{eq:def_precAbsBas} we have:

      \begin{eqnarray*}
        \forall i \in I: O_i \subseteq (\lub_{j \in J} b'_j
  \chi_{O'_j})^{-1}(\wayAboves{b_i}) & \implies & \forall i \in I:
        \forall x \in O_i:  b_i \ll (\lub_{j \in J} b'_j
                                                  \chi_{O'_j}) (x)\\
        & \implies & \forall i \in I: \forall x \in O_i:  b_i \sqsubseteq (\lub_{j \in J} b'_j
                     \chi_{O'_j}) (x)\\
        & \implies & \forall i \in I: b_i \chi_{O_i} \sqsubseteq \lub_{j \in J} b'_j
                     \chi_{O'_j}\\
        & \implies & \lub_{i \in I}b_i \chi_{O_i} \sqsubseteq \lub_{j \in J} b'_j
                     \chi_{O'_j}.
      \end{eqnarray*}

    \item By item~\ref{item:prop:prec_implies-sqsubseteq} and
      assumption of
      $(\theta \precAbsBas \theta'') \wedge (\theta' \precAbsBas
      \theta'')$, $\theta''$ is an upper bound of
      $\set{\theta, \theta'}$. Hence, since $D$ is bounded-complete,
      the function $\lub \set{\theta, \theta'}$ is
      well-defined. Without loss of generality, we may assume that the
      index sets $I$ and $J$ are disjoint, let $A \defeq I \cup J$,
      and define:
  \begin{equation}
    \label{eq:def_Y_alpha_d_alpha_general}
    \forall \alpha \in A: \quad Y_\alpha \defeq
    \left\{
      \begin{array}{ll}
        O_\alpha,& \text{if } \alpha\in I,\\
        O'_\alpha, & \text{if } \alpha\in J,
      \end{array}
    \right.
    \qquad 
    d_\alpha \defeq
    \left\{
      \begin{array}{ll}
        b_\alpha,& \text{if } \alpha\in I,\\
        b'_\alpha, & \text{if } \alpha\in J.
      \end{array}
    \right.
  \end{equation}
  It is straightfoward to verify that
  $\lub \set{\theta, \theta'} = \lub_{\alpha \in A} d_\alpha
  \chi_{Y_\alpha}$, which entails that $\lub \set{\theta, \theta'}$ is
  indeed a step function in
  $\mathcal{B}$. From~\eqref{eq:def_Y_alpha_d_alpha_general} and the
  assumption
  $(\theta \precAbsBas \theta'') \wedge (\theta' \precAbsBas
  \theta'')$, we deduce that
  $\forall \alpha \in A : Y_\alpha \subseteq
  (\theta'')^{-1}(\wayAboves{d_\alpha})$, which implies that
  $\lub \set{\theta, \theta'} \precAbsBas \theta''$.
\end{enumerate}
\end{proof}

\begin{proposition}
  \label{prop:prec_interpolation_single_case}
  For any $\theta_1, \theta_2 \in \mathcal{B}$ satisfying
  $\theta_1 \precAbsBas \theta_2$, there exists a step function
  interpolating them, {\ie},
  $\exists \hat{\theta} \in \mathcal{B}: \theta_1 \precAbsBas
  \hat{\theta} \precAbsBas \theta_2$.
\end{proposition}

\begin{proof}
  We first prove the claim for the case where $\theta_1$ is a
  single-step function. Hence, assume that $\theta_1 = b \chi_O$ and
  $\theta_2 = \lub_{j \in J} b'_j \chi_{O'_j}$. As before, for each
  $x \in O$, we let $J_x \defeq \setbarTall{j \in J}{x \in O'_j}$. As
  $J$ is a finite index set, then the collection
  $C \defeq \setbarTall{J_x}{x \in O}$ must be finite. If
  $b = \bot_D$, then by taking $\hat{\theta} \defeq \bot_{[X \to D]}$,
  the result follows. If $b \neq \bot_D$, then we have
  $C \neq \emptyset$ and we may write
  $C = \set{J_{x_1}, \ldots, J_{x_k}}$, for some $k \geq 1$ and
  $x_1, \ldots, x_k \in O$. For each $1 \leq \ell \leq k$, we define
  $\Omega_\ell \defeq \bigcap_{j \in J_{x_\ell}} O'_j$. If
  $b \neq \bot_D$, then we must have:
  \begin{equation}
    \label{eq:O_subset_Omega_ell}
  O \subseteq \bigcup_{1 \leq \ell \leq k} \Omega_\ell.    
  \end{equation}
  \noindent
  The assumption $b\chi_O \precAbsBas \theta_2$ entails that
  $\forall \ell \in \set{1, \ldots, k}: b \ll \theta_2(x_\ell) =
  \lub_{j \in J_{x_\ell}} b'_j$. Using the interpolation property of
  continuous domains
  (Lemma~\ref{lemma:interpolation_property_cont_domains}), for each
  $\ell \in \set{1, \ldots, k}$ we choose an element
  $b''_\ell \in B(D)$ satisfying:
  \begin{equation}
    \label{eq:b_ll_bpp_ell_interpolate}
    b \ll b''_\ell \ll \theta_2(x_\ell).    
  \end{equation}
  We now define
  $\hat{\theta} \defeq \lub_{1 \leq \ell \leq k} b''_\ell
  \chi_{\Omega_\ell}$. The fact that
  $b\chi_O \precAbsBas \hat{\theta}$ follows
  from~\eqref{eq:O_subset_Omega_ell}
  and~\eqref{eq:b_ll_bpp_ell_interpolate}. On the other hand, for any
  $1 \leq \ell \leq k$, we have
  $\forall x \in \Omega_\ell: \theta_2(x_\ell) \sqsubseteq
  \theta_2(x)$, which, combined
  with~\eqref{eq:b_ll_bpp_ell_interpolate} implies that
  $\forall x \in \Omega_\ell: b''_\ell \ll \theta_2(x)$. Thus,
  $\Omega_\ell \subseteq
  (\theta_2)^{-1}(\wayAboves{b''_\ell})$. Hence,
  $\hat{\theta} \precAbsBas \theta_2$.

  Now, we generalize the proof to any step function $\theta_1$. Thus,
  assume that $\theta_1 = \lub_{i\in I}b_i \chi_{O_i}$. For each
  $i \in I$, we obtain an interpolating step function $\hat{\theta}_i$
  such that
  $b_i \chi_{O_i} \precAbsBas \hat{\theta}_i \precAbsBas \theta_2$. By
  using
  Proposition~\ref{prop:prec_step_funs_basic_properties}~\ref{item:prop:step_fun_sup_prec},
  we deduce that
  $\theta_1 = \lub_{i \in I}b_i \chi_{O_i} \precAbsBas \lub_{i \in I}
  \hat{\theta}_i \precAbsBas \theta_2$.
\end{proof}

\begin{lemma}
  The pair $(\mathcal{B}, \precAbsBas)$ forms an abstract basis.
\end{lemma}

\begin{proof}
  To prove transitivity, let us assume that
  $\theta_1 \precAbsBas \theta_2 \precAbsBas \theta_3$. By
  Proposition~\ref{prop:prec_step_funs_basic_properties}~\ref{item:prop:prec_implies-sqsubseteq},
  from $\theta_2 \precAbsBas \theta_3$ we deduce that
  $\theta_2 \sqsubseteq \theta_3$. Hence, we have
  $\theta_1 \precAbsBas \theta_2 \sqsubseteq \theta_3$. From
  Proposition~\ref{prop:prec_step_funs_basic_properties}~\ref{item:prop:prec_sqsubseteq_implies_prec},
  we obtain $\theta_1 \precAbsBas \theta_3$.

  Next, we must prove the interpolation property for
  $(\mathcal{B}, \precAbsBas)$. Assume that
  $\theta_1 \precAbsBas \theta$ and $\theta_2 \precAbsBas \theta$. If
  we define $\theta_3 = \lub \set{\theta_1, \theta_2}$, then by
  Proposition~\ref{prop:prec_step_funs_basic_properties}~\ref{item:prop:step_fun_sup_prec}
  we have $\theta_3 \precAbsBas \theta$. By
  Proposition~\ref{prop:prec_interpolation_single_case}, we obtain
  another step function $\hat{\theta}$ such that
  $\theta_3 \precAbsBas \hat{\theta} \precAbsBas \theta$. Once again,
  from Proposition~\ref{prop:prec_step_funs_basic_properties} we
  deduce that $\theta_1 \precAbsBas \hat{\theta} \precAbsBas \theta$
  and $\theta_2 \precAbsBas \hat{\theta} \precAbsBas \theta$.  
\end{proof}

\begin{definition}[$\mathcal{W}$]
  \label{def:W}
  Let $\mathcal{W}$ denote the rounded ideal completion of
  $(\mathcal{B}, \precAbsBas)$.
\end{definition}
Let $i: {\mathcal{B}} \to {\mathcal{W}}$ be the map
$i(b) \defeq \setbarNormal{b' \in {\mathcal{B}}}{ b' \precAbsBas
  b}$. By~\cite[Proposition~2.2.22]{AbramskyJung94-DT}, we know that
${\mathcal{W}}$ is a continuous domain, for which $i(\mathcal{B})$
forms a basis. When $(X, \Omega(X))$ is second countable and
$(D, \sqsubseteq_D)$ is $\omega$-continuous, we can choose the bases
$B(X)$ and $B(D)$ to be countable. In this case, $\mathcal{B}$ is also
countable, and $\mathcal{W}$ is an $\omega$-continuous domain with
$i(\mathcal{B})$ as a countable basis.

\subsection{Relationship between $\mathcal{W}$ and $[X \to D]$}

In this section we demonstrate that, regardless of whether $X$ is
core-compact or not, the function space $[X \to D]$ is tightly linked
with the continuous domain $\mathcal{W}$ via an adjunction, also known
as a Galois connection. We briefly recall the concept of Galois
connection, but for a more comprehensive account, the reader may refer
to, {\eg},~\cite[Section~3.1.3]{AbramskyJung94-DT}.

\begin{definition}[Category $\Po$, Galois connection $F \dashv G$]
  We let $\Po$ denote the category of \acp{POSET} and monotonic
  maps. A Galois connection in the category $\Po$ between two
  \acp{POSET} $(C, \sqsubseteq_C)$ and $(D, \sqsubseteq_D)$ is a pair
  of monotonic maps:
  \begin{equation*}
    \begin{tikzcd}[column sep = large]
      D \arrow[r, "G", yshift = 1.1ex] & C \arrow[l, "F", "\top"', yshift =
      -1.1ex]
    \end{tikzcd}
  \end{equation*}
  such that:
  \begin{equation*}
  \forall x \in C: \forall y \in D: x \sqsubseteq_C G(y) \iff F(x)
  \sqsubseteq_D y.  
\end{equation*}
In this case, we call $F: C \to D$ the left adjoint and $G: D \to C$
the right adjoint, and write $F \dashv G$.
\end{definition}

We extend the order $\precAbsBas$ that was defined
in~\eqref{eq:def_precAbsBas} over step functions to all functions from
$X$ to $D$ as follows:
\begin{equation}
  \label{eq:def_separation_X_to_D}
  \forall f,g : X \to D: \quad f \precAbsBas g \iff \exists \theta_1,
  \theta_2 \in \mathcal{B}: f \sqsubseteq \theta_1 \precAbsBas
  \theta_2 \sqsubseteq g.
\end{equation}
Intuitively, $f \precAbsBas g$ if and only if $f$ and $g$ are
separated by two step functions $\theta_1$ and $\theta_2$ satisfying
$\theta_1 \precAbsBas \theta_2$. Based on this intuition, we refer to
the relation $\precAbsBas$ as the \emph{separation order}. We point
out that~\eqref{eq:def_separation_X_to_D} is indeed a consistent
extension of~\eqref{eq:def_precAbsBas}. Specifically, assume that $f$
and $g$ are two step functions in $\mathcal{B}$:
\begin{itemize}
\item If $f \precAbsBas g$ according to~\eqref{eq:def_precAbsBas},
  then by taking $\theta_1 = f$ and $\theta_2 = g$, the separation
  order of~\eqref{eq:def_separation_X_to_D} will also be satisfied.

\item If $f \precAbsBas g$ according
  to~\eqref{eq:def_separation_X_to_D}, then by referring to
  Proposition~\ref{prop:prec_step_funs_basic_properties}, it can be
  verified that~\eqref{eq:def_precAbsBas} also holds for $f$ and $g$.
\end{itemize}

Let us briefly recall the concept of a monotone section-retraction
pair. For a more detailed account the reader may refer
to~\parencite[Section~3.1.1]{AbramskyJung94-DT}. Assume that $D$ and
$E$ are two posets. A pair of maps $s: D \to E$ and $r: E \to D$ is
called a monotone section-retraction pair if $s$ and $r$ are monotone
and $r \circ s = \id_D$. In this case, $D$ is said to be a monotone
retract of $E$. It is straightforward to verify that if $s$ and $r$
form a section-retraction pair, then $s$ must be injective and $r$
must be surjective.

Let $X$ be an arbitrary topological space. For every
$f \in [X \to D]$, we define:
\begin{equation}
  \label{eq:def_f_star_X_to_D}
  f_* \defeq \setbarNormal{b \in \mathcal{B}}{b \precAbsBas f}.
\end{equation}

\begin{proposition}
  \label{prop:f_star_rounded_ideal_X_to_D}
  For every $f \in [X \to D]$, the set $f_*$ is a rounded ideal in
  $\mathcal{B}$.
\end{proposition}

\begin{proof}
  Consider a function $f \in [X \to D]$. Then:
  \begin{enumerate}[label=(\roman*)]
  \item $f_*$ is non-empty, because $\bot \in f_*$.
  \item $f_*$ is a lower set, because $\forall g,g' \in \mathcal{B}: g
    \sqsubseteq g' \precAbsBas f \implies g \precAbsBas f$.

  \item $f_*$ is directed: Assume that $g_1, g_2 \in \mathcal{B}$
    satisfy $g_1 \precAbsBas f$ and $g_2 \precAbsBas
    f$. By~\eqref{eq:def_separation_X_to_D}, there must exist
    $g'_1, g'_2 \in \mathcal{B}$ satisfying
    $g_1 \precAbsBas g'_1 \sqsubseteq f$ and
    $g_2 \precAbsBas g'_2 \sqsubseteq f$. Let us define
    $g, g' \in \mathcal{B}$ as $g \defeq \lub \set{g_1, g_2}$ and
    $g' \defeq \lub \set{g'_1, g'_2}$. By
    Proposition~\ref{prop:prec_step_funs_basic_properties}, we have
    $g \precAbsBas g' \sqsubseteq f$. Since
    $(\mathcal{B}, \precAbsBas)$ is an abstract basis, then it must
    have the interpolation property. Hence, there must exist a step
    function $h \in \mathcal{B}$ satisfying
    $g \precAbsBas h \precAbsBas g'$. Again, by
    Proposition~\ref{prop:prec_step_funs_basic_properties}, we have
    $g_1 \precAbsBas h \precAbsBas f$ and
    $g_2 \precAbsBas h \precAbsBas f$.
   \end{enumerate}
   Thus, the set $f_*$ is a rounded ideal.
 \end{proof}

 \begin{proposition}
   \label{prop:f_equals_lub_f_star_X_to_D}
  For every $f \in [X \to D]$, we have $f = \lub f_*$.
\end{proposition}

\begin{proof}
  It is clear that $f \sqsupseteq \lub f_*$. To prove the
  $\sqsubseteq$ direction, choose any elements $x \in X$ and
  $b \in (\wayBelows{f(x)} \cap B(D))$. By the interpolation property
  of continuous domains
  (Lemma~\ref{lemma:interpolation_property_cont_domains}), we choose
  another basis element $b' \in B(D)$ such that $b \ll b' \ll
  f(x)$. By Proposition~\ref{prop:Scott_wayaboves}, we know that
  $\wayAboves{b'}$ is Scott open. Hence,
  $f^{-1}(\wayAboves{b'}) \in \Omega(X)$. As $B(X)$ is assumed to be a
  base for $\Omega(X)$, there must exist an open set $O \in B(X)$ for
  which we have $x \in O \subseteq f^{-1}(\wayAboves{b'})$. Clearly,
  we have $b \chi_O \in \mathcal{B}$ and $b' \chi_O \in
  \mathcal{B}$. Furthermore,
  $b \chi_O \precAbsBas b' \chi_O \sqsubseteq f$.

Since $f(x) = \lub (\wayBelows{f(x)} \cap B(D))$ and $b$ was chosen
arbitrarily, we have:
\begin{equation*}
f(x) = \lub \setbarTall{b \chi_O (x)}{ b \chi_O \in f_*},  
\end{equation*}
which implies that
$f(x) \sqsubseteq \lub \setbarTall{g(x)}{g \in f_*}$. Finally, as $x$ was chosen arbitrarily, then we must
have $f \sqsubseteq \lub f_*$.
\end{proof}

In the other direction, for every rounded ideal
$\phi \in \mathcal{W}$, we define $\phi^* : X \to D$ by:
\begin{equation}
  \label{eq:def_phi_star_X_to_D}
  \phi^* \defeq \lub_{b \in \phi} b.
\end{equation}
\noindent
Recall that every rounded ideal is directed. Thus, as $[X \to D]$ is a
\ac{DCPO}, then $\phi^*$ is well-defined and continuous.

\begin{lemma}
  \label{lemma:sr_pair_X_to_D}
  The pair $(\cdot)_*$ and $(\cdot)^*$ form a monotone
  section-retraction pair between $[X \to D]$ and $\mathcal{W}$.
\end{lemma}

\begin{proof}
  By Proposition~\ref{prop:f_star_rounded_ideal_X_to_D}, the map
  $(\cdot)_*$ is a function from $[X \to D]$ to
  $\mathcal{W}$. Monotonicity of $(\cdot)_*$ follows directly
  from~\eqref{eq:def_f_star_X_to_D}. In the other direction, since
  $[X \to D]$ is a \ac{DCPO}, for each $\phi \in \mathcal{W}$, we have
  $\phi^* \in [X \to D]$. Monotonicity of $(\cdot)^*$ follows directly
  from~\eqref{eq:def_phi_star_X_to_D}. Finally, the fact that the two
  maps form a monotone section-retraction pair is a consequence of
  Proposition~\ref{prop:f_equals_lub_f_star_X_to_D}, because for each
  $f \in [X \to D]$, we have $f = (f_*)^*$. 
\end{proof}

\begin{theorem}[Galois connection]
  \label{thm:Galois_connection_X_to_D}
  The maps $(\cdot)^*$ and $(\cdot)_*$ form a Galois connection:

  \begin{equation*}
    \begin{tikzcd}[column sep = large]
      [X \to D] \arrow[r, "(\cdot)_*", yshift = 1.1ex] & \mathcal{W} \arrow[l, "(\cdot)^*", "\top"', yshift =
      -1.1ex]
    \end{tikzcd}
  \end{equation*}
  in the category $\Po$, in which, $(\cdot)_*$ is the right adjoint,
  and $(\cdot)^*$ is the left adjoint. Furthermore:
  \begin{enumerate}[label=(\roman*)]

  \item \label{item:epi_mono_X_to_D} The map $(\cdot)^*$ is an epimorphism,
    and $(\cdot)_*$ is a monomorphism.
    
  \item \label{item:Galois_id_X_to_D} $(\cdot)^* \circ (\cdot)_* = \id_{[X
      \to D]}$, {\ie},
    $\forall f \in [X \to D]: (f_*)^* = f$.


  \item \label{item:left_adjoint_Scott_cont_X_to_D} The left adjoint
    $(\cdot)^*$ is Scott continuous.
  \end{enumerate}
\end{theorem}

\begin{proof}
  To prove that the maps $(\cdot)^*$ and $(\cdot)_*$ form a Galois
  connection, we must show that:
  \begin{equation*}
    \forall \phi \in \mathcal{W}, f \in [X \to D]: \quad \phi \subseteq f_* \iff \phi^* \sqsubseteq f,
  \end{equation*}
  which is a straightforward consequence of the definitions of
  $(\cdot)^*$ and $(\cdot)_*$.

  \begin{itemize}
  \item Claims~\ref{item:epi_mono_X_to_D}
    and~\ref{item:Galois_id_X_to_D} follow from
    Lemma~\ref{lemma:sr_pair_X_to_D}.

\item For any adjunction between two \acp{DCPO}, the left adjoint is
  Scott
  continuous~\cite[Proposition~3.1.14]{AbramskyJung94-DT}. Claim~\ref{item:left_adjoint_Scott_cont_X_to_D}
  now follows from the fact that both $[X \to D]$ and $\mathcal{W}$
  are \acp{DCPO}.  
  \end{itemize}  
\end{proof}

\begin{corollary}
  \label{cor:right_adjoint_Scott_X_core_compact}
  If the right adjoint $(\cdot)_*$ is Scott continuous and $D$ is not
  a singleton, then $X$ must be core-compact.
\end{corollary}

\begin{proof}
  By Lemma~\ref{lemma:sr_pair_X_to_D} and
  Theorem~\ref{thm:Galois_connection_X_to_D}, the pair
  $\left( (\cdot)_*, (\cdot)^* \right)$ forms a monotone
  section-retraction, with a Scott continuous retraction map
  $(\cdot)^*$. When the section $(\cdot)_*$ is also Scott continuous,
  the \ac{DCPO} $[X \to D]$ becomes a continuous retract of the
  continuous domain
  $\mathcal{W}$. By~\cite[Theorem~3.1.4]{AbramskyJung94-DT}, any
  continuous retract of a continuous domain is also a continuous
  domain, hence $[X \to D]$ must be a continuous domain. By
  Theorem~\ref{thm:core_compact}, this means that $X$ must be
  core-compact.
\end{proof}

\subsection{Core-Compact $X$}
\label{subsec:core_compact_X}

By Theorem~\ref{thm:Galois_connection_X_to_D}, the continuous domain
$\mathcal{W}$ is a suitable substitute for $[X \to D]$ when $X$ is not
core-compact, which is the main focus of the current article. In this
section, however, we briefly explore the relationship between
$\mathcal{W}$ and $[X \to D]$ when $X$ is core-compact. In particular,
we show that the converse of
Corollary~\ref{cor:right_adjoint_Scott_X_core_compact} is not true in
general.

When $X$ is core-compact, by
Theorem~\ref{thm:Galois_connection_X_to_D}, $\mathcal{W}$ contains
$[X \to D]$ as a sub-poset. First, we show that, when $X$ is stable
(Definition~\ref{def:stable}) the separation order over step functions
in $\mathcal{B}$ is finer then the way-below relation over
$[X \to D]$:

\begin{proposition}
  \label{prop:step_fun_way_below_precAbsBas}
  Assume that $X$ is a stable core-compact space, and let
  $\lub_{i \in I} b_i \chi_{O_i}$ and
  $\lub_{j \in J} b'_j \chi_{O'_j}$ be two step functions in
  $\mathcal{B}$. Then:
  \begin{equation}
    \label{eq:step_fun_way_below_implies_prec}
    \lub_{i \in I} b_i \chi_{O_i}  \ll \lub_{j \in J} b'_j \chi_{O'_j}
    \implies \lub_{i \in I} b_i \chi_{O_i}  \precAbsBas \lub_{j \in J} b'_j \chi_{O'_j}.
  \end{equation}
\end{proposition}

\begin{proof}
  Claim~\eqref{eq:step_fun_way_below_implies_prec} follows
  from~\eqref{eq:step_fun_way_below_formulation}
  and~\eqref{eq:def_precAbsBas}.
\end{proof}

Next, we show that the way-below relation over $\mathcal{W}$ is finer
than that over $[X \to D]$:

\begin{lemma}
  \label{lemma:waybelow_comparison}
  Assume that $X$ is a stable core-compact space. Then, for any
  $f,g \in [X \to D]$, we have:
  \begin{equation*}
    f \ll g \implies f_* \ll g_*.   
  \end{equation*}
\end{lemma}

\begin{proof}
  By~\cite[Proposition~2]{Erker_et_al:way_below:1998}, the set of step
  functions in $\mathcal{B}$ forms a basis for the continuous domain
  $[X \to D]$. By applying the interpolation property of continuous
  domains (Lemma~\ref{lemma:interpolation_property_cont_domains})
  twice, we obtain two step functions $\lub_{i \in I} b_i \chi_{O_i}$
  and $\lub_{j \in J} b'_j \chi_{O'_j}$ which satisfy:
  \begin{equation}
    \label{eq:lemma_waybelow_comparison_step_funs_interpolating}
      f \sqsubseteq \lub_{i \in I} b_i \chi_{O_i} \ll
      \lub_{j \in J} b'_j \chi_{O'_j} \sqsubseteq g.
  \end{equation}
  According to Proposition~\ref{prop:step_fun_way_below_precAbsBas}, we
  must have:
  \begin{equation}
    \label{eq:lemma_waybelow_comparison_step_funs}
    \lub_{i \in I} b_i \chi_{O_i} \precAbsBas \lub_{j \in J} b'_j
    \chi_{O'_j}.
  \end{equation}
  From~\eqref{eq:lemma_waybelow_comparison_step_funs_interpolating} we
  deduce that:
  \begin{equation}
    \label{eq:lemma_waybelow_comparison_step_funs_subseteq}
    f_* \subseteq (\lub_{i \in I} b_i \chi_{O_i})_* \subseteq
    (\lub_{j \in J} b'_j \chi_{O'_j})_* \subseteq g_*.
  \end{equation}
  By~\cite[Proposition~2.2.22]{AbramskyJung94-DT}, the
  relations~\eqref{eq:lemma_waybelow_comparison_step_funs}
  and~\eqref{eq:lemma_waybelow_comparison_step_funs_subseteq} imply
  $f_* \ll g_*$.
\end{proof}

The converse of Lemma~\ref{lemma:waybelow_comparison}, however, is not
true in general. In other words, the way-below relation over
$\mathcal{W}$ can be strictly finer than that over $[X \to D]$.

\begin{example}
  \label{example:way_below_finer}
  Assume that $X=[-2,2]$ under the Euclidean topology and
  $D = \intvaldom[\R\lift]$. As such, $X$ is core-compact and
  stable. Let $O = (-1,1)$, $b = [1,3]$, and $b' = [2,2]$. Then, we
  have $b \chi_O \precAbsBas b' \chi_O$, which implies
  $(b \chi_O)_* \ll (b' \chi_O)_*$. But,
  by~\eqref{eq:step_fun_way_below_formulation},
  $b \chi_O \not \ll b' \chi_O$.
\end{example}

In particular, Example~\ref{example:way_below_finer} shows that the
left adjoint $(\cdot)^*$ does not preserve the order of
approximation. While $(b \chi_O)_* \ll (b' \chi_O)_*$ holds, we have
$\left( (b \chi_O)_* \right)^* = b \chi_O \not \ll b' \chi_O = \left(
  (b' \chi_O)_*
\right)^*$. By~\cite[Proposition~3.1.14]{AbramskyJung94-DT}, this
means that the right adjoint is not Scott continuous, and the converse
of Corollary~\ref{cor:right_adjoint_Scott_X_core_compact} is not true in
general, even when $X$ is core-compact.

Nevertheless, in some cases, the converses of
Corollary~\ref{cor:right_adjoint_Scott_X_core_compact} and
Lemma~\ref{lemma:waybelow_comparison} do hold. Of course, this is true
for some trivial cases, {\eg}, when $X$ is a singleton. The converses,
however, hold even for non-trivial cases. In fact, the Galois
connection of Theorem~\ref{thm:Galois_connection_X_to_D} can reduce to
an isomorphism:

\begin{example}
  Assume that $X = D = \N \cup \set{+\infty}$ under the Scott
  topology. We take the collection of sets of the form
  $O_n \defeq \setbarNormal{k \in \N}{k \geq n} \cup \set{+\infty}$,
  for $n \in \N$, as the base for the Scott topology on $X$, and take
  $\N$ as the (domain-theoretic) basis for $D$. In this case, the
  Galois connection of Theorem~\ref{thm:Galois_connection_X_to_D}
  reduces to an isomorphism. The reason is that, over Scott open
  subsets of $X$, the way-below relation $\ll$ and the subset relation
  coincide. It is straightforward to verify that $X$ is core-compact
  and stable.
\end{example}

\subsection{Relationship between $\mathcal{W}$ and Algebraic
  Completion of $[X \to D]$}

  Let us take the set $\mathcal{B}$, but instead of the relation
  $\precAbsBas$, we order the step functions in $\mathcal{B}$ under
  the order $\sqsubseteq$. Then, the rounded ideal completion of
  $(\mathcal{B}, \sqsubseteq)$ results in an algebraic domain
  $\mathcal{W}_{\mathrm{alg}}$~\cite[Proposition~2.2.22]{AbramskyJung94-DT}.

We define the maps $u: \mathcal{W} \to \mathcal{W}_{\mathrm{alg}}$ and
$\ell: \mathcal{W}_{\mathrm{alg}} \to \mathcal{W}$ as follows:
\begin{equation*}
  \left\{
    \begin{array}{ll}
      \forall \phi \in \mathcal{W}: & u(\phi) \defeq \bigcup
                                      \setbarNormal{\lowerSet{b}}{(b
                                      \in \mathcal{B}) \wedge (b_*
                                      \subseteq \phi)},\\
      \forall \psi \in \mathcal{W}_{\mathrm{alg}}: & \ell(\psi) \defeq \bigcup
                                      \setbarNormal{b_*}{(b
                                      \in \mathcal{B}) \wedge (\lowerSet{b}
                                      \subseteq \psi)},\\
    \end{array}
  \right.
\end{equation*}
It is straightforward to verify that $\ell \dashv u$, {\ie}, they form
a Galois connection as follows:

  \begin{equation*}
    \begin{tikzcd}[column sep = large]
      \mathcal{W} \arrow[r, "u", yshift = 1.1ex] & \mathcal{W}_{\mathrm{alg}} \arrow[l, "\ell", "\top"', yshift =
      -1.1ex],
    \end{tikzcd}
  \end{equation*}
  \noindent
  and $\ell$ is surjective. In general, the \ac{DCPO} $\mathcal{W}$ is
  not algebraic. Hence, in general, $\ell$ does not preserve the order
  of approximation, and $u$ is not Scott continuous.

  We know from~\cite[Proposition~3.1.6]{AbramskyJung94-DT} that, when
  $X$ is core-compact, $[X \to D]$ is a retract of
  $\mathcal{W}_{\mathrm{alg}}$. Similar to
  Lemma~\ref{lemma:waybelow_comparison}, one may prove that the
  way-below relation over $\mathcal{W}_{\mathrm{alg}}$ is finer than
  that over $\mathcal{W}$. As the \ac{DCPO} $\mathcal{W}$ is not
  always algebraic, the way-below relation over
  $\mathcal{W}_{\mathrm{alg}}$ can be strictly finer than that over
  $\mathcal{W}$. Hence, the continuous domain $\mathcal{W}$ is always
  in between $[X \to D]$ and its algebraic completion
  $\mathcal{W}_{\mathrm{alg}}$, and the inclusions can be strict.


\section{A Domain for Temporal Discretization}
\label{sec:domain_temporal_discretization}

Consider the differential equation $y'(t) = f(y(t))$ from the
\ac{IVP}~\eqref{eq:main_ivp}. By integrating both sides, we obtain
$y( t+h) = y(t) + \int_t^{t+h} f(y(\tau)) \md \tau$, for all
$t \in [0,a]$ and $h \in [0,a-t]$. This can be written as:
\begin{equation}
  \label{eq:y_t_h_i}
  y(t+h) =
  y(t) + i( t, h),
\end{equation}
in which the integral $i( t, h)$ represents the dynamics of the solution
from $t$ to $t+h$. Thus, a general schema for validated solution of
\ac{IVP}~\eqref{eq:main_ivp} may be envisaged as follows:

\begin{enumerate}[label=(\roman*)]

\item For some $k \geq 1$, consider the partition
  $Q=(q_0, \ldots, q_k)$ of the interval $[0,a]$.

\item Let $Y(0) \defeq (0, \ldots,0)$.

\item \label{item:gen_schema_iteration} For each
  $j \in \set{0, \ldots, k-1}$ and $h \in (0,q_{j+1} - q_j]$:
  \begin{equation}
    \label{eq:Y_h_I_j}
    Y(q_j+h) \defeq Y(q_j) + I( q_j, h),
  \end{equation}
  where $I( q_j, h)$ is an interval enclosure of the integral factor
  $i( q_j, h)$ from equation~\eqref{eq:y_t_h_i}. The operator $I$, in
  general, depends on several parameters, including (enclosures of)
  the vector field and its derivatives, the enclosure $Y(q_j)$, the
  index $j$, etc.
  
\end{enumerate}
In~\eqref{eq:Y_h_I_j}, the operator `$+$' denotes interval addition,
and for the method to be validated, the term $I( q_j, h)$ must account
for all the inaccuracies, {\eg}, floating-point error, truncation
error, etc.

The schema is indeed a general one which encompasses various validated
approaches to \ac{IVP} solving in the literature, most notably, Euler
methods
of~\parencite{EdalatPattinson2006-Euler-PARA,Edalat_Farjudian_Mohammadian_Pattinson:2nd_Order_Euler:2020:Conf},
and Runge-Kutta methods
of~\parencite{Marciniak:Selected_Interval_Methods:2009,AlexandreDitSandretto:Validated_Runge_Kutta:2016}. 

In step~\ref{item:gen_schema_iteration} of the schema, the solver
moves forward in time, from $q_j$ to $q_{j+1}$. This requires keeping
the state, {\ie}, the solution up to the partition point $q_j$, and
referring to this state in iteration~$j$. As such, the schema has an
imperative style. This is in contrast with the functional style
adopted in language design for real number computation. For instance,
the languages designed
in~\parencite{Escardo96-tcs,Farjudian:Shrad:2007,DiGianantonio_Edalat-PCDF:2013}
for computation over real numbers and real functions are functional
languages based on lambda calculus, with their denotational semantics
provided by domain models.

In a functional framework, the solution of the
\ac{IVP}~\eqref{eq:main_ivp} is obtained as the fixpoint of a
higher-order operator. Domain models are particularly suitable for
fixpoint computations of this type. For Picard method of \ac{IVP}
solving, fixpoint formulations have been
obtained~\parencite{Edalat_Lieutier:Domain_Calculus_One_Var:MSCS:2004,Edalat_Pattinson2007-LMS_Picard}. For
Euler and Runge-Kutta methods, however, such fixpoint formulations do
not exist in the literature. This is because the commonly used domain
models in real number computation are not suitable for temporal
discretization of differential equations. Let us briefly expand on
this claim.

A straightfoward way of obtaining a fixpoint formulation for the above
general schema is to define a functional $\Phi$ over interval
functions as follows:
\begin{equation*}
  \Phi(Y)(x) \defeq \left\{
      \arrayoptions{2ex}{1.3}
      \begin{array}{ll}
        ( 0, \ldots, 0), & \text{if } x = 0,\\
        Y( q_j) + I( q_j, x-q_j), & \text{if } q_j < x \leq q_{j+1}.
      \end{array}
      \right.
\end{equation*}
The fixpoint of this operator (if it exists) will be the right
choice. The problem is that, the enclosures obtained by applying
$\Phi$ do not have upper (respectively, lower) semi-continuous upper
(respectively, lower) bounds and, by
Proposition~\ref{prop:D_0_n:Upper_Lower_SC}, the domain models used
in,
{\eg},~\parencite{Edalat_Lieutier:Domain_Calculus_One_Var:MSCS:2004,Edalat_Pattinson2007-LMS_Picard},
are not applicable. As a result, for a functional definition of Euler
and Runge-Kutta operators, we need a new domain model which is
different from, {\eg}, $D_n^{(0)}([0,a])$. To that end, we take the
following observations as general guidelines:

\begin{enumerate}[label=(\arabic*)]

\item The bounds of enclosures generated by $\Phi$ may lose their
  semi-continuity only over the partition points $q_0, \ldots, q_k$.

\item The integral operator $I$ typically generates enclosures with
  bounds that are continuous within each half-open interval
  $(q_j,q_{j+1}]$. This is true of the relevant validated methods of
  the literature, {\eg}, the Euler operators
  of~\parencite{EdalatPattinson2006-Euler-PARA,Edalat_Farjudian_Mohammadian_Pattinson:2nd_Order_Euler:2020:Conf},
  and Runge-Kutta operators
  of~\parencite{Marciniak:Selected_Interval_Methods:2009,AlexandreDitSandretto:Validated_Runge_Kutta:2016}.
\end{enumerate}

In essence, we must relax the semi-continuity requirement on the
bounds of enclosures, and may only require \emph{left} semi-continuity
at partition points. This relaxation of the requirement necessitates a
novel approach. The reason is that, while the \ac{POSET}
$D_n^{(0)}([0,a])$ of semi-continuous enclosures forms an
$\omega$-continuous domain, the \ac{POSET} of left semi-continuous
enclosures is not even continuous (Corollary~\ref{cor:D_Q_not_cont}).

\subsection{Left Semi-Continuous Maps and Enclosures}
\label{subsec:left_semi_cont_enclosures}

Let $[-K,K]^\uparrow$ denote the \ac{POSET} with carrier set $[-K,K]$
ordered by $\forall x, y \in [-K,K]: x \sqsubseteq y \iff x \leq
y$. Similarly, let $[-K,K]^\downarrow$ denote the \ac{POSET} with
carrier set $[-K,K]$ ordered by
$\forall x, y \in [-K,K]: x \sqsubseteq y \iff x \geq y$. Both
$[-K,K]^\uparrow$ and $[-K,K]^\downarrow$ are $\omega$-continuous
domains, which are non-algebraic when $K > 0$. In both cases, the set
$\Q \cap [-K,K]$ is a basis.

\begin{proposition}
  Assume that $(X, \Omega(X))$ is a topological space. Then,
  $f: X \to [-K,K]$ is:
  \begin{enumerate}[label=(\roman*)]
  \item upper semi-continuous $\iff f \in [ X \to [-K,K]^\downarrow]$.

  \item lower semi-continuous $\iff f \in [ X \to [-K,K]^\uparrow]$.
    
  \end{enumerate}
\end{proposition}

\begin{proof}
  The Scott open subsets of $[-K,K]^\downarrow$ are $[-K,K]$ and the
  collection $\setbarNormal{[-K,x)}{-K \leq x \leq K}$. Similarly, the
  Scott open subsets of $[-K,K]^\uparrow$ are $[-K,K]$ and the
  collection $\setbarNormal{(x,K]}{-K \leq x \leq K}$. The proof now
  follows from definition of upper/lower semi-continuity and
  Proposition~\ref{prop:Scott_wayaboves}.
\end{proof}

To solve the \ac{IVP}~\eqref{eq:main_ivp}, we consider interval
functions of type $f: [0,a] \to \intvaldom[{[-K,K]^n}]$. For Picard
method, it suffices to consider the Euclidean topology over
$[0,a]$~\parencite{Edalat_Lieutier:Domain_Calculus_One_Var:MSCS:2004,Edalat_Pattinson2007-LMS_Picard}. Under
the Euclidean topology, the interval $[0,a]$ is a core-compact
space. According to Proposition~\ref{prop:D_0_n:Upper_Lower_SC}, by
considering the Euclidean topology over $[0,a]$, we obtain enclosures
with upper and lower semi-continuous bounds. Based on our previous
explanations, however, for methods that proceed based on temporal
discretization, we may only require semi-continuity from the left.

Consider the set $O \defeq \setbarNormal{(a,b]}{ a, b \in \R}$ of left
half-open intervals of real numbers. The collection $O$ forms a base
for the so-called upper limit topology over
$\R$~\parencite{Sorgenfrey:1947:paracompact}. At an abstract level,
this topology is sufficient to capture left semi-continuity. As we are
laying down the foundation for an effective framework, however, we
work with a coarser variant of the upper limit topology. We consider
the set $O_{\Q} \defeq \setbarNormal{(a,b]}{ a, b \in \Q}$ of left
half-open intervals with rational end-points. The collection $O_{\Q}$
forms a base for what we refer to as the \emph{rational upper limit
  topology}.

Let $\R_{(\Q]}$ denote the topological space with $\R$ as the carrier
set under the rational upper limit topology. For any $X \subseteq \R$,
we let $X_{(\Q]} \defeq (X , \tau_{(\Q]})$ denote the topological
space with carrier set $X$ and the topology $\tau_{(\Q]}$ inherited by
$X$ as a subspace of $\R_{(\Q]}$. In contrast to the upper limit
topology, the rational upper limit topology is second-countable,
hence, strictly coarser. In fact, given any interval $[x,y]$, with
$x < y$, and an irrational point $r \in (x,y)$, the half-open interval
$(x,r]$ is open in the upper limit topology over $[x,y]$, but not in
the rational upper limit topology.

\begin{definition}
  \label{def:left_semi_continuity}
  We say that $f: X \to \R$ is (rational) left upper (respectively,
  lower) semi-continuous at $x_0 \in X$ if it is upper (respectively,
  lower) semi-continuous at $x_0$ with respect to the topology
  $\tau_{(\Q]}$. We drop the qualifier `rational' for brevity, and
  simply write left upper/lower semi-continuous. In particular, we say
  that a function $f: X \to \R$ is left upper (respectively, lower)
  semi-continuous if and only if it is left upper (respectively,
  lower) semi-continuous at every point $x_0 \in X$.
\end{definition}

We define:
\begin{equation*}
  \left\{
    \begin{array}{l}
      \mathcal{U}_{\Q} \defeq [ [0,a]_{(\Q]} \to [-K,K]^\downarrow],\\
      \mathcal{L}_{\Q} \defeq [ [0,a]_{(\Q]} \to [-K,K]^\uparrow].
    \end{array}
  \right.  
\end{equation*}
\noindent
It is straightfoward to prove
the following:
\begin{proposition}
  \label{prop:left_up_low_sc:domain_formulation}
  Assume that $f: [0,a] \to [-K,K]$. Then:
  \begin{enumerate}[label=(\roman*)]
  \item $f$ is left upper semi-continuous $\iff$
    $f \in \mathcal{U}_{\Q}$.

  \item $f$ is left lower semi-continuous $\iff$
    $f \in \mathcal{L}_{\Q}$.

  \end{enumerate}
\end{proposition}

To solve the \ac{IVP}~\eqref{eq:main_ivp} using temporal
discretization, we consider the following function space:
\begin{equation}
  \label{eq:def_D_Q}
  {\mathcal{D}}_{\Q} \defeq [ [0,a]_{(\Q]} \to \intvaldom[{[-K,K]^n}]],
\end{equation}
which we refer to as the \emph{\ac{POSET} of left semi-continuous
  enclosures}. We first observe a counterpart of
Proposition~\ref{prop:D_0_n:Upper_Lower_SC}:
\begin{proposition}
  A function
  $f \equiv (f_1, \ldots, f_n ): [0,a] \to \intvaldom[{[-K,K]^n}]$ is
  in ${\mathcal{D}}_{\Q}$ if and only if:
  \begin{equation*}
    \forall j \in \set{1, \ldots, n}: \left( \uep[f_j] \in
      \mathcal{U}_{\Q} \right) \wedge  \left( \lep[f_j] \in
      \mathcal{L}_{\Q} \right).
  \end{equation*}
\end{proposition}
\noindent
Next, we prove that none of the \acp{DCPO} ${\mathcal{D}}_{\Q}$,
$\mathcal{U}_{\Q}$, or $\mathcal{L}_{\Q}$ is continuous:

\begin{proposition}
  \label{prop:0_a_not_core_compact}
  The topological space $[0,a]_{(\Q]}$ is not core-compact.
\end{proposition}

\begin{proof}
  It suffices to prove the following claim:
  \begin{equation}
    \label{eq:not_core_compact_single_half_intervals}
    \forall (x,y], (s,t] \subseteq [0,a]: \quad (x,y] \ll (s,t]
    \implies (x,y] = \emptyset,
  \end{equation}
  in which $x,y,s,t \in \Q$. To
  prove~\eqref{eq:not_core_compact_single_half_intervals}, let us
  assume that $(x,y] \neq \emptyset$, with $0 \leq x < y \leq a$. We
  let $z$ be a rational number in $(x,y)$ and define the following
  increasing chain of open sets:
  \begin{equation*}
    \forall n \in \N: \quad A_n \defeq (s,z] \cup ( z +
    \frac{t-z}{2^{n+1}}, t].
  \end{equation*}
  We have $(s,t] \subseteq \bigcup \setbarNormal{A_n}{n \in \N}$, but
  $\forall n \in \N: (x,y] \not \subseteq A_n$. Thus, $(x,y] \not \ll
  (s,z]$, which is a contradiction.
\end{proof}

\begin{remark}
  It can be shown, with a similar proof, that the upper limit topology
  on $[0,a]$ is not core-compact
  either. See~\parencite[Example~5.2.14]{Goubault-Larrecq:Non_Hausdorff_topology:2013}
  for more examples.
\end{remark}

\begin{corollary}
  \label{cor:D_Q_not_cont}
  None of the \acp{POSET} ${\mathcal{D}}_{\Q}$, $\mathcal{U}_{\Q}$, or
  $\mathcal{L}_{\Q}$ is continuous.
\end{corollary}

\begin{proof}
  This follows from Proposition~\ref{prop:0_a_not_core_compact} and
  Theorem~\ref{thm:core_compact}.
\end{proof}

As a result, we follow the construction of the previous section using
abstract bases. To that end, we consider the following countable
bases:

\begin{itemize}
\item the base $O_{\Q}$ of left half-open intervals with rational
  end-points for the rational upper limit topology on $[0,a]$;
\item the basis $B_{\intvaldom[{[-K,K]^n}]}$ of hyper-rectangles with
  rational coordinates for $\intvaldom[{[-K,K]^n}]$;
\item the basis $\Q \cap [-K,K]$ of rational numbers in $[-K,K]$ for
  both $[-K,K]^\downarrow$ and $[-K,K]^\uparrow$.
\end{itemize}
Using these bases, we obtain the abstract bases of step functions
${\mathcal{B}}_{\mathcal{D}}$, ${\mathcal{B}}_{\mathcal{U}}$, and
${\mathcal{B}}_{\mathcal{L}}$, corresponding to the \acp{DCPO}
${\mathcal{D}}_{\Q}$, $\mathcal{U}_{\Q}$, and $\mathcal{L}_{\Q}$,
respectively, according to~\eqref{eq:def_B}
and~\eqref{eq:def_precAbsBas}. The abstract bases
${\mathcal{B}}_{\mathcal{U}}$ and ${\mathcal{B}}_{\mathcal{L}}$,
although ordered under different---in fact, opposite---orders, have
the same carrier set, which we denote by ${\mathcal{A}}_{\Q}$.

\begin{remark}
  In the construction of step functions in
  ${\mathcal{B}}_{\mathcal{D}}$, ${\mathcal{B}}_{\mathcal{U}}$, and
  ${\mathcal{B}}_{\mathcal{L}}$, we restrict the parameters to
  rational numbers to obtain an effective structure. Similar results
  can be obtained by replacing $\Q$ with any other countable dense
  subset of $\R$ which has the necessary effective structure, {\eg},
  is effectively enumerable, has a decidable order $\leq$, etc. For
  instance, the set $\D$ of dyadic numbers is an important case that
  is indeed used in our experiments (Section~\ref{sec:Experiments})
  which are implemented using arbitrary-precision interval
  arithmetic. Nonetheless, to keep the presentation simple, we stay
  focused on rational numbers.
\end{remark}

The set ${\mathcal{A}}_{\Q}$ of step functions has a fairly simple
description. Assume that $P = (p_0, \ldots, p_k)$ is a partition of
$[0,a]$ with rational partition points, {\ie},
$P \in {\mathcal{P}}_{\Q}$. We say that $f: [0,a] \to \Q$ is a
\emph{rational} $P$-function if, for some constants
$\set{c_0, \ldots, c_k} \subseteq \Q$:
\begin{equation}
  \label{eq:rational_P_fun}
  f(0) = c_0 \wedge \forall i \in \set{1, \ldots, k}: \restrictTo{f}{(p_{i-1},p_i]}  =
  \lambda x. c_i.
\end{equation}
\noindent
Thus, $f$ is left-continuous, but does not have to be right-continuous
at the partition points, and over $(p_{i-1},p_i]$, it is a constant
function. The set ${\mathcal{A}}_{\Q}$ consists of all rational
$P$-functions, with $P$ ranging over ${\mathcal{P}}_{\Q}$:
\begin{equation*}
  {\mathcal{A}}_{\Q} = \setbarTall{f: [0,a] \to [-K,K]}{\exists P \in {\mathcal{P}}_{\Q}: f \text{ is a rational $P$-function}}.
\end{equation*}

For each function $h \in {\mathcal{A}}_{\Q}$, the partition points are
rational numbers. It can also be verified that the functions in
$\mathcal{U}_{\Q}$ and $\mathcal{L}_{\Q}$ cannot have jumps at
irrational numbers.

\begin{example}
  Let $x_0$ be an irrational number in $(0,1)$, {\eg},
  $x_0 \defeq 1/\sqrt{2}$. Assume that $f: [0,1]_{(\Q]} \to [-1,1]^\downarrow$ is defined
  as:
  \begin{equation*}
    \forall x \in [0,1]: \quad f(x) \defeq
    \left\{
      \begin{array}{ll}
        0, & \text{ if } x \leq x_0,\\
        1, & \text{ if }  x_0 < x.\\
      \end{array}
    \right.
  \end{equation*}
  The set $[-1,1)$ is Scott open in $[-1,1]^\downarrow$, but
  $f^{-1}([-1,1)) = [0,x_0]$, which is not open in
  $[0,1]_{(\Q]}$. Hence, $f \not \in \mathcal{U}_{\Q}$. In other
  words, although $f$ is clearly upper semi-continuous with respect to
  the upper limit topology on $[0,1]$, it is not upper semi-continuous
  with respect to the rational upper limit topology.
\end{example}

\subsection{The $\omega$-Continuous Domain
  ${\mathcal{W}}_{\mathcal{D}}$}
\label{subsec:construction_W_Q}

For the bases ${\mathcal{B}}_{\mathcal{D}}$,
${\mathcal{B}}_{\mathcal{U}}$, and ${\mathcal{B}}_{\mathcal{L}}$, we
denote the separation order $\precAbsBas$ as
$\precAbsBas_{\mathcal{D}}$, $\precAbsBas_{\mathcal{U}}$, and
$\precAbsBas_{\mathcal{L}}$, respectively. The separation order
of~\eqref{eq:def_precAbsBas} and~\eqref{eq:def_separation_X_to_D} can
be reformulated, for (say) $\precAbsBas_{\mathcal{L}}$, as follows:
For any pair of functions $\theta_1, \theta_2: [0,a] \to [-K,K]$ in
${\mathcal{A}}_{\Q}$, we have:
  \begin{equation}
    \label{eq:prec_L_reformulated}
    \theta_1 \precAbsBas_{\mathcal{L}} \theta_2 \iff  \exists \delta
    > 0: \forall
    t \in [0,a]: \left( \theta_1(t) = -K \right) \vee \left( \theta_1(t) \leq
      \theta_2(t) - \delta \right).
\end{equation}
For any pair of (arbitrary) functions $f,g: [0,a] \to [-K,K]$, we
obtain:
\begin{equation*}
  f \precAbsBas_{\mathcal{L}} g \iff \exists \theta_1, \theta_2 \in  {\mathcal{A}}_{\Q}:
  f \leq \theta_1 \precAbsBas_{\mathcal{L}} \theta_2 \leq g.
\end{equation*}
\noindent
Thus, $f \precAbsBas_{\mathcal{L}} g$ if and only if $f \leq g$, and
the two functions are separated by two rational $P$-functions
$\theta_1$ and $\theta_2$ satisfying
$\theta_1 \precAbsBas_{\mathcal{L}} \theta_2$. Similarly, we have
$f \precAbsBas_{\mathcal{U}} g$ if and only if $f \geq g$, and the two
functions are separated by two rational $P$-functions. In fact, we
observe the following relations among the three separation orders:
\begin{equation}
  \label{eq:prec_D_prec_U_prec_L}
  \forall f,g : [0,a] \to [-K,K]^n: \quad f \precAbsBas_{\mathcal{D}} g \iff
  \forall j \in \set{ 1, \ldots, n}: (\uep[f_j]
  \precAbsBas_{\mathcal{U}} \uep[g_j]) \wedge (\lep[f_j]
  \precAbsBas_{\mathcal{L}} \lep[g_j]). 
\end{equation}
In what follows, for convenience, we state most of our results for
$\precAbsBas_{\mathcal{L}}$. The corresponding results can be stated
for $\precAbsBas_{\mathcal{U}}$ in a straightfoward manner.

Let $\theta_1, \theta_2 \in {\mathcal{A}}_{\Q}$ and assume that
$\theta_1 \precAbsBas_{\mathcal{L}} \theta_2$. We define
$\delta_0 \defeq \inf \setbarTall{\theta_2(t) - \theta_1(t)}{t \in
  [0,a]}$. From~\eqref{eq:prec_L_reformulated}, we know that
$\forall t \in [0,a]: \theta_1(t) < K$. If $\theta_1$ does not touch
the lower endpoint of the interval $[-K,K]$ either, then
$\delta_0 > 0$. Formally:
\begin{equation*}
  \left(  \forall t \in [0,a]: -K < \theta_1(t) \right)
  \implies \delta_0 > 0.
\end{equation*}
\noindent
If at some points $t \in [0,a]$ we have $\theta_1(t) = -K$, we may
still obtain the following useful inequalities by clipping the values
inside the $[-K,K]$ range, as long as we avoid $\theta_2(t) = -K$:
\begin{align}
  & \exists \delta > 0:
    \forall t \in [0,a]: \theta_1(t) \leq \max( -K, \theta_2(t) - \delta), \label{eq:theta_1_less_K_case}\\
  &\forall t \in [0,a]: -K < \theta_2(t) \implies \exists \delta > 0:
    \forall t \in [0,a]: \theta_1(t) \leq
    \theta_2(t) - \delta. \label{eq:theta_2_greater_minus_K_case}
\end{align}

For arbitrary functions $f,g: [0,a] \to [-K,K]$, however, the above do
not hold. Even if we restrict to left lower semi-continuous functions,
the above do not hold. For instance, assume that $f$ is the constant
function $f(t) \defeq -K$, while $g$ is defined as follows:
\begin{equation*}
  g(t) \defeq
  \left\{
    \begin{array}{ll}
      0, & \text{if } t = 0,\\
      \frac{K}{a}t - K, & \text{if } t \in (0,a].
    \end{array}
  \right.
\end{equation*}
Although $f \precAbsBas_{\mathcal{L}} g$ and
$\forall t \in [0,a]: -K < g(t)$, it is not true that
$\exists \delta > 0: \forall t \in [0,a]: f(t) \leq g(t) - \delta$. In
other words, \eqref{eq:theta_2_greater_minus_K_case} does not
hold. The reason is that $g$ can take values which are arbitrarily
close to the lower bound of the interval $[-K,K]$. If we exclude such
cases, then we obtain the counterparts
of~\eqref{eq:theta_1_less_K_case}
and~\eqref{eq:theta_2_greater_minus_K_case}:

\begin{proposition}
  \label{prop:rat_Q_funs_sep_inequalities}
  Assume that $f,g: [0,a] \to [-K,K]$ and $f \precAbsBas_{\mathcal{L}} g$. Then:
  \begin{enumerate}[label=(\roman*)]
  \item \label{item:f_less_K_eps}
    $\exists \delta > 0: \forall t \in [0,a]: f(t) \leq \max(
    -K, g(t) - \delta)$.
  \item \label{item:g_greater_minus_K_eps}
    $\left[\exists \epsilon > 0: \forall t \in [0,a]: -K + \epsilon <
      g(t) \right]
    \implies \exists \delta > 0: \forall t \in [0,a]: f(t) \leq g(t) -
    \delta$.
  \end{enumerate}
\end{proposition}
\begin{proof}
  As $f \precAbsBas_{\mathcal{L}} g$, for two rational $P$-functions
  $\theta_1, \theta_2 \in {\mathcal{A}}_{\Q}$, we must have
  $f \leq \theta_1 \precAbsBas_{\mathcal{L}} \theta_2 \leq g$. The
  proof of~\ref{item:f_less_K_eps} now follows
  from~\eqref{eq:theta_1_less_K_case}.
  
  To prove~\ref{item:g_greater_minus_K_eps}, we take a rational number
  $\hat{\epsilon} \in (0, \epsilon)$. We define
  $\hat{\theta}_2(t) \defeq \max( \theta_2(t), -K +
  \hat{\epsilon})$. Then, we must have
  $f \leq {\theta}_1 \precAbsBas_{\mathcal{L}} \hat{\theta}_2 \leq
  g$. Furthermore, it is clear that
  $\forall t \in [0,a]: -K < \hat{\theta}_2(t)$. Therefore,
  by~\eqref{eq:theta_2_greater_minus_K_case}, we have:
  $\exists \delta > 0: \forall t \in [0,a]: \theta_1(t) \leq
  \hat{\theta}_2(t) - \delta$, which implies that
  $f(t) \leq g(t) - \delta$.
\end{proof}

We will also need a kind of inverse of the above results, which is
formulated as follows:

\begin{proposition}
  \label{prop:rational_Q_fun_sep_arbit_f}
  Assume that $\theta \in {\mathcal{A}}_{\Q}$ and $f: [0,a] \to [-K,K]$
  is an arbitrary function. Then:

  \begin{enumerate}[label=(\roman*)]
  \item \label{item:theta_sep_f}
    $\left[ \exists \delta > 0: \forall t \in [0,a]: \theta(t) \leq \max( -K, f(t)
    - \delta) \right] \implies \theta \precAbsBas_{\mathcal{L}} f$.

  \item \label{item:f_sep_theta}
    $\left[\exists \delta > 0: \forall t \in [0,a]: \min( f(t) + \delta, K)
    \leq \theta(t)\right] \implies \theta \precAbsBas_{\mathcal{U}} f$.
  \end{enumerate}
\end{proposition}

\begin{proof}
  To prove~\ref{item:theta_sep_f}, let $\delta_0 \in (0, \delta)$ be a
  rational number. We define $\theta_2: [0,a] \to [-K,K]$ as follows:
  \begin{equation*}
    \forall t \in [0,a]: \quad \theta_2(t) \defeq
    \left\{
    \begin{array}{ll}
       \theta(t) + \delta_0, & \text{if } \theta(t) > -K,\\
      -K, &  \text{if } \theta(t) = -K.
    \end{array}
    \right.
  \end{equation*}
  It is straightfoward to verify that
  $\theta_2 \in {\mathcal{A}}_{\Q}$. Thus, we have
  $\theta \precAbsBas_{\mathcal{L}} \theta_2 \leq f$, which implies
  that $\theta \precAbsBas_{\mathcal{L}} f$. The proof
  of~\ref{item:f_sep_theta} is similar.
\end{proof}

We point out that Proposition~\ref{prop:rational_Q_fun_sep_arbit_f}
does not hold if $\theta$ is replaced with an arbitrary function
$g: [0,a] \to [-K,K]$. In fact, it does not hold even for left lower
semi-continuous functions. For instance, let $f,g: [0,1] \to [-3,3]$
be defined as follows:
\begin{equation*}
  \forall t \in [0,1]: \quad g(t) \defeq 
    \left\{
    \begin{array}{ll}
      -2, & \text{if $t \in (2^{-(2n+1)},2^{-2n}] $ for some $n \in
            \N$},\\
      0, & \text{if $t \in (2^{-(2n+2)},2^{-(2n+1)}] $ for some $n \in \N$},\\      
      0, & \text{if $t = 0$},
    \end{array}
    \right.
  \end{equation*}
  and $\forall t\in [0,1]: f(t) \defeq g(t) + 1$. Then, we have
  $\forall t \in [0,1]: g(t) \leq f(t) - 1$, but
  $g \precAbsBas_{\mathcal{L}} f$ does not hold, because the two
  cannot be separated by rational $P$-functions.

\begin{definition}[Non-degenerate enclosure]
  \label{def:Non-degenerate_enclosure}
  We call an enclosure $f: [0,a] \to \intvaldom[{[-K,K]^n}]$
  \emph{non-degenerate} if:
  \begin{equation*}
    \exists \epsilon > 0 : \forall j \in \set{ 1, \ldots, n}: \forall t
    \in [0,a]: \left( \lep[f_j](t) < K - \epsilon \right) \wedge \left( -K + \epsilon <
      \uep[f_j](t)  \right).
  \end{equation*}
\end{definition}

\begin{proposition}
  \label{prop:U_L:non-degen}
  Assume that $f,g : [0,a] \to [-K,K]$ and $K > 0$. Then:
  \begin{enumerate}[label=(\roman*)]
  \item $f \precAbsBas_{\mathcal{L}} g \implies \exists \epsilon > 0:
    \forall t \in [0,a]: f(t) < K-\epsilon$.
  \item $f \precAbsBas_{\mathcal{U}} g \implies \exists \epsilon > 0:
    \forall t \in [0,a]: -K+\epsilon < f(t)$.
  \end{enumerate}
\end{proposition}

\begin{proof}
  The results are straightfoward consequences of
  Proposition~\ref{prop:rat_Q_funs_sep_inequalities}.
\end{proof}

\begin{corollary}
  \label{corollary:prec_D_non-degen}
  Assume that $K > 0$ and $f, g: [0,a] \to \intvaldom[{[-K,K]^n}]$. If
  $f \precAbsBas_{\mathcal{D}} g$, then $f$ must be non-degenerate.
  \end{corollary}

  \begin{proof}
    This follows from Proposition~\ref{prop:U_L:non-degen}
    and~\eqref{eq:prec_D_prec_U_prec_L}.
  \end{proof}
  
Let $n \in \N \setminus \set{0}$ and let $\PS(\R^n)$ be the set of all
subsets of $\R^n$. For any real $K \geq 0$, we define the operator
$T_{\negthinspace K,n}: \PS(\R^n) \to \PS(\R^n)$ by:
\begin{equation}
  \label{eq:T_K_n}
  \forall X \in \PS(\R^n): \quad  T_{\negthinspace K,n}(X) \defeq X \cap [-K,K]^n.
\end{equation}
Based on Proposition~\ref{prop:rat_Q_funs_sep_inequalities} and
Corollary~\ref{corollary:prec_D_non-degen}, we obtain the following
result, which provides another perspective on the separation relation:

\begin{proposition}
  \label{prop:precAbsBas_symExpand}
  Assume that $f, g: [0,a] \to \intvaldom[{[-K,K]^n}]$. Then:
  \begin{equation*}
    f \precAbsBas_{\mathcal{D}} g \implies \exists \delta > 0: \forall t \in [0,a]: f(t) \sqsubseteq
    T_{\negthinspace K,n}( g(t) \symExpand \delta),
  \end{equation*}
  in which $\symExpand$ is the symmetric expansion of Definition~\ref{def:symmetric_expansion}.
\end{proposition}

The set ${\mathcal{B}}_{\mathcal{D}}$ can be effectively enumerated,
and the relation $\precAbsBas_{\mathcal{D}}$ is decidable over
${\mathcal{B}}_{\mathcal{D}}$. Thus, to obtain an effective framework,
we designate
$({\mathcal{B}}_{\mathcal{D}}, \precAbsBas_{\mathcal{D}})$ as an
abstract basis, over which we construct our effective domain model:

\begin{definition}[${\mathcal{W}}_{\mathcal{D}}$]
  \label{def:W_Q}
  Let ${\mathcal{W}}_{\mathcal{D}}$ denote the rounded ideal completion of
  $({\mathcal{B}}_{\mathcal{D}}, \precAbsBas_{\mathcal{D}})$.
\end{definition}

Let
$i_{\mathcal{D}}: {\mathcal{B}}_{\mathcal{D}} \to
{\mathcal{W}}_{\mathcal{D}}$ be the map
$i_{\mathcal{D}}(b) \defeq \setbarNormal{x \in
  {\mathcal{B}}_{\mathcal{D}}}{ x \precAbsBas_{\mathcal{D}}
  b}$. By~\cite[Proposition~2.2.22]{AbramskyJung94-DT}, we know that
${\mathcal{W}}_{\mathcal{D}}$ is an $\omega$-continuous domain, for
which $i_{\mathcal{D}}({\mathcal{B}}_{\mathcal{D}})$ forms a countable
basis. From Theorem~\ref{thm:Galois_connection_X_to_D}, we deduce that
${\mathcal{W}}_{\mathcal{D}}$ and the non-continuous \ac{DCPO}
${\mathcal{D}}_{\Q}$ are related via the following Galois connection:
\begin{equation}
  \label{eq:Galois_DQ_WD}
    \begin{tikzcd}[column sep = large]
      {\mathcal{D}}_{\Q} \arrow[r, "(\cdot)_*", yshift = 1.1ex] & {\mathcal{W}}_{\mathcal{D}} \arrow[l, "(\cdot)^*", "\top"', yshift =
      -1.1ex]
    \end{tikzcd}
  \end{equation}

  With the above Galois connection, we may delegate the \ac{IVP}
  solving from ${\mathcal{W}}_{\mathcal{D}}$ to
  ${\mathcal{D}}_{\Q}$---which is more suitable for performing the
  main computations---and then return to
  ${\mathcal{W}}_{\mathcal{D}}$. First, we formulate the following
  property:
  
\begin{definition}[\Ac{UIC}]
  We say that a function
  $\Phi: {\mathcal{D}}_{\Q} \to {\mathcal{D}}_{\Q}$ has the \ac{UIC}
  property if it is monotone and satisfies:
  \begin{equation*}
    \forall \phi \in {\mathcal{W}}_{\mathcal{D}}: \forall \delta > 0: \exists b_\delta \in
    \phi: \forall t \in [0,a]: \quad T_{\negthinspace K,n}\left( \Phi(\phi^*)(t) \symExpand
   \delta \right) \sqsubseteq
          \Phi(b_{\delta})(t).
  \end{equation*}
\end{definition}
Informally, for any given ideal $\phi \in {\mathcal{W}}_{\mathcal{D}}$ and accuracy
$\delta > 0$, there exists an element $b_\delta \in \phi$ for which
$\Phi(b_{\delta})$ approximates $\Phi(\phi^*)$ uniformly over $[0,a]$
to within $\delta$ accuracy. This condition is satisfied by the
operators which appear in Definition~\ref{def:Phi} (for second-order
Euler) and Definition~\ref{def:Phi_R} (for Runge-Kutta Euler)
methods. We expect this property to hold for any similar operator
defined according to the general schema. With that in mind, the
following lemma presents a procedure for obtaining Scott-continuous
operators to be used in fixpoint formulations:

\begin{lemma}
  \label{lem:F_Phi_general_result}
  Assume that $\Phi: {\mathcal{D}}_{\Q} \to {\mathcal{D}}_{\Q}$ has the \ac{UIC}
  property and define $F: {\mathcal{W}}_{\mathcal{D}} \to {\mathcal{W}}_{\mathcal{D}}$ by
  $\forall \phi \in {\mathcal{W}}_{\mathcal{D}}: F(\phi) \defeq \left(\Phi(\phi^*)
  \right)_*$. Then:
  \begin{enumerate}[label=(\roman*)]
  \item \label{item:F_Phi_general_alternative_formulation} $\forall \phi \in {\mathcal{W}}_{\mathcal{D}}: F(\phi) = \bigcup_{b \in \phi}
    \left( \Phi(b)\right)_*$.
  \item \label{item:F_Phi_general_Scott} $F: {\mathcal{W}}_{\mathcal{D}} \to {\mathcal{W}}_{\mathcal{D}}$ is Scott-continuous.
  \end{enumerate}
\end{lemma}

\begin{proof}  
  \begin{enumerate}[label=(\roman*)]
    \item We must
  prove that:

  \begin{equation*}
    \forall \phi \in {\mathcal{W}}_{\mathcal{D}}: \quad
      \left( \Phi(\phi^*)\right)_* = \bigcup_{b \in \phi}
      \left( \Phi(b)\right)_* .
    \end{equation*}
    The proof of the $\supseteq$ direction is relatively
    straightforward. For any $b \in \phi$, we have
    $b \sqsubseteq \phi^*$. By monotonicity of $\Phi$, we obtain
    $\Phi(b) \sqsubseteq \Phi(\phi^*)$, which, in turn, entails that
    $\left( \Phi(b)\right)_* \subseteq \left( \Phi(\phi^*)\right)_*$.
    
    Next, we prove the $\subseteq$ direction. Let us take an arbitrary
    $\hat{b} \in \left( \Phi(\phi^*)\right)_*$, which must satisfy
    $\hat{b} \precAbsBas_{\mathcal{D}} \Phi(\phi^*)$, and by
    Corollary~\ref{corollary:prec_D_non-degen}, is non-degenerate. By
    Proposition~\ref{prop:precAbsBas_symExpand}, we have:
      \begin{equation} 
       \label{eq:d_b_delta_hat}
        \exists \hat{\delta}>0: \forall t \in [0,a]:  \quad
        \hat{b}(t) \sqsubseteq  T_{\negthinspace K,n}(\Phi(\phi^*)(t) \symExpand \hat{\delta}).
      \end{equation}
      From the \ac{UIC} property of $\Phi$, we have:
      \begin{equation}
        \label{eq:d_Phi_delta}
        \forall \delta > 0: \exists b_{\delta} \in \phi: \forall t \in
        [0,a]: \quad
 T_{\negthinspace K,n}\left( \Phi(\phi^*)(t) \symExpand
   \delta \right) \sqsubseteq
          \Phi(b_{\delta})(t).
      \end{equation}
      If, in~\eqref{eq:d_Phi_delta}, we take
      $\delta = \hat{\delta}/2$, we obtain:
        \begin{equation*}
          \forall t \in
          [0,a]: \quad
          T_{\negthinspace K,n}\left( \Phi(\phi^*)(t) \symExpand
            \hat{\delta}/2 \right) \sqsubseteq
          \Phi(b_{\hat{\delta}/2})(t).
        \end{equation*}
        By symmetrically expanding both sides by $\hat{\delta}/2$,
        and clipping the results using $T_{\negthinspace K,n}$, we
        obtain:
        \begin{equation}
          \label{eq:expanded_clipped}
          \forall t \in
          [0,a]: \quad
          T_{\negthinspace K,n}\left( \Phi(\phi^*)(t) \symExpand
            \hat{\delta} \right) \sqsubseteq T_{\negthinspace K,n} \left(
        \Phi(b_{\hat{\delta}/2})(t) \symExpand
          \hat{\delta}/2 \right).
      \end{equation}
      From~\eqref{eq:expanded_clipped} and~\eqref{eq:d_b_delta_hat},
      we obtain:
      \begin{eqnarray*}
        & & \forall t \in [0,a]: 
        \hat{b}(t) \sqsubseteq T_{\negthinspace K,n} \left(
        \Phi(b_{\hat{\delta}/2})(t) \symExpand
          \hat{\delta}/2 \right)\\
        (\text{by Proposition~\ref{prop:rational_Q_fun_sep_arbit_f}}) &
                                                                                           \implies &
                                                                                                         \hat{b}
                                                                                                         \precAbsBas_{\mathcal{D}}
                                                                                                         \Phi(b_{\hat{\delta}/2})
        \\
        (\text{by~\eqref{eq:def_f_star_X_to_D}}) 
                                                                                         &
                                                                                           \implies
                                                                                                       &
                                                                                                         \hat{b}
                                                                                                         \in
                                                                                                         \left(
                                                                                                         \Phi(b_{\hat{\delta}/2})\right)_*.
      \end{eqnarray*}

    \item Monotonicity of $F$ follows
      from~\ref{item:F_Phi_general_alternative_formulation} and the
      monotonicity of union. As ${\mathcal{W}}_{\mathcal{D}}$ is
      $\omega$-continuous,
      by~\cite[Proposition~2.2.14]{AbramskyJung94-DT}, it suffices to
      prove that for any chain
      $\phi_0 \subseteq \phi_1 \subseteq \ldots \subseteq \phi_k
      \subseteq \ldots$, we have
      $F \left( \bigcup_{k \in \N} \phi_k\right) = \bigcup_{k \in \N}F
      (\phi_k)$. This is also a consequence
      of~\ref{item:F_Phi_general_alternative_formulation}, because:
    \begin{equation*}
    F \left( \bigcup_{k \in \N}
      \phi_k\right) =  \bigcup_{b \in \cup_{k \in \N} \phi_k}
      \left( \Phi(b)\right)_* = \bigcup_{k \in \N} \bigcup_{b \in \phi_k}
      \left( \Phi(b)\right)_* = \bigcup_{k \in \N} F(\phi_k).  
    \end{equation*}

%
  \end{enumerate}
\end{proof}

The relationship between the operators $F$ and $\Phi$ from
Lemma~\ref{lem:F_Phi_general_result}, and the left and right adjoints
$(\cdot)^*$ and $(\cdot)_*$ of the Galois connection from
Theorem~\ref{thm:Galois_connection_X_to_D} is depicted in the
following commutative diagram, in which
\begin{tikzcd}
  {} \arrow[r, hookrightarrow] &  {}
\end{tikzcd}
denotes a monomorphism, and
\begin{tikzcd}
  {} \arrow[r, twoheadrightarrow] &  {}
\end{tikzcd}
denotes an epimorphism:

\begin{equation*}
  \begin{tikzcd}[row sep = large, column sep = large]
    {\mathcal{D}}_{\Q} \arrow[r, hookrightarrow, "(\cdot)_*"]    
    & {\mathcal{W}}_{\mathcal{D}}\\
    {\mathcal{D}}_{\Q} \arrow[u, "\Phi(\cdot)"] &
    {\mathcal{W}}_{\mathcal{D}} \arrow[u, swap, "F(\cdot)"] \arrow[l,
    twoheadrightarrow,"(\cdot)^*" ]
  \end{tikzcd}
\end{equation*}

\section{Second-Order Euler Operator}
\label{susec:2nd_Order_Euler_Operator}

Based on the foundation laid so far, we derive a functional
formulation of the \emph{second-order Euler operator}~$E^2$ that was
first introduced---with an imperative
formulation---in~\parencite{Edalat_Farjudian_Mohammadian_Pattinson:2nd_Order_Euler:2020:Conf}. We
recall the definition of the operator~$E^2$ as appears
in~\parencite[Definition~3.7]{Edalat_Farjudian_Mohammadian_Pattinson:2nd_Order_Euler:2020:Conf}:

\begin{definition}[Second-Order Euler Operator:
  $E^2$]
  \label{def:2nd_order_Euler}

  Let $a \in (0, \frac{K}{M ( 1+ nM')}]$. The second-order Euler
  operator
  $E^2 : {\mathcal{P}}_1 \times {\mathcal{V}}^1 \to [ 0, a] \to
  \intvaldom[{[-K,K]^n}]$ is defined as follows: for a given partition
  $Q \equiv ( q_0, \ldots, q_k)$ of $[ 0, a]$ satisfying $|Q| \leq 1$,
  and a given pair $( u, u') \in {\mathcal{V}}^1$:
  \begin{equation}
    \label{eq:2nd_order_Euler}
    y( x) \defeq \left\{
      \arrayoptions{2ex}{1.3}
      \begin{array}{ll}
        ( 0, \ldots, 0), & \text{if } x = 0,\\
        y( q_i) + \int_{q_i}^x u \left( y(q_i)\right) + ( t - q_i) (
        u' \cdot u) \left( y( q_i) \symExpand \Delta q_i M \right) dt    , & \text{if } q_i < x \leq q_{i+1},
      \end{array}
      \right.
    \end{equation}    
    \noindent
    in which:
    \begin{itemize}
    \item $y \equiv E^2_{( u, u')}(Q)$;
    \item $\Delta q_i \defeq q_{i+1} - q_i$;
    \item $(u' \cdot u)( \cdot)$
    denotes the product of the interval matrix $u'( \cdot)$ with the
    interval vector $u( \cdot )$.
    \end{itemize}

\end{definition}

In~\parencite{Edalat_Farjudian_Mohammadian_Pattinson:2nd_Order_Euler:2020:Conf},
some basic results regarding the operator $E^2$ have been derived. We
will, in particular, refer to the following:

\begin{lemma}
  \label{lem:E2_Lipschitz}
  Assume that $Q = ( q_0, \ldots, q_k) \in {\cal P}_1$, and let
  $y_j = [ \lep[y_j], \uep[y_j]] \defeq E_{(u,u')}^2( Q)_j$ be the $j$-th
  component of $E_{(u,u')}^2( Q)$, for every
  $j \in \set{1, \ldots, n}$. Then, both $\lep[y_j]$ and $\uep[y_j]$ are
  Lipschitz continuous with Lipschitz constant:
  \begin{equation*}
    \Lambda_Q \defeq M ( 1 + |Q| nM').
  \end{equation*}
 \noindent
  In particular:
  \begin{equation}
    \label{eq:E2_well_defined}
    \forall x \in [ 0, a]: \quad  E_{(u,u')}^2( Q)(x) \in \intvaldom[{[-K,K]^n}].
  \end{equation}
\end{lemma}

\begin{proof}
  For completeness, we present the proof as given
  in~\parencite[Lemma~3.9]{Edalat_Farjudian_Mohammadian_Pattinson:2nd_Order_Euler:2020:Conf}. We
  prove the lemma for $\uep[y_j]$. The proof for $\lep[y_j]$ is almost
  identical. Hence, our aim is to show that:
  \begin{equation}
    \label{eq:uepf_Lambda_Q_Lipschitz}
    \forall x, x' \in [ 0, a]: \quad \absn{\uep[y_j](x') - \uep[y_j](x)} \leq
    \Lambda_Q \absn{x'-x}.    
  \end{equation}
  Indeed, it suffices to prove~(\ref{eq:uepf_Lambda_Q_Lipschitz}) for
  the special case of $q_i \leq x \leq x' \leq q_{i+1}$, for some
  $0 \leq i \leq k-1$. Referring to~(\ref{eq:2nd_order_Euler}), note
  that:
  \begin{itemize}
  \item The interval entries in the vector $u$ and matrix $u'$ are
    bounded by $M$ and $M'$, respectively. As $u'$ is an $n \times n$
    matrix, and $u$ is an $n \times 1$ vector, each interval
    component of the vector $u' \cdot u$ is bounded by $nMM'$.
  \item $\forall t \in [q_i, q_{i+1}]: t - q_i \leq |Q|$.
  \end{itemize}
  As such, we obtain:
  \begin{equation*}
    \absn{\uep[y_j](x') - \uep[y_j](x)} \leq \int_x^{x'} \left( M + |Q| nMM'
    \right) dt = M (1 + |Q| nM') (x' - x),
  \end{equation*}
  which proves~(\ref{eq:uepf_Lambda_Q_Lipschitz}). Extending the proof
  to all pairs $x, x' \in [0,a]$ is straightforward. Finally,
  claim~(\ref{eq:E2_well_defined}) also follows from the assumptions
  that $a \leq \frac{K}{M(1+nM')}$ and $|Q| \leq 1$. Thus, the proof
  is complete.
\end{proof}

  To derive a functional formulation of $E^2$, we define the following
  auxiliary operator:
\begin{definition}[Operator: $\Phi$]
  \label{def:Phi}
  Let $a \in (0, \frac{K}{M ( 1+ nM_1)}]$, with $K$, $M$, and $M_1$ as
  in Definition~\ref{def:V1}. For a given partition
  $Q \equiv ( q_0, \ldots, q_k) \in {\mathcal{P}}_{\Q}$, a given pair
  $( u, u') \in {\mathcal{V}}^1$, and $\phi \in {\mathcal{D}}_{\Q}$,
  we define:
  \begin{equation}
    \label{eq:Phi}
    y_\phi( x) \defeq \left\{
      \arrayoptions{2ex}{1.3}
      \begin{array}{ll}
        ( 0, \ldots, 0), & \text{if } x = 0,\\
        T_{\negthinspace K,n} \left[ \phi( q_j) + \int_{q_j}^x u \left( \phi(q_j)\right) + ( t - q_j) (
        u' \cdot u) \left(  T_{\negthinspace K,n} ( G_j (\phi) )
        \right) \md t \right]   , & \text{if } q_j < x \leq q_{j+1},
      \end{array}
      \right.
    \end{equation}    
    \noindent
    in which:
    \begin{itemize}
    \item $T_{\negthinspace K,n}$ is as defined in~\eqref{eq:T_K_n}.

    \item
      $G_j(\phi) \defeq \phi( q_j) \symExpand M \Delta_j$
    and $\Delta_j \defeq q_{j+1} - q_j$.
  \item $(u' \cdot u)( \cdot)$ denotes the product of the interval
    matrix $u'( \cdot)$ with the interval vector $u( \cdot )$.
    \end{itemize}
    The operator
    $\Phi :{\mathcal{P}}_{\Q} \times {\mathcal{V}}^1 \to
    {\mathcal{D}}_{\Q} \to \left([0,a] \to \intvaldom[{[-K,K]^n}]
    \right)$ is defined by:
    \begin{equation*}
      \Phi_{(u,u')}(Q)(\phi) \defeq y_\phi. 
    \end{equation*}
  \end{definition}
  
  \begin{proposition}
    \label{prop:Phi_D_QL_unif_ideal_prop}
    Assume that $a \in (0, \frac{K}{M ( 1+ nM_1)}]$,
    $Q \equiv ( q_0, \ldots, q_k) \in {\mathcal{P}}_{\Q}$, and
    $( u, u') \in {\mathcal{V}}^1$. Then:
    \begin{equation}
      \label{eq:Phi_range_D_QL}
      \forall \phi \in {\mathcal{D}}_{\Q}: \quad
      \Phi_{(u,u')}(Q)(\phi) \in {\mathcal{D}}_{\Q}.
    \end{equation}
    Furthermore,
    $\Phi_{(u,u')}(Q): {\mathcal{D}}_{\Q} \to {\mathcal{D}}_{\Q}$ has
    the \ac{UIC} property.
  \end{proposition}

  \begin{proof}

    In~\eqref{eq:Phi}, let us denote the unclipped term as:
    \begin{equation*}
      \hat{y}_\phi( x) = 
    \phi( q_j) + \int_{q_j}^x u \left( \phi(q_j)\right) + ( t - q_j) (
        u' \cdot u) \left(  T_{\negthinspace K,n} ( G_j (\phi) )
        \right) \md t      , \quad \text{if } q_j < x \leq q_{j+1}.
      \end{equation*}
      We note that $\phi$ is applied only at $q_j$. Hence, the terms
      $\phi( q_j)$, $u \left( \phi(q_j)\right)$, and
      $( u' \cdot u) \left( T_{\negthinspace K,n} ( G_j (\phi) )
      \right)$, are all constant over the interval $(q_j,
      q_{j+1}]$. Let us assume that:
      \begin{equation*}
        \left\{
          \begin{array}{l}
            \phi( q_j) = (\alpha_{j,1}, \ldots, \alpha_{j,n}) \in \intvaldom[\R^n\lift],\\
            u \left( \phi(q_j)\right) = (\beta_{j,1}, \ldots,
            \beta_{j,n}) \in \intvaldom[\R^n\lift],\\
            ( u' \cdot u) \left( T_{\negthinspace K,n} ( G_j (\phi) )
      \right) = ( \gamma_{j,1}, \ldots, \gamma_{j,n}) \in \intvaldom[\R^n\lift].
          \end{array}
          \right.
        \end{equation*}
        As such, for each $1 \leq i \leq n$, we obtain:
        \begin{equation*}
          \left\{
            \begin{array}{l}
              \uep[(\hat{y}_\phi( x))_i] = \uep[\alpha_{j,i}] +
              \uep[\beta_{j,i}] (x - q_j) + \uep[\gamma_{j,i}] (x-q_j)^2/2, \\
              \lep[(\hat{y}_\phi( x))_i] = \lep[\alpha_{j,i}] +
              \lep[\beta_{j,i}] (x - q_j) + \lep[\gamma_{j,i}]
              (x-q_j)^2/2.
            \end{array}
          \right.
        \end{equation*}
        Therefore, the upper and lower bounds of
        $(\hat{y}_\phi( x))_i$ are quadratic---hence continuous---over
        the interval $(q_j, q_{j+1}]$. After applying
        $T_{\negthinspace K,n}$, these bounds remain
        continuous. Hence, for each $1 \leq i \leq n$, the bounds
        $\lep[(y_\phi)_i]$ and $\uep[(y_\phi)_i]$ are continuous over
        each $(q_j,q_{j+1}]$. This, together with the assumption
        $Q \in {\mathcal{P}}_{\Q}$, implies~\eqref{eq:Phi_range_D_QL}.
        
        To prove the \ac{UIC} property, first we note that, as all the
        operations in~\eqref{eq:Phi} are monotonic, so is
        $\Phi_{(u,u')}(Q)$. Furthermore, for any given
        $\psi \in {\mathcal{W}}_{\mathcal{D}}$ and $\epsilon > 0$, for
        each $i \in \set{ 0, \ldots, k}$,
        using~\eqref{eq:def_phi_star_X_to_D}, we find
        $b^i_{\epsilon} \in \psi$ such that
        $T_{\negthinspace K,n} (\psi^*(q_i) \symExpand \epsilon)
        \sqsubseteq b^i_{\epsilon}(q_i)$. As the set
        $\setbarNormal{b^i_{\epsilon}}{i \in \set{0, \ldots, k}}$ is
        finite and $\psi$ is rounded, there exists an element
        $b_{\epsilon} \in \psi$ which satisfies
        $\forall i \in \set{0, \ldots, k}: b^i_{\epsilon}
        \precAbsBas_{\mathcal{D}} b_{\epsilon}$. From
        Proposition~\ref{prop:prec_step_funs_basic_properties}~\ref{item:prop:prec_implies-sqsubseteq},
        we obtain
        $\forall i \in \set{0, \ldots, k}: b^i_{\epsilon} \sqsubseteq
        b_{\epsilon}$. Thus:
      \begin{equation}
        \label{eq:d_b_eps_phi_star}
        \forall \epsilon > 0: \exists b_{\epsilon} \in \psi: \forall i
        \in \set{ 0, \ldots, k}: \quad T_{\negthinspace K,n}
        (\psi^*(q_i) \symExpand \epsilon) 
        \sqsubseteq b_{\epsilon}(q_i).
      \end{equation}
      Over the interval $(q_j, q_{j+1}]$, the dependence of $y_\psi$
      (as defined in~\eqref{eq:Phi}) on $\psi$ is solely based on
      $\psi(q_j)$. As a result, from~\eqref{eq:d_b_eps_phi_star},
      combined with Scott continuity of $u$, $u'$, and integration, we
      deduce:
      \begin{equation*}
        \forall \psi \in {\mathcal{W}}_{\mathcal{D}}: \forall \delta > 0: \exists b_{\delta} \in \psi: \forall t \in
        [0,a]: \quad
        T_{\negthinspace K,n}\left( \Phi_{(u,u')}(Q)(\psi^*)(t) \symExpand
          \delta \right) \sqsubseteq
        \Phi_{(u,u')}(Q)(b_{\delta})(t).
      \end{equation*}
  \end{proof}

  \begin{definition}[Euler Operator: $F$]
  \label{def:F}

  Let $a \in (0, \frac{K}{M ( 1+ nM_1)}]$. The parametric second-order
  Euler operator:
  \begin{equation*}
    F : {\mathcal{P}}_{\Q} \times {\mathcal{V}}^1 \to {\mathcal{W}}_{\mathcal{D}} \to {\mathcal{W}}_{\mathcal{D}}  
\end{equation*}
is defined as follows: for any given
$Q \equiv ( q_0, \ldots, q_k) \in {\mathcal{P}}_{\Q}$,
$( u, u') \in {\mathcal{V}}^1$, and $\phi \in {\mathcal{W}}_{\mathcal{D}}$:
\begin{equation}
  \label{eq:def_F}
  F_{\negthinspace (u,u')}(Q)(\phi) \defeq
  \left( \Phi_{(u,u')}(Q)(\phi^*)\right)_* .
    \end{equation}
  \end{definition}

  The formula~\eqref{eq:def_F} can be depicted as the following
  commutative diagram:

\begin{equation*}
  \begin{tikzcd}[row sep = large, column sep = large]
    {\mathcal{D}}_{\Q} \arrow[r, hookrightarrow, "(\cdot)_*"]    
    & {\mathcal{W}}_{\mathcal{D}}\\
    {\mathcal{D}}_{\Q} \arrow[u, "\Phi_{(u,u')}(Q)(\cdot)"] &
    {\mathcal{W}}_{\mathcal{D}} \arrow[u, swap, "F_{\negthinspace (u,u')}(Q)(\cdot)"] \arrow[l,
    twoheadrightarrow,"(\cdot)^*" ]
  \end{tikzcd}
\end{equation*}
  
\begin{lemma}
  \label{lemma:F_Scott_Continuous}
  For every partition $Q \in {\mathcal{P}}_{\Q}$ and
  $( u, u') \in {\mathcal{V}}^1$, the function
  $F_{\negthinspace (u,u')}(Q): {\mathcal{W}}_{\mathcal{D}} \to {\mathcal{W}}_{\mathcal{D}}$ is Scott
  continuous.
\end{lemma}
\begin{proof}
  This follows from Proposition~\ref{prop:Phi_D_QL_unif_ideal_prop}
  and Lemma~\ref{lem:F_Phi_general_result}.
\end{proof}

\begin{theorem}[Second-Order Euler Operator: $E^2$]
  \label{thm:E_2_MFPS_Ideal_Compl}
  Assume that $a \in (0, \frac{K}{M ( 1+ nM_1)}]$,
  $Q \equiv ( q_0, \ldots, q_k) \in {\mathcal{P}}_{1,\Q}$ is a
  partition of $[0,a]$,
  $( u, u') \in {\mathcal{V}}^1$, and $E^2$ is the
  second-order Euler operator of
  Definition~\ref{def:2nd_order_Euler}. If we denote the bottom
  element of ${\mathcal{W}}_{\mathcal{D}}$ by $\bot$, then we have
  $E^2_{(u,u')}(Q) = \psi^*$ where:
  \begin{equation*}
    \psi \defeq \fix F_{\negthinspace (u,u')}(Q) = \bigsqcup_{m \in \N} \left(
      F_{\negthinspace (u,u')}(Q) \right)^m (\bot) = \left(
      F_{\negthinspace (u,u')}(Q) \right)^{k+1} (\bot).  
  \end{equation*}  
\end{theorem}

\begin{proof}
    For each $j \in \N$, we define
    $\psi_{[j]} \defeq \left( F_{\negthinspace (u,u')}(Q) \right)^j
    (\bot)$ and $Y_{[j]} \defeq y_{\psi^*_{[j]}}$, in which
    $y_{\psi^*_{[j]}}$ is as defined in~\eqref{eq:Phi}. In particular,
    we have,
    $ \forall j \in \N: \psi_{[j+1]} = \left( Y_{[j]}
    \right)_*$. Furthermore, from
    Theorem~\ref{thm:Galois_connection_X_to_D}~\ref{item:Galois_id_X_to_D},
    we obtain:
    \begin{equation}
      \label{eq:psi_star_j_1_Y_j}
      \forall j \in \N: \quad \psi^*_{[j+1]} = \left( \left (Y_{[j]}
        \right)_* \right)^* = Y_{[j]}.
    \end{equation}
    Using induction, we prove that:
  \begin{equation}
    \label{eq:E2_Induction}
    \forall j \in \set{0, \ldots, k}: \forall x \in [0,q_j]: E^2_{( u,
      u')}(Q)(x) = Y_{[j]}(x).
  \end{equation}
  The base case of $j=0$ is immediate as both sides evaluate to
  $( 0, \ldots, 0)$. Next, assume that~\eqref{eq:E2_Induction} holds
  up to some $j \in \set{0, \ldots, k-1}$. In particular, we have
  $E^2_{( u, u')}(Q)(q_j) = Y_{[j]}(q_j)$. We define
  $\alpha_j \defeq E^2_{( u, u')}(Q)(q_j) = Y_{[j]}(q_j)$ and obtain:
  \begin{equation}
    \label{eq:Y_j_1-with_T_k_n}
      \forall x \in (q_j, q_{j+1}]: \quad Y_{[j+1]}( x) =
      y_{\psi^*_{[j+1]}}(x) = T_{\negthinspace K,n} \left[ \alpha_j +
        \int_{q_j}^x u \left( \alpha_j \right) + ( t - q_j) (
        u' \cdot u) \left(  T_{\negthinspace K,n} ( \alpha_j \symExpand M \Delta_j        \right) \md t \right] .
  \end{equation}
  By using~\eqref{eq:E2_well_defined}, we can remove applications of
  $T_{\negthinspace K,n}$ in~\eqref{eq:Y_j_1-with_T_k_n}, and obtain:
    \begin{equation*}
    \forall x \in (q_j, q_{j+1}]: \quad Y_{[j+1]}( x) = \alpha_j +
    \int_{q_j}^x u \left( \alpha_j \right) + ( t - q_j) (
    u' \cdot u) ( \alpha_j \symExpand M \Delta_j)         \md t  ,
  \end{equation*}
  which, combined with~\eqref{eq:2nd_order_Euler}, implies
  $\forall x \in (q_j, q_{j+1}]: Y_{[j+1]}( x) = E^2_{( u,
    u')}(Q)(x)$. Thus, we have
  proven~\eqref{eq:E2_Induction}. From~\eqref{eq:psi_star_j_1_Y_j}
  and~\eqref{eq:E2_Induction}, we obtain:
    \begin{equation*}
    \forall j \in \set{0, \ldots, k}: \forall x \in [0,q_j]: E^2_{( u,
      u')}(Q)(x) = \psi^*_{[j+1]}(x),
  \end{equation*}
  which completes the proof of the theorem.    
  \end{proof}

\begin{remark}
  In the statement of Theorem~\ref{thm:E_2_MFPS_Ideal_Compl}, we
  required $Q \in {\mathcal{P}}_{1,\Q}$ instead of
  $Q \in {\mathcal{P}}_{\Q}$. This is because the proof of
  Lemma~\ref{lem:E2_Lipschitz} requires this stronger assumption.
\end{remark}

\begin{remark}
    To obtain an enclosure of the solution of the
    \ac{IVP}~\eqref{eq:main_ivp}, in
    Theorem~\ref{thm:E_2_MFPS_Ideal_Compl}, one must require
    $( u, u') \in {\mathcal{V}}_{\negmedspace f}^1$, in which $f$ is
    the vector field of \ac{IVP}~\eqref{eq:main_ivp}.
  \end{remark}

\subsection{Computability}
\label{subsec:computability}

As pointed out previously, the Galois connection
of~\eqref{eq:Galois_DQ_WD} allows us to delegate \ac{IVP} solving from
${\mathcal{W}}_{\mathcal{D}}$ to ${\mathcal{D}}_{\Q}$ and then return
to ${\mathcal{W}}_{\mathcal{D}}$, a fact which has also been reflected
in the formulation of the Euler operator $F$ in~\eqref{eq:def_F}. This
raises a valid question on the purpose of constructing, and the
utility of, the $\omega$-continuous domain
${\mathcal{W}}_{\mathcal{D}}$. The answer lies mainly in computable
analysis. As the \ac{POSET} ${\mathcal{D}}_{\Q}$ is not continuous, it
does not admit an effective structure, while the $\omega$-continuous
domain ${\mathcal{W}}_{\mathcal{D}}$ does. In fact, to study
computability of the Euler operator, an effective structure is needed
on all the underlying domains. To that end, we use the concept of
effectively given domains~\parencite{Smyth77-effgd}. Specifically, we
follow the approach taken
in~\cite[Section~3]{Edalat-Sunderhauf-dt-computability-on-reals99}.

\begin{definition}[Effectively given domain]
Assume that $(D, \sqsubseteq)$ is
an $\omega$-continuous domain, with a countable basis $B$ that is
enumerated as follows:
    \begin{equation}
      \label{eq:basis_enum}
      B = \set{ b_0 = \bot, b_1, \ldots, b_n, \ldots}.
    \end{equation}
    We say that the domain $D$ is \emph{effectively given}, with
    respect to the enumeration~\eqref{eq:basis_enum}, if the set
    $\setbarNormal{(i,j) \in \N \times \N}{ b_i \ll b_j}$ is
    recursively enumerable, in which $\ll$ is the way-below relation
    on $D$.
  \end{definition}

    In the following proposition, computability is to be understood
    according to \acl{TTE}~\parencite{Weihrauch2000:book}.
    
    \begin{proposition}[Computable elements and functions]
      \label{prop:comp_elem_fun}
    Let $D$ and $E$ be two effectively given domains, with enumerated
    bases $B_1 = \set{d_0, d_1, \ldots, d_n, \ldots}$ and
    $B_2 = \set{e_0, e_1, \ldots, e_n, \ldots}$, respectively:

    \begin{enumerate}[label=(\roman*)]
      \item An
    element $x \in D$ is computable $\iff$ the set
    $\setbarNormal{i \in \N}{d_i \ll x}$ is recursively enumerable.

    \item \label{item:egd_fun_computable}
    A
    map $f: D \to E$ is computable $\iff$
    $\setbarNormal{(i,j) \in \N \times \N}{ e_i \ll f(d_j)}$ is
    recursively enumerable.

        \end{enumerate}
  \end{proposition}

  \begin{proof}
    See ~\parencite[Proposition~3 and Theorem~9]{Edalat-Sunderhauf-dt-computability-on-reals99}.
  \end{proof}
  
  The domain ${\mathcal{W}}_{\mathcal{D}}$ is $\omega$-continuous, with
  $i_{\mathcal{D}}({\mathcal{B}}_{\mathcal{D}})$ as a countable basis. We have:
  \begin{equation}
    \label{eq:waybelow_lessthan_B_Q}
    \forall b,b' \in {\mathcal{B}}_{\mathcal{D}}: \quad b \precAbsBas b' \iff i_{\mathcal{D}}(b)
    \ll i_{\mathcal{D}}(b').
  \end{equation}

\begin{proposition}
  Assume that ${\mathcal{B}}_{\mathcal{D}}$ is enumerated as
  ${\mathcal{B}}_{\mathcal{D}} = \set{b_0, b_1, b_2, \ldots }$, and
  $i_{\mathcal{D}}({\mathcal{B}}_{\mathcal{D}})$ is enumerated as:  
\begin{equation}
  \label{eq:i_B_Q:enumerated}
  i_{\mathcal{D}}({\mathcal{B}}_{\mathcal{D}}) = \set{i_{\mathcal{D}}(b_0), i_{\mathcal{D}}(b_1), i_{\mathcal{D}}(b_2), \ldots }.
\end{equation}
Then, the way-below relation $\ll$ over $i_{\mathcal{D}}({\mathcal{B}}_{\mathcal{D}})$ is
recursive with respect to the enumeration~\eqref{eq:i_B_Q:enumerated}.
\end{proposition}
\begin{proof}
  This follows from~\eqref{eq:waybelow_lessthan_B_Q}, and the fact
  that only rational numbers are used in representation of the
  elements of ${\mathcal{B}}_{\mathcal{D}}$.
\end{proof}

Thus, we obtain an effective structure over
${\mathcal{W}}_{\mathcal{D}}$, with respect to the
enumeration~\eqref{eq:i_B_Q:enumerated}. The Euler operator takes
input from the domain ${\mathcal{V}}^1$ as well. The set of rational
hyper-rectangles $B_{\intvaldom[\R^n\lift]}$ is a basis for
$\intvaldom[\R^n\lift]$. It is straightforward to construct an
effective structure over $\intvaldom[\R^n\lift]$ with respect to any
reasonable enumeration of $B_{\intvaldom[\R^n\lift]}$. Moreover, using
rational hyper-rectangles, we can construct an effective structure
over ${\mathcal{V}}^1$ with respect to an enumeration of a basis:
\begin{equation}
  \label{eq:BV1_enum}
B_{{\mathcal{V}}^1} = \set{(u_0, u'_0), (u_1, u'_1), \ldots, (u_n, u'_n), \ldots}.  
\end{equation}

\begin{proposition}
  \label{prop:b_BQ_ideal}
  Assume that $Q \in {\mathcal{P}}_{\Q}$ and $( u, u') \in {\mathcal{V}}^1$. Then, for any $b \in {\mathcal{B}}_{\mathcal{D}}$, we have:
  \begin{equation}
    \label{eq:F_Phi_b}
    F_{\negthinspace (u,u')}(Q)(i_{\mathcal{D}}(b)) = \left(
      \Phi_{(u,u')}(Q)(b)\right)_* .
  \end{equation}
\end{proposition}

\begin{proof}
  By Theorem~\ref{thm:Galois_connection_X_to_D}~\ref{item:Galois_id_X_to_D}, we
  have $(i_{\mathcal{D}}(b))^* = (b_*)^* = b$. Hence, the
  equality~\eqref{eq:F_Phi_b} follows from~\eqref{eq:def_F}.
\end{proof}

\begin{corollary}
  \label{cor:F_recursive_relation}
  Assume that $Q \in {\mathcal{P}}_{\Q}$. Then, the relation:
  \begin{equation*}
    i_{\mathcal{D}}(b_i) \ll F_{\negthinspace (u_k,u'_k)}(Q)(i_{\mathcal{D}}(b_j))  
  \end{equation*}
  is recursive with respect to the
  enumerations~\eqref{eq:i_B_Q:enumerated} and~\eqref{eq:BV1_enum}.
\end{corollary}

\begin{proof}
  As $Q \in {\mathcal{P}}_{\Q}$,
  $( u_k, u'_k) \in B_{{\mathcal{V}}^1}$, and
  $b_j \in {\mathcal{B}}_{\mathcal{D}}$, then
  $\Phi_{(u_k,u'_k)}(Q)(b_j)$ is a piecewise quadratic enclosure with
  rational coefficients. By Proposition~\ref{prop:b_BQ_ideal},
  deciding
  $i_{\mathcal{D}}(b_i) \ll F_{\negthinspace
    (u_k,u'_k)}(Q)(i_{\mathcal{D}}(b_j))$ reduces to deciding
  $b_i \precAbsBas \Phi_{(u_k,u'_k)}(Q)(b_j)$, which, in turn, reduces
  to deciding inequalities of the form
  $\alpha x^2 + \beta x + \gamma < \delta$ over an interval $[p,q]$,
  with $\alpha, \beta, \gamma, \delta \in \Q$, and with $p$ and $q$
  computable, which is semi-decidable.
\end{proof}

From Proposition~\ref{prop:comp_elem_fun}~\ref{item:egd_fun_computable}
and Corollary~\ref{cor:F_recursive_relation}, we obtain the main
result of this section:

\begin{theorem}[Computability]
  \label{thm:computability_F}
  For any fixed $Q \in {\mathcal{P}}_{\Q}$, the map
  $F_{\negthinspace ( \cdot , \cdot)}(Q)(\cdot): {\mathcal{V}}^1 \times
  {\mathcal{W}}_{\mathcal{D}} \to {\mathcal{W}}_{\mathcal{D}}$ is computable.
\end{theorem}

  \subsection{Convergence Analysis}
  \label{subsec:convergence_analysis}
  
In this section, we demonstrate that the operator $E^2$ is indeed
second-order. In classical numerical analysis, one may find the
following criteria that determine whether a method is of order $p$
(see, {\eg},~\cite[Page 8]{Iserles:2009}):

\begin{enumerate}[label=(C\arabic*)]
\item \label{item:n_th_order_criterion_polynomial} A numerical scheme
  for solving \acp{IVP} is said to be of order $p$ if, whenever the
  solution of an \ac{IVP} is a polynomial of degree at most $p$, the
  solution can be obtained exactly by the scheme.

\item \label{item:n_th_order_criterion_error_order} Alternatively, a
  numerical scheme is said to be of order $p$ if, for every step size
  $h$, the local error incurred is of order $O(h^{p+1})$. For the
  Euler method, this entails that the global error must be of order
  $O(h^p)$.
\end{enumerate}

The bounds obtained for the width of $E^2_{(u,u')}(Q)$
in~\cite[Section~4]{Edalat_Farjudian_Mohammadian_Pattinson:2nd_Order_Euler:2020:Conf}
are too conservative and do not reflect the second-order nature of the
method. For instance,
in~\cite[Corollary~4.3]{Edalat_Farjudian_Mohammadian_Pattinson:2nd_Order_Euler:2020:Conf},
the following bound is obtained for equidistant partitions:
\begin{equation}
  \label{eq:conv_speed_ENTCS_equidistant_Q}
  \widthOf\left(E^2_{(u,u')}(Q)\right)
  \leq \frac{1}{2} |Q|M \left( \me^{aL} - 1 \right).
\end{equation}
This is of the same order as the bounds obtained
in~\parencite{EdalatPattinson2006-Euler-PARA} for the first-order
Euler method. Here, we present a more accurate convergence analysis
which demonstrates that $E^2$ is indeed second-order, according to
criterion~\ref{item:n_th_order_criterion_error_order}. The basic idea
is to adopt the midpoint-width representation of intervals introduced
by~\textcite{Moore_Jones:Safe_Starting_Regions:1977}.

\begin{definition}[$m(A)$] For every interval
  $I = [a,b]$, we define $m(I) \defeq (a+b)/2$. For every interval
  matrix $A = [A_{ij}]_{m \times n}$, we define
  $m(A) \defeq [m(A_{ij})]_{m \times n}$.
\end{definition}
\begin{proposition}[Midpoint-width representation]
  Assume that $A = [A_{ij}]_{m \times n}$ is an interval matrix. Then:
  \begin{equation*}
    A = m(A) + W,
  \end{equation*}
  in which $W$ is an $m \times n$ interval matrix with entries $W_{ij}
  = \frac{1}{2} [-1,1] \widthOf(A_{ij})$.
\end{proposition}
\begin{proof}
This is a straightforward generalization of the fact that any interval
$I = [a,b]$ may be written as:
\begin{equation*}
I = [a,b] = \frac{a+b}{2} + \left[ \frac{a-b}{2}, \frac{b-a}{2} \right]
  = m(I) + \frac{1}{2} [-1,1] \widthOf(I).
\end{equation*}
\end{proof}
  
To simplify the arguments that follow, we reiterate an assumption
from~\parencite{Edalat_Farjudian_Mohammadian_Pattinson:2nd_Order_Euler:2020:Conf}:

\begin{assumption}
  In the sequel, we assume that $u$ and $u'$ are interval extensions
  of the vector field $f$ and its $\overline{L}$-derivative $\overline{L}(f)$,
  respectively, and satisfy the following additional condition:
\begin{equation}
  \label{eq:L1_Lipschitz_Condition}
\forall x \in \intvaldom[{[-K,K]^n}]: \quad  \widthOf(u(x)) \leq
\norm{u'(x)}_\infty  \widthOf(x),
\end{equation}
where $\norm{\cdot}_{\infty}$ is the interval matrix norm of
Definition~\ref{def:interval_norms}.
\end{assumption}
\noindent
As stated
in~\cite[Corollary~3.6]{Edalat_Farjudian_Mohammadian_Pattinson:2nd_Order_Euler:2020:Conf},
if $u$ and $u'$ are the canonical interval extensions of the classical
vector field $f: [-K,K]^n \to [-M,M]^n$ and its
$\overline{L}$-derivative $\overline{L}(f)$, respectively, then, $u$
and $u'$ do satisfy condition~(\ref{eq:L1_Lipschitz_Condition}).

Theorem~\ref{thm:E_2_MFPS_Ideal_Compl} entails that, if
$Q \equiv ( q_0, \ldots, q_k) \in {\mathcal{P}}_{1,\Q}$, and if we let
$y \defeq E^2_{( u, u')}(Q)$, then:
  \begin{equation}
    \label{eq:2nd_order_Euler_E}
    y( x) = \left\{
      \arrayoptions{2ex}{1.3}
      \begin{array}{ll}
        ( 0, \ldots, 0), & \text{if } x = 0,\\
        y( q_i) + \int_{q_i}^x u \left( y(q_i)\right) + ( t - q_i) (
        u' \cdot u) \left( y( q_i) \symExpand M \Delta_i \right) \md
        t    , & \text{if } q_i < x \leq q_{i+1}.
      \end{array}
      \right.
    \end{equation}    
    For each $i \in \set{0, \ldots, k-1}$, let us define:
    \begin{equation}
      \label{eq:def_A_i}
      A_i \defeq y( q_i) \symExpand M \Delta_i .
    \end{equation}
    By the midpoint-width representation, we write:
    \begin{equation}
      \label{eq:midpoint-width_u_A_i}
      \left\{
        \arrayoptions{0.5ex}{1.3}
        \begin{array}{lcl}
          u(A_i)& = & m(u(A_i)) + W_i,\\
          u'(A_i) & = & m(u'(A_i)) + W'_i.
        \end{array}
      \right.
    \end{equation}
    \noindent
    Note that $u(A_i)$, $m(u(A_i))$, and $W_i$ are $n \times 1$
    vectors, whereas $u'(A_i)$, $m(u'(A_i))$, and $W'_i$ are $n \times
    n$ matrices. Next, we define the following constants:
    \begin{equation}
      \label{eq:def_omega_omega_prime}
      \omega_Q \defeq \frac{1}{2} \max_{0 \leq i \leq k-1} \widthOf(W_i),\quad
      \omega'_Q \defeq \frac{1}{2} \max_{0 \leq i \leq k-1} \widthOf(W'_i).
    \end{equation}
    When the partition $Q$ is clear from the context, we drop the
    subscripts and simply write $\omega$ and $\omega'$.

\begin{proposition}
  \label{prop:speed_of_conv_2nd_order}
  Assume that $Q \equiv ( q_0, \ldots, q_k) \in {\mathcal{P}}_{1,\Q}$,
  $L \defeq n M_1$, $\omega$ and $\omega'$ are as
  in~\eqref{eq:def_omega_omega_prime}, and define:
  \begin{equation}
    \label{eq:def_rho_Q}
    \rho_Q \defeq L \omega + nM  \omega' + n \omega \omega' .
  \end{equation}
  Then, $\forall i \in \set{0, \ldots, k-1}: \forall x \in [ q_i, q_{i+1}]: \widthOf\left( E^2_{(u,u')}(
      Q) (x) \right) \leq  \widthOf\left( E^2_{(u,u')}(
      Q) (q_i) \right)( 1 + |Q|L) + \frac{|Q|^2}{2} \rho_Q$.
\end{proposition}
  \begin{proof}
    Let $y \defeq E^2_{(u,u')}(Q)$. From~\eqref{eq:2nd_order_Euler_E},
    we obtain:
    \arrayoptions{0.5ex}{1.3}
    \begin{eqnarray*}
      \widthOf(y(x)) & \leq & \widthOf(y( q_i)) + \int_{q_i}^x \widthOf
                              \left( u \left( y(q_i)\right) \right) +
                              \widthOf \left( ( t - q_i) (
                              u' \cdot u) \left( y( q_i) \symExpand M \Delta_i \right) \right) \md t \\
      \text{(by~(\ref{eq:L1_Lipschitz_Condition}) and \eqref{eq:def_A_i})} & \leq &                                                                      \widthOf(y(q_i))
                                                                                    +
                                                                                    \int_{q_i}^x
                                                                                    \norm{u'(y(q_i))}_\infty
                                                                                    \widthOf(y(q_i))
                                                                                    \md
                                                                                    t
                                                                                    +
                                                                                    \int_{q_i}^x \widthOf \left( ( t - q_i) (
                                                                                    u'
                                                                                    \cdot
                                                                                    u)
                                                                                    \left(
                                                                                    A_i
                                                                                    \right)
                                                                                    \right)
                                                                                    \md
                                                                                    t\\
                     & \leq &  \widthOf\left( y( q_i) \right) +
                              \widthOf\left( y( q_i) \right) (x-q_i) L
                              + \frac{1}{2}(x-q_i)^2 \widthOf \left( u'(A_i)
                              u(A_i) \right)\\
      (\text{as $x-q_i \leq |Q|$})
& \leq &  \widthOf\left( y( q_i) \right) (1 + |Q| L)
                              + \frac{|Q|^2}{2}\widthOf\left(u'(A_i)
                              u(A_i) \right).
    \end{eqnarray*}
    It remains to show that $\widthOf(u'(A_i) u(A_i)) \leq
    \rho_Q$. From the midpoint-width
    representation~\eqref{eq:midpoint-width_u_A_i}, we obtain:
    \begin{align*}
      u'(A_i) u(A_i) &= \Big( m(u'(A_i)) + W'_i \Big) \Big(
                       m(u(A_i)) + W_i \Big)\\
 &= m(u'(A_i)) m(u(A_i)) + m(u'(A_i)) W_i + W'_i m(u(A_i)) + W'_i W_i.
    \end{align*}
    A term-by-term analysis of the width of the last expression shows that:
    \begin{itemize}
    \item The first term, {\ie}, $m(u'(A_i)) m(u(A_i))$, has width
      zero.
    \item The width of the second term $m(u'(A_i)) W_i$ is bounded by
      the product of
      $\norm{m(u'(A_i))}_\infty$ and $\omega$. Hence, it is bounded by
      $L \omega$.
    \item Similarly, the width of the third term $W'_i m(u(A_i))$ is
      bounded by the product of $\omega'$ and
      $\norm{m(u(A_i))}_1$. Hence, it is bounded by $\omega' nM$.
    \item The width of the fourth term $W'_i W_i$ is also clearly
      bounded by $n\omega' \omega$.
    \end{itemize}
    Thus, we obtain $\widthOf(u'(A_i) u(A_i)) \leq
    \rho_Q$.    
\end{proof}

\begin{corollary}[Speed of convergence]
  \label{cor:speed_convergence}
  Assume that $\rho_Q$ is as defined in~\eqref{eq:def_rho_Q} and
  $L > 0$. For any partition $Q \in {\mathcal{P}}_{1,\Q}$, we have:
  \begin{equation}
    \label{eq:conv_speed_arbit_Q}
    \widthOf\left(E^2_{(u,u')}(Q)\right)
    \leq \frac{|Q|\ \rho_Q}{2L}  \left( \me^{a r_Q L} - 1 \right).     
  \end{equation}
\noindent
In particular,
  when $Q$ is equidistant:
  \begin{equation}
    \label{eq:conv_speed_equidistant_Q}
    \widthOf\left(E^2_{(u,u')}(Q)\right)
    \leq \frac{|Q|\ \rho_Q}{2L}  \left( \me^{a L} - 1 \right).     
  \end{equation}
\end{corollary}

\begin{proof}
  Assume that $Q = ( q_0, \ldots, q_k)$, and let $ c \defeq
  \frac{|Q|^2\ \rho_Q}{2}$, $d = 1 + |Q|L$. We prove by induction on $i$
  that:

  \begin{equation*}
    \forall x \in [q_i, q_{i+1}]: \quad \widthOf\left(
      E^2_{(u,u')}(Q)(x) \right)
    \leq c \sum_{j=0}^i d^j.
  \end{equation*}
  The case of $i=0$ is immediate from
  Proposition~\ref{prop:speed_of_conv_2nd_order} and the fact that
  $\widthOf(y(q_0)) = 0$. For $i>0$, again, by
  Proposition~\ref{prop:speed_of_conv_2nd_order}, we have: {
    \arrayoptions{0.5ex}{1.3}
  \begin{eqnarray*}
    \widthOf\left( E^2_{(u,u')}(Q)(x) \right) & \leq &  \widthOf\left(
                                                       E^2_{(u,u')}(Q)(q_i)
                                                       \right) \cdot d
                                                       + c\\
    (\text{by induction hypothesis}) & \leq & \left( c \sum_{j=0}^{i-1} d^j \right)
                                              d + c   =  c \sum_{j=0}^{i} d^j.
  \end{eqnarray*}}
Thus, we obtain:
{  \arrayoptions{0.5ex}{1.3}
  \begin{eqnarray*}
    \widthOf\left( E^2_{(u,u')}(Q) \right) & \leq & c \sum_{j=0}^{k-1} d^j
                                                    = c \frac{d^k -
                                                    1}{d - 1}
    \\
                                           & = & \frac{|Q|\ \rho_Q}{2L}  \left(  (1 + |Q| L)^k - 1\right)\\
                                           & = & \frac{|Q|\ \rho_Q}{2L}
                                                    \left(  (1 + m(Q) r_Q  L)^k - 1\right)\\
    & \leq & \frac{|Q|\ \rho_Q}{2L} \left(  \left(1 + \frac{a}{k} r_Q  L
             \right)^k - 1\right)\\
                                           & \leq &  \frac{|Q|\ \rho_Q}{2L}
                                                    \left( \me^{ar_Q L} - 1 \right),                                                    
  \end{eqnarray*}}
\noindent
which
proves~(\ref{eq:conv_speed_arbit_Q}). Inequality~(\ref{eq:conv_speed_equidistant_Q})
now follows, because for an equidistant $Q$, we have $r_Q = 1$.
\end{proof} 

According to Corollary~\ref{cor:speed_convergence}, if $\rho_Q$ is of
order $O(|Q|)$, then the convergence must be second-order, {\ie},
$O(|Q|^2)$. This can be guaranteed if the vector field is continuously
differentiable:
\begin{theorem}[Second-order convergence]
  \label{thm:Second_Order_Convergence_int_Lip}
  Assume that $u$ and $u'$ are interval extensions of the classical
  vector field $f$ and its $\overline{L}$-derivative $\overline{L}(f)$,
  respectively. If $u'$ is interval Lipschitz, then $E^2$ has
  second-order convergence.
\end{theorem}

\begin{proof}
  By~\cite[Corollary~4.3]{Edalat_Farjudian_Mohammadian_Pattinson:2nd_Order_Euler:2020:Conf},
  we know that $\widthOf\left(E^2_{(u,u')}(Q)\right)$ is
  $O(|Q|)$. Hence, by~\eqref{eq:def_A_i}, $\widthOf(A_i)$ is
  $O(|Q|)$. As $u'$ is bounded, then $u$ is interval
  Lipschitz. Therefore, $\widthOf(u(A_i))$ is $O(|Q|)$, and
  by~\eqref{eq:midpoint-width_u_A_i}
  and~\eqref{eq:def_omega_omega_prime}, we must have
  $\omega \in O(|Q|)$. If we also assume that $u'$ is interval
  Lipschitz, then by a similar argument we can conclude that
  $\omega' \in O(|Q|)$.

  Hence, from~\eqref{eq:def_rho_Q}, and the fact that both $\omega$
  and $\omega'$ are $O(|Q|)$, we deduce that $\rho_Q \in
  O(|Q|)$. This, together with Corollary~\ref{cor:speed_convergence},
  implies that $\widthOf\left(E^2_{(u,u')}(Q)\right)$ is indeed
  $O(|Q|^2)$.
\end{proof}

\subsection{Further Extensions}

In Theorem~\ref{thm:Second_Order_Convergence_int_Lip}, the assumption
that $u'$ is interval Lipschitz was imposed. By
Remark~\ref{rem:intval_Lipschitz_total}, an interval Lipschitz $u'$
must map maximal elements to maximal elements. Since $u'$ is Scott
continuous, assuming that $u'$ is interval Lipschitz entails that $u$
is continuously differentiable. This is adequate for a qualitative
convergence analysis, but in practice, these assumptions are
restrictive. Below, we discuss two relaxations of these assumptions
that are important in practical applications.

\subsubsection{Non-Differentiable Vector Fields}
\label{subsubsec:non_differentiable_fields}

Assuming that $u'$ is interval Lipschitz ensures that
$\omega' \in O(|Q|)$, and as a consequence, that the term $nM \omega'$
in~\eqref{eq:def_rho_Q} also is in $O(|Q|)$. By inspecting the proof
of Proposition~\ref{prop:speed_of_conv_2nd_order}, it can be seen that
if $u'$ maps some maximal elements to non-maximal elements, it may
hinder the $O(|Q|^2)$ convergence of the operators, but only over the
points where it takes non-maximal values.

As we will see in Figure~\ref{fig:id_versus_absSin}, our experiments
show that the second-order convergence is retained for a
non-differentiable vector field. We conjecture that this is true in case $u'$
takes non-maximal values on isolated points of the domain ({\eg}, over
a finite subset of the domain). In fact, by Rademacher's theorem, if a
function is Lipschitz continuous over an open subset $U$ of $\R^n$,
then it is (Fr{\'e}chet) differentiable almost everywhere (with
respect to the Lebesgue measure) over
$U$~\cite[Corollary~4.19]{Clarke_et_al:Nonsmooth_Control:Book:1998}. As
such, it is plausible that the second-order convergence is retained,
even when the vector field is merely Lipschitz continuous.

\subsubsection{Approximations of the Vector Field}

The assumption that $u'$ is interval Lipschitz is restrictive, even
when the vector field is continuously differentiable. In
implementations, finitely-representable objects must be used to
approximate ideal elements. For instance, one may use piecewise
constant enclosures with rational coordinates to approximate
continuous functions in $C([0,1])$. For a function such as
$f(x) = \exp(x): [0,1] \to \R$, it is impossible to find such a
finitely-representable enclosure which has a width of zero over the
entire domain, because $\exp$ is a transcendental function.

To cope with this issue, we consider sequences $(u_n,u'_n)_{n \in \N}$
satisfying:
\begin{equation*}
\lim_{n \to \infty} d(u_n, u) = 0 \quad  \text{and} \quad \lim_{n \to \infty} d(u'_n, u') = 0,
\end{equation*}
in which $d(\cdot, \cdot)$ is the interval distance of
Definition~\ref{def:interval_distance}. With this added layer of
approximation, we can still obtain a second-order convergence, as long
as the sequence $(u_n,u'_n)_{n \in \N}$ converges to $(u,u')$---under
the distance $d(\cdot, \cdot)$---fast enough. Similar analyses may be
found in
~\parencite{EdalatPattinson2006-Euler-PARA,FarjudianKonecny2008:wollic-lnai,Edalat_Farjudian_Mohammadian_Pattinson:2nd_Order_Euler:2020:Conf}. Hence,
we do not present the details here and just state the main theorem.

\begin{theorem}
  \label{thm:speed_of_convergence_field_approx}
  Let $(Q_n)_{n \in \N}$ be a sequence of partitions in
  ${\mathcal{P}}_{1,\Q}$ satisfying $\lim_{n \to \infty} |Q_n| = 0$
  and assume that $u'$ is interval Lipschitz. Furthermore, assume that
  $( u, u') = \bigsqcup_{n \in \N} ( u_n, u'_n)$, and for some
  constant $C_1 > 0$, we have $d( u_n, u) \leq C_1 |Q_n|^2$ and
  $d( u'_n, u') \leq C_1 |Q_n|^2$. By letting
  $y_n \defeq E^2_{( u_n, u'_n)}$, we obtain:
  \begin{enumerate}[label=(\roman*)]
  \item
    $\exists C_2 \geq 0, \forall n \in \N: \quad \widthOf( y_n) \leq
    C_2 |Q_n|^2$.
  \item $\bigsqcup_{n \in \N} y_n$ is real-valued and is a solution of~\eqref{eq:main_ivp}.
  \end{enumerate}
  
\end{theorem}

\begin{proof}
  The proof is a straightforward modification of the proof
  of~\cite[Theorem~4.11]{Edalat_Farjudian_Mohammadian_Pattinson:2nd_Order_Euler:2020:Conf},
  using the midpoint-width analysis.
\end{proof}

Our implementation of the Euler operators is based on the MPFI
library~\parencite{Revol_Rouillier:MPFI:05}, which uses arbitrary precision
dyadic numbers as end-points of intervals. As we will see in
Section~\ref{sec:Experiments}, under this representation, we still
obtain quadratic convergence, even for some non-differentiable vector
fields.

\subsection{Monotonicity with Respect to Refinement of Partitions}

\label{subsec:monotonicity_refinement_partitions}

A desirable property for an Euler operator is monotonicity with
respect to refinement of partitions. If $Q$ and $Q'$ are two
partitions of $[0,a]$ satisfying $Q \sqsubseteq Q'$, it would be
desirable to have
$\forall (u,u') \in {\mathcal{V}}_{\negmedspace f}^1: E_{(u,u')}(Q)
\sqsubseteq E_{(u,u')}(Q')$. The operator $E^2$ does not satisfy this
property, even though we have proven its convergence. In other words,
if
$Q_1 \sqsubseteq Q_2 \sqsubseteq \ldots \sqsubseteq Q_n \sqsubseteq
\ldots$ is a sequence of partitions which satisfies
$\lim_{n \to \infty} |Q_n| = 0$, then it is guaranteed that the
enclosures produced by $E^2_{(u,u')}(Q_n)$ converge to a limit, but
the sequence is not a shrinking chain.

The main reason for this lack of monotonicity is that, although the
Lipschitz constant of the solution $z$ of the main
\ac{IVP}~\eqref{eq:main_ivp} is bounded by $M$, the Lipschitz constant
of the approximations obtained via application of $E^2$---or the
operator $F$ for that matter---are bounded by a larger constant
$M(1+nM')$ (see Lemma~\ref{lem:E2_Lipschitz}). Thus, to obtain
monotonicity with respect to refinement of partitions, in
Definition~\ref{def:Phi}, we can make one of the following
modifications:

\begin{enumerate}[label=(\arabic*)]
\item We can modify $G_j$ by replacing $M$ with $M(1+nM')$ and define
  $G_j(\phi) \defeq \phi( q_j) \symExpand M(1+nM') \Delta_j$. The
  resulting operator will have a slower convergence.
  
\item We can modify the
  definition~\eqref{eq:Phi} as follows:
  \begin{equation*}
    y_\phi( x) \defeq \left\{
      \arrayoptions{0.2ex}{1.3}
      \begin{array}{ll}
        ( 0, \ldots, 0), & \text{if } x = 0,\\
        T_{\negthinspace M(x-q_j),n} \left[ \phi( q_j) + \int_{q_j}^x u \left( \phi(q_j)\right) + ( t - q_j) (
        u' \cdot u) \left(  T_{\negthinspace K,n} ( G_j (\phi) )
        \right) \md t \right]   , 
                         & \text{if } q_j < x \leq q_{j+1}.
      \end{array}
      \right.
    \end{equation*}    
    This modification results in an operator with faster convergence,
    at the cost of a significant increase in clutter for presentation
    purposes.
\end{enumerate}
Due to these considerations, we opted for keeping the formulation
presented in Definition~\ref{def:Phi}.

\section{Runge-Kutta Euler Operator}
\label{sec:Runge_Kutta_Euler_Operator}

In this section, we demonstrate how the framework that we have
developed can be used for domain-theoretic formulation of Runge-Kutta
methods. We begin by a brief reminder of the Runge-Kutta theory,
tailored to the autonomous \ac{IVP}~\eqref{eq:main_ivp}. More
comprehensive accounts may be found in classical textbooks on
numerical analysis,
{\eg},~\parencite{Lambert:Numerical_Methods:Book:1991,Quarteroni_et_al:Numerical_Mathematics:Book:2000}. The
explicit $m$-stage Runge-Kutta method for solving the
\ac{IVP}~\eqref{eq:main_ivp} proceeds as follows:

\begin{enumerate}[label=(\roman*)]
\item The interval $[0,a]$ is partitioned as $Q = (q_0, \ldots, q_k)$,
  for some $k \geq 1$.

\item The initial value is assigned as $y_0 = (0, \ldots, 0)$.

\item For each $j \in \set{0, \ldots, k-1}$, we let $h \defeq q_{j+1} -
  q_j$, and then: 
  \begin{equation}
    \label{eq:RK_steps_general}
    y_{j+1} \defeq y_j + h \sum_{i=1}^m w_i k_{ij},
  \end{equation}
  in which, for every $i \in \set{1, \ldots, m}$:
  \begin{equation}
    \label{eq:kij_RK}
    k_{ij} \defeq f \left( y_j + h \sum_{\ell = 1}^{i-1}
      a_{i \ell} k_{\ell j} \right).
  \end{equation}
\end{enumerate}

The coefficients $w_i$ and $a_{i \ell}$ are parameters that are
specific to each variant of the Runge-Kutta method. For any
Runge-Kutta method of order $p$, the local truncation error at step
$j+1$ is expressed as follows:
\begin{align}
  r_{j+1}(h) & =  y( q_j + h) - \left( y( q_j) + h \sum_{i=1}^m w_i
                   k_{ij}(h) \right)  \nonumber \\
  & =  \psi( y(q_j)) h^{p+1} +  O(h^{p+2}), \label{eq:trunc_err_general}
\end{align}
in which, $y( q_j + h)$ and $y( q_j)$ denote the exact solutions at
$q_j+h$ and $q_j$, respectively, and $k_{ij}(h)$ is obtained
from~\eqref{eq:kij_RK} for the exact value of $y(q_j)$ substituted for
$y_j$.

In developing validated Runge-Kutta methods, one of the main tasks
involves deriving explicit bounds for the truncation error. An
extensive catalogue of validated Runge-Kutta methods may be found
in~\parencite{Marciniak:Selected_Interval_Methods:2009}, together with
explicit formulation of their truncation error. There is no uniform
formulation of the truncation error that is parametrized by (say) the
order of the method. In simple terms, for each variant, the explicit
formulation must be obtained separately. As may be seen
from~\parencite{Marciniak:Selected_Interval_Methods:2009}, the explicit
formulation of the truncation error, especially of higher-order
variants, can become very lengthy and even span pages.

As a result, we consider one of the Runge-Kutta variants---commonly
referred to as the Euler method---which results in relatively shorter
formulas. To simplify the presentation further, we focus on the
\emph{autonomous scalar} case. In contrast with the Euler method of
Section~\ref{susec:2nd_Order_Euler_Operator}, here we follow the
Runge-Kutta theory, and demonstrate that our domain-theoretic
framework is applicable for temporal discretization in the context of
Runge-Kutta theory as well.

\subsection{Scalar Euler: Runge-Kutta Formulation}
\label{subsec:scalar_Euler_Runge_Kutta}

We consider the scalar case of \ac{IVP}~\eqref{eq:main_ivp}, {\ie},
with $n=1$. To begin with, we assume that the vector field $f$ is
continuously differentiable of any order that appears in the upcoming
formulas. For the scalar case of a Runge-Kutta method of order $p$,
Taylor's theorem is applicable, and the truncation error
of~\eqref{eq:trunc_err_general} is commonly expressed as:
\begin{align*}
  r_{j+1}(h) &= \psi( y(q_j)) h^{p+1} +  O(h^{p+2})\\
  &= r^{(p+1)}_{j+1}(0) \frac{h^{p+1}}{(p+1)!} + r^{(p+2)}_{j+1}(
    \theta h) \frac{h^{p+2}}{(p+2)!},
\end{align*}
for some $\theta \in (0,1)$. This follows from the fact that, as the
method is assumed to be of order $p$, we must have
$\forall \ell \in \set{ 0, \ldots, p}: r^{(\ell)}_{j+1}(0) = 0$.

For the autonomous scalar case of the Euler method, the formulation
of~\eqref{eq:RK_steps_general} can be presented as follows:
\begin{equation*}
    y_{j+1} \defeq y_j + f(y_j) h.
  \end{equation*}
  As the method is first-order, the local truncation error is given
  by:
  \begin{equation}
    \label{eq:psi_r_3_Euler}
    r_{j+1}(h) = \psi( y_j) h^2 +  r^{(3)}_{j+1}(
    \theta h) \frac{h^{3}}{6},
  \end{equation}
  for some $\theta \in (0,1)$. The following provides explicit
  formulations of $\psi$ and $r^{(3)}_{j+1}$ solely based on the vector field
  and its derivatives:
  \begin{align}
    \psi(y_j) &= \frac{f'(y_j)f(y_j)}{2}, \nonumber \\
    r^{(3)}_{j+1}(\theta h) &= f''(y_j + \theta h) \left( f(y_j +
                              \theta h) \right)^2 + \left( f'(y_j + \theta h) \right)^2
                              f(y_j + \theta h). \label{eq:r3_Euler_in_terms_of_f_derivs}
  \end{align}

  In the analysis of Euler method according to Runge-Kutta theory, the
  vector field $f : [-K,K] \to [-M,M]$ is required to be twice continuously
  differentiable, which we write as $f \in C^2([-K,K])$. As the domain
  $[-K,K]$ is compact, there are two positive constants $M_1$ and
  $M_2$ satisfying:
  \begin{equation}
    \label{def:Euler_deriv_bounds_C2}
    \forall x \in [-K,K]: \quad \absn{f'(x)} \leq M_1 \wedge \absn{f''(x)} \leq M_2.
  \end{equation}
  By considering equations~\eqref{eq:psi_r_3_Euler}
  and~\eqref{eq:r3_Euler_in_terms_of_f_derivs}, we introduce the
  constant $\alpha$ as follows:
  \begin{equation*}
    \alpha \defeq \frac{(M_2 M + M_1^2 )M}{6}.
  \end{equation*}
  To obtain a validated version, let $Y$, $F$, and $F'$ be interval
  extensions of $y$, $f$, and $f'$, respectively. Then, we obtain an
  interval extension $\Psi$ of $\psi$ by defining:
  \begin{equation*}
    \Psi(x) \defeq \frac{F'(x)F(x)}{2}.
  \end{equation*}
  \noindent
  In the general schema presented in
  Section~\ref{sec:domain_temporal_discretization}, the
  formula~\eqref{eq:Y_h_I_j} now admits the following explicit
  form:
  \begin{equation*}
    Y(q_j + h) = Y(q_j) + F( Y(q_j)) h + \Psi(Y(q_j)) h^2
    + [-\alpha, \alpha] h^3 .
  \end{equation*}
  
  Although the classical Runge-Kutta theory requires
  $f \in C^2([-K,K])$, it is possible to relax the condition slightly
  for the validated version. To that end, we recall the following
  generalization of classical Taylor's theorem:
  
\begin{lemma}
  \label{lemma:Taylor_Lipschitz_OneVar}
  Assume that $p \in \N$. Let $f: [\alpha, \beta] \to \R$ be a
  function that is $p$-times continuously differentiable and suppose
  that the $p$-th derivative $f^{(p)}$ of $f$ is Lipschitz in
  $(\alpha, \beta)$. Then:
  \begin{equation*}
    f(\beta) \in  \sum_{j = 0}^p  f^{(j)}(\alpha) \cdot \frac{(\beta-\alpha)^j}{j!} + 
    L(f^{(p)})(\theta) \cdot \frac{(\beta-\alpha)^{p+1}}{(p+1)!},
  \end{equation*}
  for some $\theta \in (\alpha, \beta)$.
\end{lemma}

\begin{proof}
  See~\cite[Lemma~2.13]{Edalat_Farjudian_Mohammadian_Pattinson:2nd_Order_Euler:2020:Conf}.
\end{proof}

Thus, we may just require $f$ to be in $C^1([-K,K])$, with $f'$
Lipschitz continuous, and replace~\eqref{def:Euler_deriv_bounds_C2}
with the following:
  \begin{equation*}
    \forall x \in [-K,K]: \quad \absn{f'(x)} \leq M_1 \wedge \norm{L(f')(x)} \leq M_2,
  \end{equation*}
  in which $L(f')$ is the $L$-derivative of $f'$, and $\norm{\cdot}$
  is the interval norm from Definition~\ref{def:interval_norms}.

  We now proceed to define the counterpart $\Phi^{\mathrm{R}}$ of the
  operator $\Phi$ from Definition~\ref{def:Phi}:\footnote{The
    superscript `$\mathrm{R}$' stands for Runge-Kutta.}

\begin{definition}[Operator: $\Phi^{\mathrm{R}}$]
  \label{def:Phi_R}
  For a given partition
  $Q \equiv ( q_0, \ldots, q_k) \in {\mathcal{P}}_{\Q}$, a given
  triple $\left( u, u',u''\right) \in \hat{D}^{(2)}$, and
  $\phi \in {\mathcal{D}}_{\Q}$, we define:
  \begin{equation}
    \label{eq:Phi_R}
    y_\phi( x) \defeq \left\{
      \arrayoptions{2ex}{1.3}
      \begin{array}{ll}
        0, & \text{if } x = 0,\\
        T_{\negthinspace K,1} [ \phi( q_j) + u(\phi( q_j))(x-q_j)
        + & \\
        \qquad (u' \cdot u)(\phi(q_j)) \frac{(x-q_j)^2}{2} + [-\alpha, \alpha](x-q_j)^3 ]   , & \text{if } q_j < x \leq q_{j+1},
      \end{array}
      \right.
    \end{equation}    
    \noindent
    where:
    \begin{itemize}

    \item $\hat{D}^{(2)}$ is the domain from
      Definition~\ref{def:domain_of_Cp_funs}.

    \item $\alpha \defeq \frac{(M_2 M + M_1^2 )M}{6}$, in which $M$,
      $M_1$, and $M_2$ are from Definition~\ref{def:domain_of_Cp_funs}.

    \item $T_{\negthinspace K,1}$ is as defined in~\eqref{eq:T_K_n}.

    \item $(u' \cdot u)( \cdot)$ denotes the product of the interval
      $u'( \cdot)$ with the interval $u( \cdot )$.
    \end{itemize}
    The operator
    $\Phi^{\mathrm{R}} :{\mathcal{P}}_{\Q} \times \hat{D}^{(2)} \to
    {\mathcal{D}}_{\Q} \to \left([0,a] \to \intvaldom[{[-K,K]}]
    \right)$ is defined by:
    \begin{equation*}
      \Phi^{\mathrm{R}}_{(u,u',u'')}(Q)(\phi) \defeq y_\phi. 
    \end{equation*}
  \end{definition}
  
  \begin{remark}
    Although $u''$ does not appear explicitly in~\eqref{eq:Phi_R}, it
    is implicitly used via the value $\alpha$ which depends on $M_2$
    (from Definition~\ref{def:domain_of_Cp_funs}), which is, in turn,
    an upper bound on the values that $u''$ takes.
  \end{remark}
  
  Having defined the operator $\Phi^{\mathrm{R}}$, we may proceed by
  taking the same steps as those taken in
  Section~\ref{susec:2nd_Order_Euler_Operator}. For instance, the
  counterpart of Proposition~\ref{prop:Phi_D_QL_unif_ideal_prop} may
  be stated as follows:
  
  \begin{proposition}
    \label{prop:Phi_R_D_QL_unif_ideal_prop}
    Assume that $Q \equiv ( q_0, \ldots, q_k) \in {\mathcal{P}}_{\Q}$ and
    $( u, u', u'') \in \hat{D}^{(2)}$. Then:
    \begin{equation}
      \label{eq:Phi_R_range_D_QL}
      \forall \phi \in {\mathcal{D}}_{\Q}: \quad
      \Phi^{\mathrm{R}}_{(u,u',u'')}(Q)(\phi) \in {\mathcal{D}}_{\Q}.
    \end{equation}
    Furthermore,
    $\Phi^{\mathrm{R}}_{(u,u',u'')}(Q): {\mathcal{D}}_{\Q} \to
    {\mathcal{D}}_{\Q}$ has the \ac{UIC} property.
  \end{proposition}

  \begin{proof}
    From~\eqref{eq:Phi_R}, by using an argument similar to that used
    in the proof of Proposition~\ref{prop:Phi_D_QL_unif_ideal_prop},
    we deduce that the bounds $\lep[y_\phi]$ and $\uep[y_\phi]$ are
    continuous over each $(q_j,q_{j+1}]$. This, together with the
    assumption $Q \in {\mathcal{P}}_{\Q}$,
    implies~\eqref{eq:Phi_R_range_D_QL}.

    To prove the \ac{UIC} property, first we note that, as all the
    operations in~\eqref{eq:Phi_R} are monotonic, so is
    $\Phi^{\mathrm{R}}_{(u,u',u'')}(Q)$. Furthermore, for any given
    $\psi \in {\mathcal{W}}_{\mathcal{D}}$, as $\psi$ is rounded,
    then, by~\eqref{eq:def_phi_star_X_to_D}, together with the fact
    that each partition may have only finitely many partition points,
    we obtain:
      \begin{equation}
        \label{eq:d_b_eps_phi_star_R}
        \forall \epsilon > 0: \exists b_{\epsilon} \in \psi: \forall i
        \in \set{ 0, \ldots, k}: \quad T_{\negthinspace K,1}
        (\psi^*(q_i) \symExpand \epsilon) 
        \sqsubseteq b_{\epsilon}(q_i).
      \end{equation}
      From~\eqref{eq:d_b_eps_phi_star_R}, combined with Scott
      continuity of $u$ and $u'$, we deduce:
      \begin{equation*}
        \forall \delta > 0: \exists b_{\delta} \in \psi: \forall t \in
        [0,a]: \quad
 T_{\negthinspace K,1}\left( \Phi^{\mathrm{R}}_{(u,u',u'')}(Q)(\psi^*)(t) \symExpand
   \delta \right) \sqsubseteq
          \Phi^{\mathrm{R}}_{(u,u',u'')}(Q)(b_{\delta})(t).
      \end{equation*}
  \end{proof}
  
  \begin{definition}[Euler Operator $F^{\mathrm{R}}$]
  \label{def:F_R}
  The scalar Runge-Kutta Euler operator:
  \begin{equation*}
    F^{\mathrm{R}} : {\mathcal{P}}_{\Q} \times \hat{D}^{(2)} \to {\mathcal{W}}_{\mathcal{D}} \to {\mathcal{W}}_{\mathcal{D}}  
\end{equation*}
is defined as follows: for any given
$Q \equiv ( q_0, \ldots, q_k) \in {\mathcal{P}}_{\Q}$,
$( u, u', u'') \in \hat{D}^{(2)}$, and $\phi \in {\mathcal{W}}_{\mathcal{D}}$:
\begin{equation}
  \label{eq:def_F_R}
  F^{\mathrm{R}}_{\negthinspace (u,u',u'')}(Q)(\phi) \defeq
  \left( \Phi^{\mathrm{R}}_{(u,u',u'')}(Q)(\phi^*)\right)_* .
    \end{equation}
  \end{definition}

  Once again, the relationship between the operators $F^{\mathrm{R}}$,
  $\Phi^{\mathrm{R}}$, and the left and right adjoints $(\cdot)^*$ and
  $(\cdot)_*$ of the Galois connection from
  Theorem~\ref{thm:Galois_connection_X_to_D} can be depicted as in the
  following commutative diagram:

\begin{equation*}
  \begin{tikzcd}[row sep = large, column sep = large]
    {\mathcal{D}}_{\Q} \arrow[r, hookrightarrow, "(\cdot)_*"]    
    & {\mathcal{W}}_{\mathcal{D}}\\
    {\mathcal{D}}_{\Q} \arrow[u, "\Phi^{\mathrm{R}}_{(u,u',u'')}(Q)(\cdot)"] &
    {\mathcal{W}}_{\mathcal{D}} \arrow[u, swap, "F^{\mathrm{R}}_{\negthinspace (u,u',u'')}(Q)(\cdot)"] \arrow[l,
    twoheadrightarrow,"(\cdot)^*" ]
  \end{tikzcd}
\end{equation*}
  
\begin{lemma}
  \label{lemma:F_R_Scott_Continuous}
  For every partition $Q \equiv ( q_0, \ldots, q_k) \in {\mathcal{P}}_{\Q}$
  and $( u, u', u'') \in \hat{D}^{(2)}$, the function
  $F^{\mathrm{R}}_{\negthinspace (u,u',u'')}(Q): {\mathcal{W}}_{\mathcal{D}} \to
  {\mathcal{W}}_{\mathcal{D}}$ is Scott continuous.
\end{lemma}
\begin{proof}
    This follows from Proposition~\ref{prop:Phi_R_D_QL_unif_ideal_prop}
  and Lemma~\ref{lem:F_Phi_general_result}.
\end{proof}

As
$F^{\mathrm{R}}_{\negthinspace (u,u',u'')}(Q): {\mathcal{W}}_{\mathcal{D}} \to {\mathcal{W}}_{\mathcal{D}}$ is Scott-continuous, we may define the Runge-Kutta Euler
operator using its fixpoint:

\begin{definition}[Runge-Kutta Euler Operator: $E^{\mathrm{R}}$]
  \label{def:E_R_Ideal_Compl}
  Assume that $Q \equiv ( q_0, \ldots, q_k) \in {\mathcal{P}}_{\Q}$ is a
  partition of~$[0,a]$, $( u, u', u'') \in \hat{D}^{(2)}$, and $\bot$
  denotes the bottom element of ${\mathcal{W}}_{\mathcal{D}}$. We define
  $ E^{\mathrm{R}}_{(u,u',u'')}(Q) \defeq \psi^*$, where:
  \begin{equation*}
    \psi \defeq \fix F^{\mathrm{R}}_{\negthinspace (u,u',u'')}(Q) = \bigsqcup_{m \in \N} \left(
      F^{\mathrm{R}}_{\negthinspace (u,u',u'')}(Q) \right)^m (\bot).
  \end{equation*}  
\end{definition}

Similar to Theorem~\ref{thm:E_2_MFPS_Ideal_Compl}, we can prove that:
\begin{equation*}
  \bigsqcup_{m \in \N} \left( F^{\mathrm{R}}_{\negthinspace (u,u',u'')}(Q) \right)^m (\bot) = \left( F^{\mathrm{R}}_{\negthinspace (u,u',u'')}(Q) \right)^{k+1} (\bot),
\end{equation*}
which provides us with a validated Euler operator based on Runge-Kutta
theory.

For any fixed $Q \in {\mathcal{P}}_{\Q}$, computability of the map
$F_{\negthinspace ( \cdot , \cdot, \cdot)}(Q)(\cdot): \hat{D}^{(2)}
\times {\mathcal{W}}_{\mathcal{D}} \to {\mathcal{W}}_{\mathcal{D}}$
can be proven using an approach similar to that taken in
Section~\ref{subsec:computability}. The key fact is that
$\alpha x^3 + \beta x^2 + \gamma x + \rho < \delta$ is semi-decidable
over an interval $[p,q]$, with
$\alpha, \beta, \gamma, \rho, \delta \in \Q$, and with $p$ and $q$
computable. In fact, using the approach of
Section~\ref{subsec:computability}, one may prove computability of any
Runge-Kutta variant as long as the bounds within each interval are
polynomials with rational coefficients. This is because the proof of
computability essentially reduces to deciding
$\sum_{i=0}^m \alpha_i x^i < \delta$ over intervals of the form
$[p,q]$, where $\alpha_0, \ldots, \alpha_m, \delta \in \Q$, and with
$p$ and $q$ computable. When $p, q \in \Q$, this is decidable
according to Tarski's
theorem~\parencite{Tarski:Decision_Method:1998}. When $p$ and $q$ are
computable, the problem becomes semi-decidable.

Soundness and completeness follow from the relevant error analysis as
presented in~\parencite{Marciniak:Selected_Interval_Methods:2009},
where error analyses of several variants of Runge-Kutta method may be
found. As for convergence analysis, it is well-known that the
Runge-Kutta formulation of the Euler operator provides a first-order
convergence
rate~\cite[Theorem~5.4]{Lambert:Numerical_Methods:Book:1991}, a
fact which will be verified by our experiments in
Section~\ref{sec:Experiments} as well (in particular,
Figure~\ref{fig:cos_Euler_RK}).

\section{Experiments}
\label{sec:Experiments}

The domain theoretic framework of the current article makes it
possible to use effective data-types, {\ie}, it is possible to
implement the operators directly on digital computers. We have indeed
implemented the following operators for \ac{IVP} solving:
\begin{enumerate}[label=(\arabic*)]
\item The first-order Euler operator $E^c$
  of~\parencite{EdalatPattinson2006-Euler-PARA};
\item The second-order Euler operator $E^2$ of
  Definition~\ref{def:2nd_order_Euler}.
\item The Runge-Kutta (Euler) operator $E^{\mathrm{R}}$ of Definition~\ref{def:E_R_Ideal_Compl}.
\end{enumerate}
The source code is available on
GitHub.\footnote{\url{https://github.com/afarjudian/IVP_MPFI}.} We
have used the arbitrary-precision interval arithmetic library
MPFI~\parencite{Revol_Rouillier:MPFI:05} for our
implementations. Specifically, we have used the C++ Boost library
implementation which provides a wrapper around the original MPFI
types.\footnote{\href{https://www.boost.org/doc/libs/1_76_0/libs/multiprecision/doc/html/boost_multiprecision/tut/interval/mpfi.html}{https://www.boost.org/doc/libs/1\_76\_0/libs/multiprecision/doc/html/boost\_multiprecision/tut/interval/mpfi.html}}
This is an ideal library for our purposes as operations such as
matrix/vector arithmetic can be implemented seamlessly over the
underlying MPFI types.

We consider the following \acp{IVP} for our experiments:

\begin{enumerate}[label=(\Alph*)]
\item To begin with, we consider the simple scalar \ac{IVP}:
  
  \begin{equation}
    \label{eq:id_ivp}
    \left\{
      \arrayoptions{0.5ex}{1.3}
      \begin{array}{r@{\hspace{0.5ex}=\hspace{0.5ex}}l}
        y'(t) & y(t),\\
        y(0) & 1.\\
      \end{array}
    \right.
  \end{equation}
  This \ac{IVP} has the closed-form solution $y'(t) = \exp(t)$ over
  the entire real line. We solve this equation over the interval
  $[0,1]$, over which we may take $M=3, M_1 = 1$, and $M_2 = 0$.

\item Next, we consider another scalar \ac{IVP}:

    \begin{equation}
  \label{eq:cos_ivp}
    \left\{
      \arrayoptions{0.5ex}{1.3}
      \begin{array}{r@{\hspace{0.5ex}=\hspace{0.5ex}}l}
        y'(t) & \cos( y(t)),\\
        y(0) & 0,\\
      \end{array}
    \right.
  \end{equation}
  which has the following closed form solution over the entire real line:
  \begin{equation*}
    y(t) = 2 \atan\left( \tanh\left( \frac{t}{2} \right) \right).
  \end{equation*}
  We solve this \ac{IVP} over the interval $[0,5]$. As the vector field is
  $f(x) = \cos(x)$, we take $M = M_1 = M_2 = 1$.

\item The next \ac{IVP} is slightly more complicated:

  \begin{equation}
  \label{eq:expSin_ivp_non_autonomous}
    \left\{
      \arrayoptions{0.5ex}{1.3}
      \begin{array}{r@{\hspace{0.5ex}=\hspace{0.5ex}}l}
        \hat{y}'(t) & 10 \cos( 10 t) \hat{y}(t),\\
        \hat{y}(0) & 1.\\
      \end{array}
    \right.
  \end{equation}
  This \ac{IVP} also has a closed-form solution
  $\hat{y}(t) = \exp( \sin(10 t))$ over the entire real line. The
  equation~\eqref{eq:expSin_ivp_non_autonomous} is non-autonomous. By
  assigning $y(t) \defeq ( t, \hat{y}(t))$, we obtain an autonomous
  non-scalar equation with the vector field
  $f( \alpha, \beta) = ( 1, 10 \cos( 10 \alpha) \beta )$ and initial
  value $y(0) = ( 0 , 1)$. As the autonomous equation is non-scalar,
  the Runge-Kutta operator is not applicable. We seek a solution in
  the interval $[0, 0.1]$ for which we take $M=30$ for the Euler
  operators.

\item We also consider an \ac{IVP} with a non-differentiable vector field:
  \begin{equation}
    \label{eq:absSin_ivp_non_autonomous}
    \left\{
      \arrayoptions{0.5ex}{1.3}
      \begin{array}{r@{\hspace{0.5ex}=\hspace{0.5ex}}l}
        \hat{y}'(t) & \absn{\sin( t+\hat{y}(t))},\\
        \hat{y}(0) & 1.\\
      \end{array}
    \right.
  \end{equation}
  We are not aware if this \ac{IVP} has a closed form
  solution. Equation~\eqref{eq:absSin_ivp_non_autonomous} is also
  non-autonomous. By assigning $y(t) \defeq ( t, \hat{y}(t))$, we
  obtain an autonomous equation with the vector field
  $f( \alpha, \beta) = ( 1, \absn{\sin(\alpha + \beta)} )$ and initial
  value $y(0) = ( 0 , 1)$. We seek a solution in the interval $[0, 5]$
  for which we take $M=1$ for the Euler operators.

  Over the interval $[0,5]$, the presence of the absolute value
  function does indeed make a difference, as shown in
  Figure~\ref{fig:absSinPlots}.

\begin{figure}[t]
  \centering
  \scalebox{0.17}[0.15]{\includegraphics{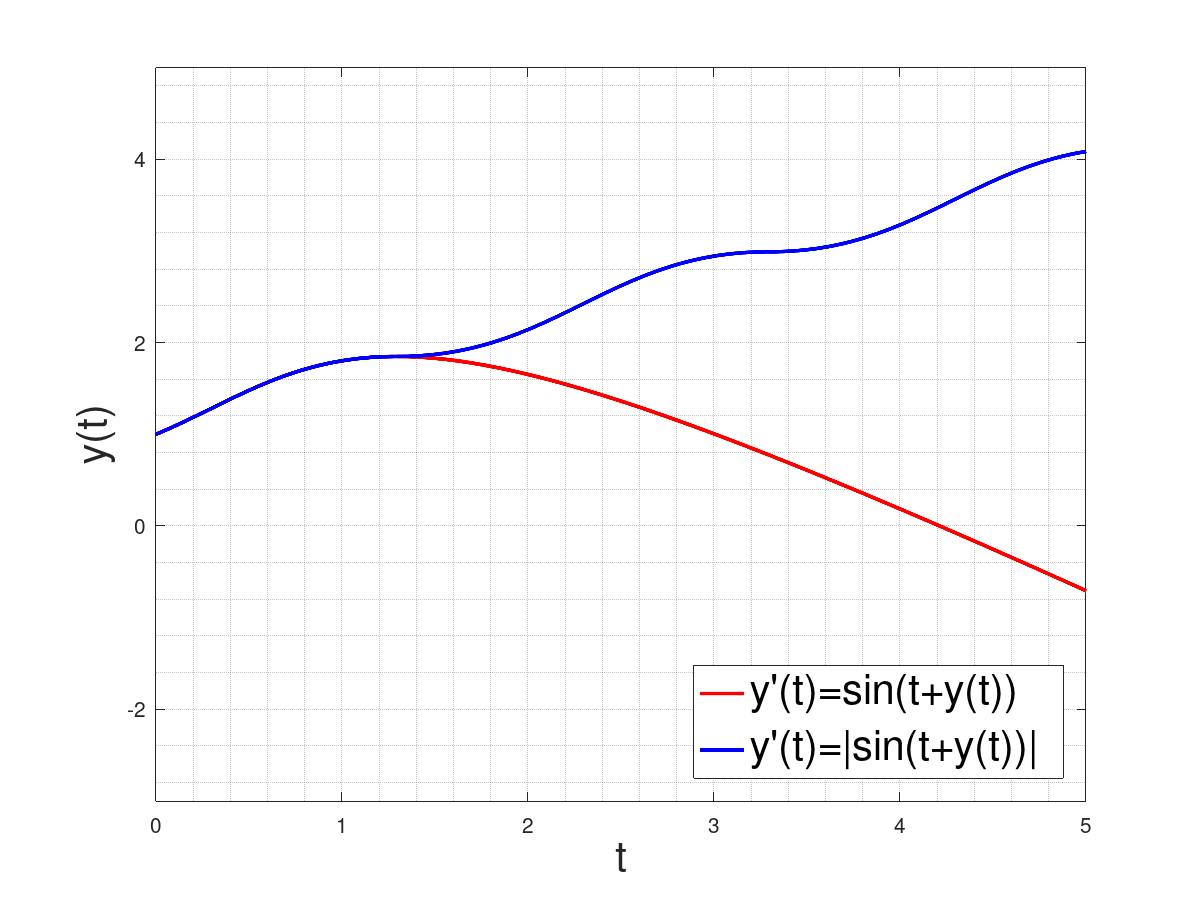}}
  \caption{The solution of $y'(t) = \sin(t+y(t))$ in
    {\color{red}{red}} versus that of $y'(t) = \absn{\sin(t+y(t))}$ in
    {\color{blue}{blue}}, over the domain $[0,5]$. Around $t=1.295$, the
    two solutions diverge.}
  \label{fig:absSinPlots}
\end{figure}

\end{enumerate}

  In our experiments, we consider a parameter \emph{depth}, which
  denotes how many times the domain $[0,a]$ must been bisected. For
  instance, when running any of the operators at depth $4$, the domain
  $[0,a]$ is discretized into $2^4=16$ equal sized subintervals. We
  ran the first-order operator $E^c$
  of~\parencite{EdalatPattinson2006-Euler-PARA} and the second-order
  Euler operator $E^2$ on all of the \acp{IVP}, and ran the
  Runge-Kutta operator $E^{\mathrm{R}}$ on the scalar
  \acp{IVP}~\eqref{eq:id_ivp} and~\eqref{eq:cos_ivp}, for a range of
  depths from $2$ to $16$.

\begin{figure}[h]
  \centering
  \scalebox{0.3}[0.18]{\includegraphics{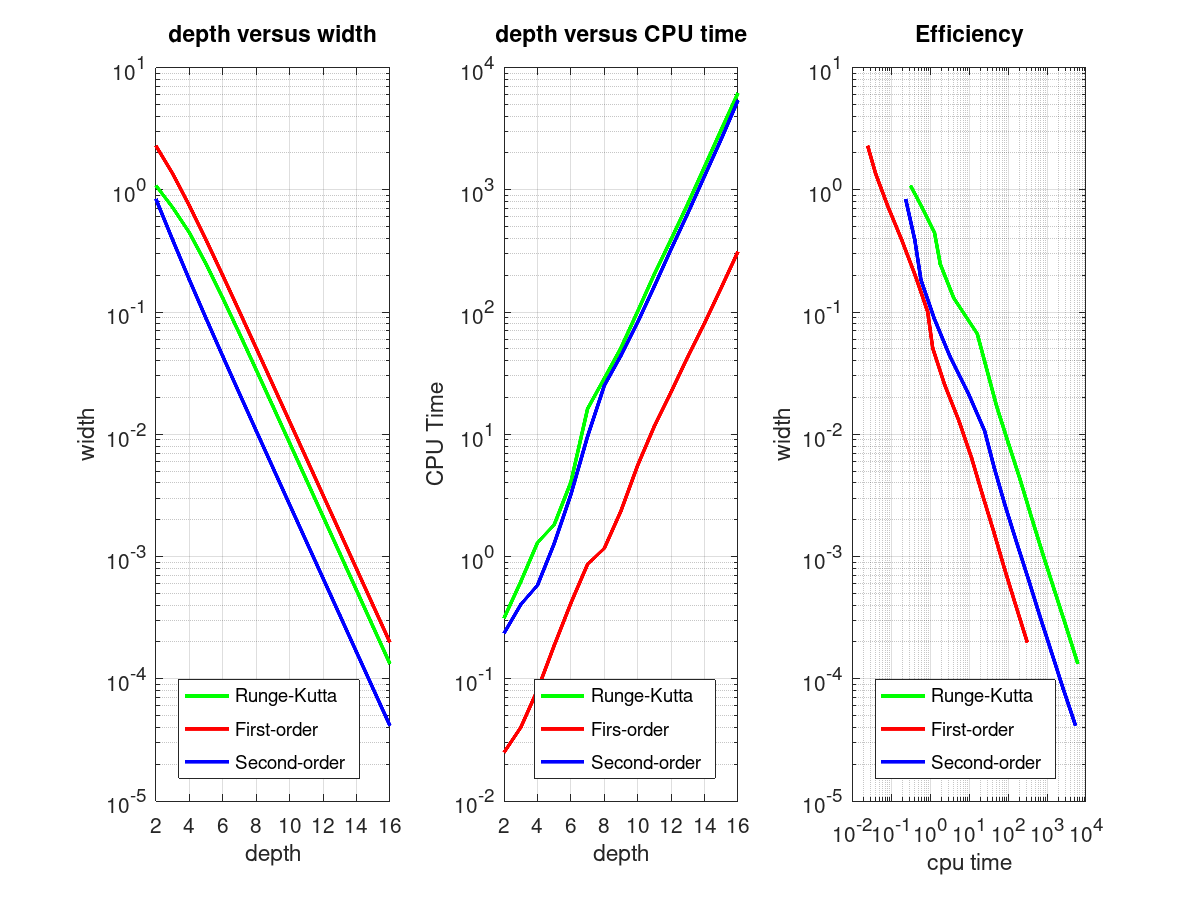}}
  \caption{Comparison of the first-order, second-order, and
    Runge-Kutta Euler methods on the \ac{IVP}~\eqref{eq:id_ivp}, with
    $y'(t) = y(t)$.}
  \label{fig:identity_Euler_RK}
\end{figure}

Figure~\ref{fig:identity_Euler_RK} contains three plots related to the
\ac{IVP}~\eqref{eq:id_ivp} with $y'(t) = y(t)$:

\begin{itemize}
  
\item The left plot shows that at each depth, the width of the
  solution obtained from $E^2$ is noticeably smaller than the widths
  obtained from the Runge-Kutta $E^{\mathrm{R}}$ and the first-order
  operator $E^c$.
  
\item The middle plot shows that at equal depths, the second-order and
  Runge-Kutta methods take more time compared with the first-order
  method. This is expected due to the overhead of handling
  derivatives.
  
\item The right plot shows that, in the particular case of the
  \ac{IVP}~\eqref{eq:id_ivp}, the overall performance of the
  first-order operator $E^c$ is better. In simple terms, by spending
  the same amount of CPU time, the first-order method provides tighter
  enclosures of the solution.
\end{itemize}


\begin{figure}[h]
  \centering
  \scalebox{0.3}[0.18]{\includegraphics{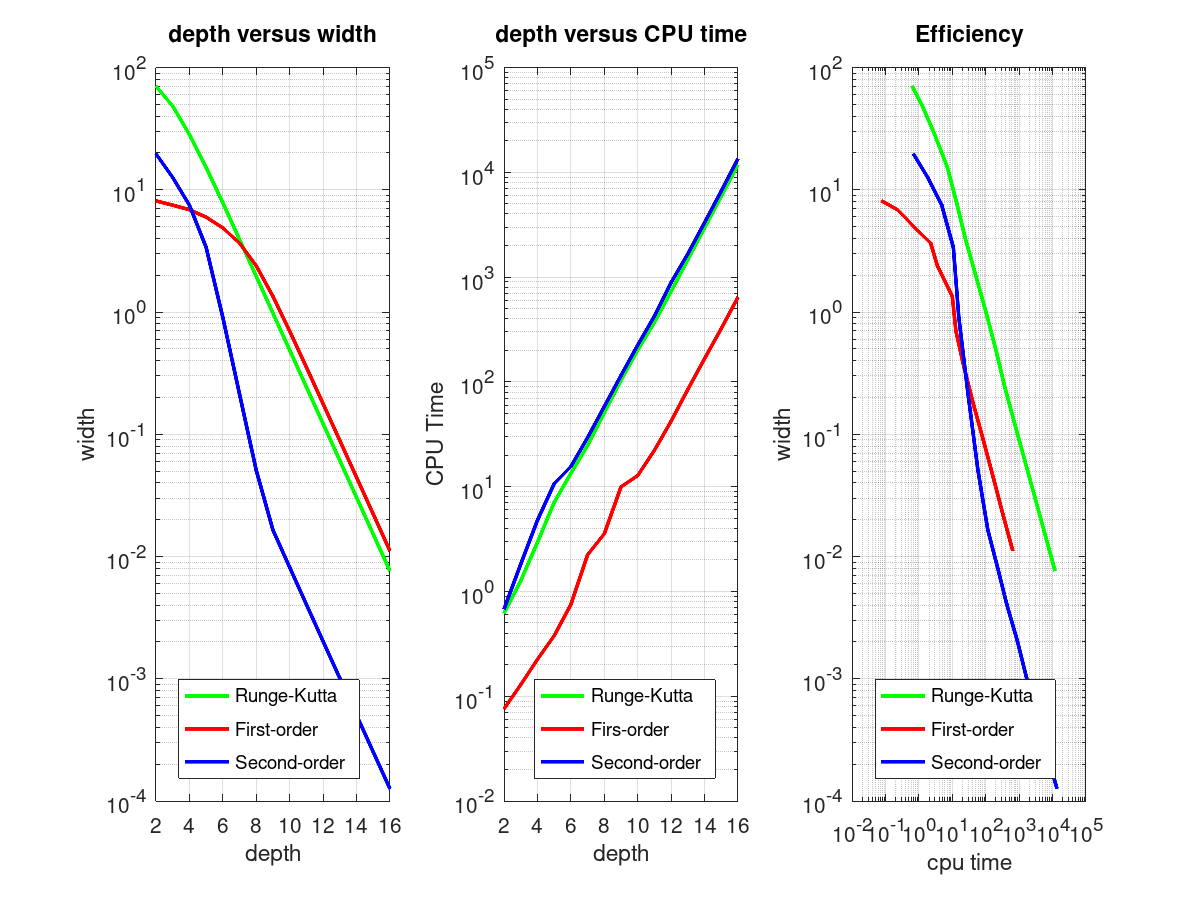}}
  \caption{Comparison of the first-order, second-order, and
    Runge-Kutta Euler methods on the \ac{IVP}~\eqref{eq:cos_ivp}, with
    $y'(t) = \cos(y(t))$.}
  \label{fig:cos_Euler_RK}
\end{figure}

Figure~\ref{fig:cos_Euler_RK} contains the plots related to the
\ac{IVP}~\eqref{eq:cos_ivp}. The right plot demonstrates that at lower
depths, the first-order operator $E^c$ is the most
efficient. Nonetheless, at higher depths, the second-order operator
outperforms both the first-order and Runge-Kutta operators due to its
superior convergence rate.

\begin{figure}[h]
  \centering
  \scalebox{0.3}[0.18]{\includegraphics{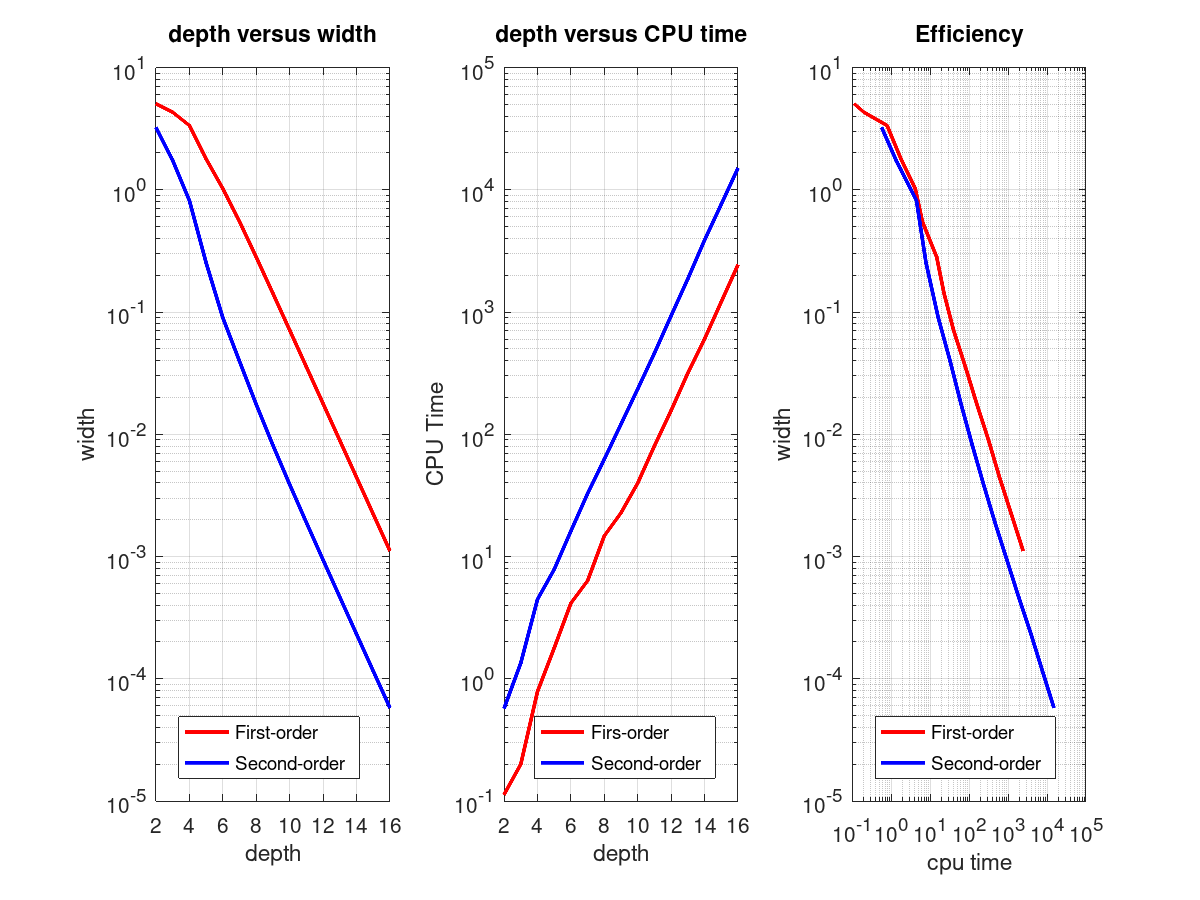}}
  \caption{Comparison of the first-order and second-order Euler
    methods on the \ac{IVP}~\eqref{eq:expSin_ivp_non_autonomous} with
    $\hat{y}'(t) =  10 \cos( 10 t) \hat{y}(t)$.}
  \label{fig:expSin}
\end{figure}

Figure~\ref{fig:expSin} contains the plots related to the
\ac{IVP}~\eqref{eq:expSin_ivp_non_autonomous}. As the \ac{IVP} is
converted to a non-scalar one, we may only apply the first-order and
second-order operators $E^c$ and $E^2$, respectively. The left and
middle plots have the same quality as those of
Figure~\ref{fig:identity_Euler_RK}. The right plot, however, shows
that, in this case, the second-order method outperforms the
first-order Euler method in overall efficiency.

\begin{figure}[h]
  \centering
  \scalebox{0.3}[0.18]{\includegraphics{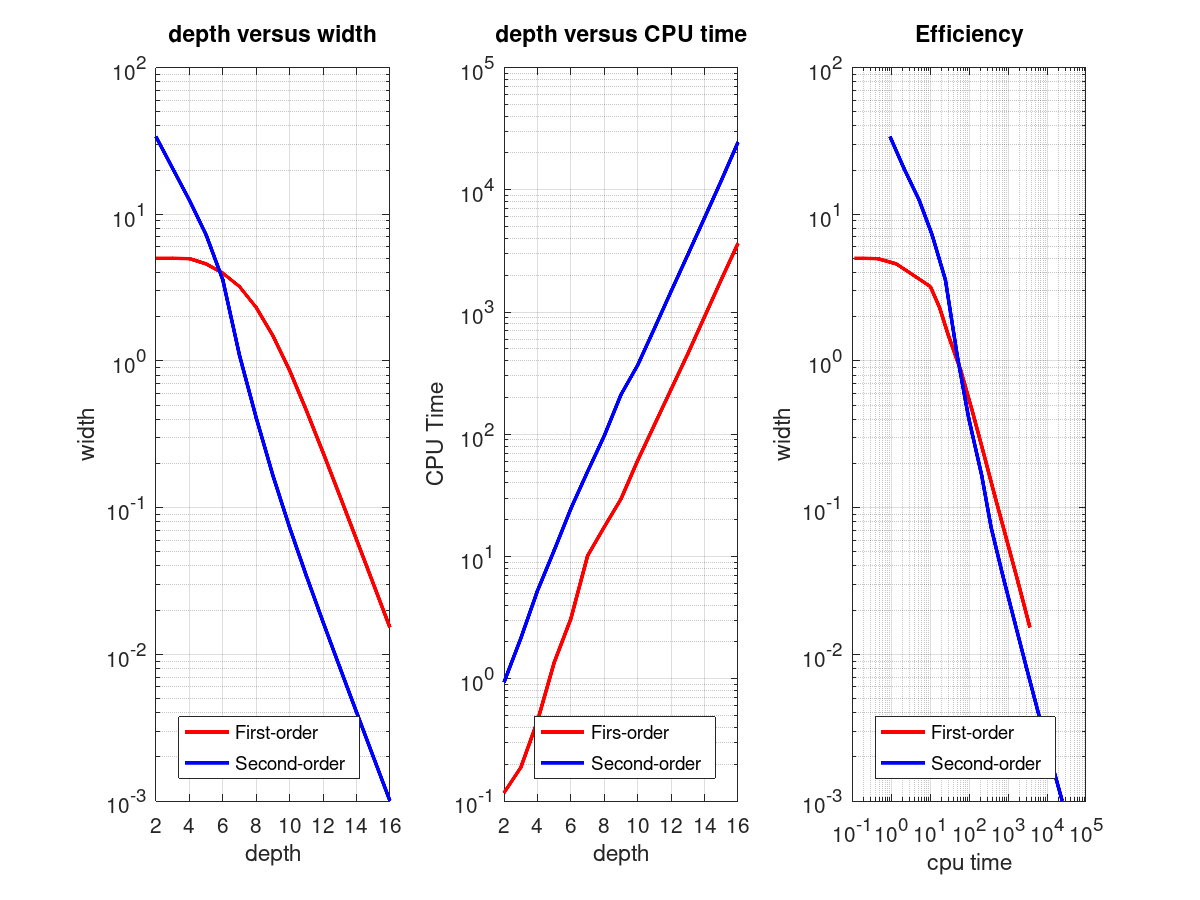}}
  \caption{Comparison of the first-order and second-order Euler
    methods on the \ac{IVP}~\eqref{eq:absSin_ivp_non_autonomous}, with
    a non-differentiable vector field: $\hat{y}'(t) = \absn{\sin( t+\hat{y}(t))}$.}
  \label{fig:absSin}
\end{figure}

Figure~\ref{fig:absSin} contains the relevant plots for the
\ac{IVP}~\eqref{eq:absSin_ivp_non_autonomous}. As the \ac{IVP} is
converted to a non-scalar one, and the vector field is not continuously
differentiable, we may only apply the first-order and second-order
operators $E^c$ and $E^2$, respectively. The middle plot is as
expected, but the left and right plots demonstrate the effect of the
overhead of handling derivatives in the second-order method. In lower
depths, the first-order Euler operator performs better, but as the
depth is increased, the faster convergence of the second-order
operator makes it overall more efficient.

\begin{figure}[h]
  \centering
  \scalebox{0.18}[0.15]{\includegraphics{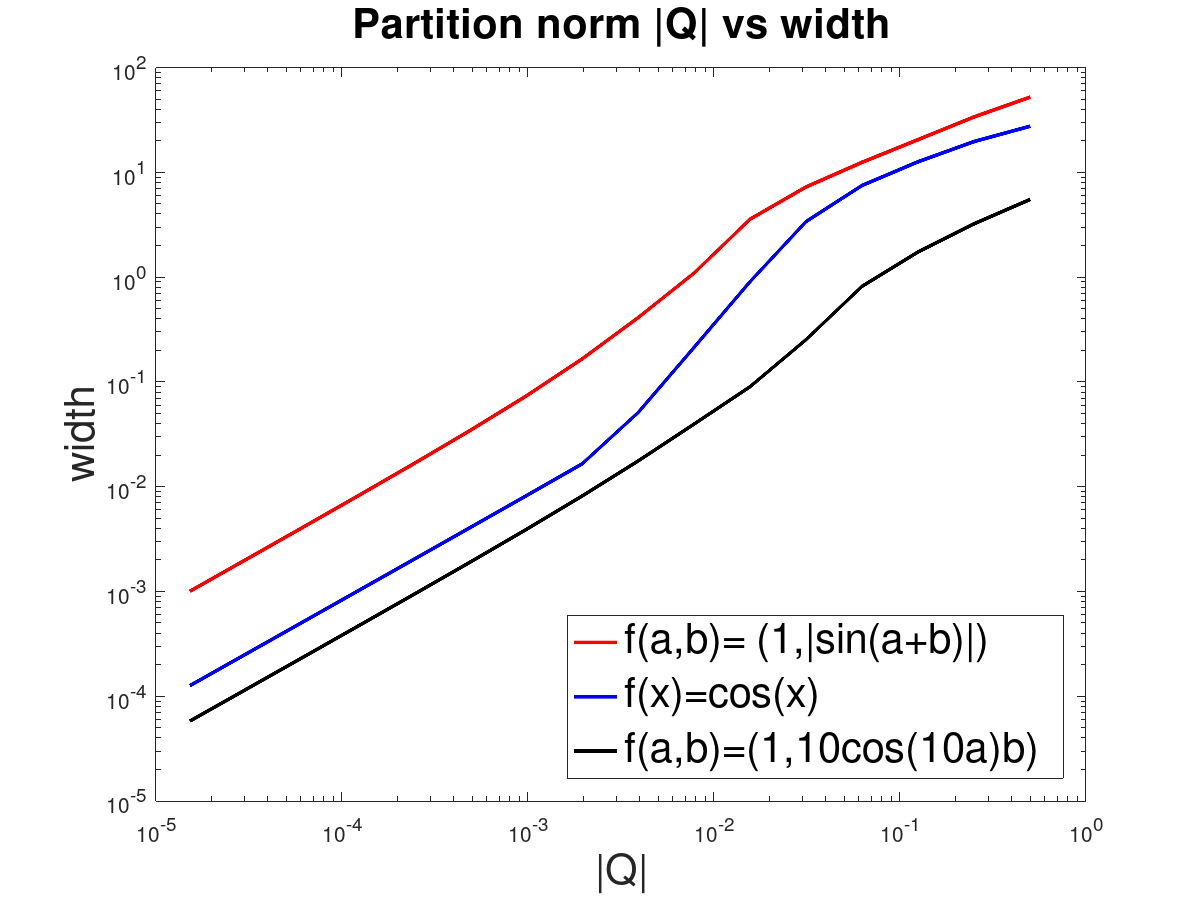}}
  \caption{Comparison of the convergence of the second-order Euler
    method on the \ac{IVP}~\eqref{eq:absSin_ivp_non_autonomous}, with
    a non-differentiable vector field (in {\color{red}{red}}), and the
    \acp{IVP}~\eqref{eq:cos_ivp}
    and~\eqref{eq:expSin_ivp_non_autonomous}, with vector fields that are
    continuously differentiable (in {\color{blue}{blue}} and
    {\color{black}{black}}). The second-order convergence seems to be
    retained for the \ac{IVP}~\eqref{eq:absSin_ivp_non_autonomous}.}
  \label{fig:id_versus_absSin}
\end{figure}

As we have already pointed out in
Section~\ref{subsubsec:non_differentiable_fields}, although
Theorems~\ref{thm:Second_Order_Convergence_int_Lip}
and~\ref{thm:speed_of_convergence_field_approx} require the vector field to
be continuously differentiable, we conjecture that the second-order
convergence is retained over non-differentiable vector fields as well. As a
witness to this conjecture, we compare the convergence of the
second-order Euler method for solving the
\ac{IVP}~\eqref{eq:absSin_ivp_non_autonomous}--with a
non-differentiable vector field---with those obtained from solving
\acp{IVP}~\eqref{eq:cos_ivp} and~\eqref{eq:expSin_ivp_non_autonomous},
with vector fields that are continuously differentiable. Figure~\ref{fig:id_versus_absSin} suggests
that the convergence has the same order in this case.

\subsection{Static Analysis}

The bound~\eqref{eq:conv_speed_ENTCS_equidistant_Q} (obtained
in~\parencite{Edalat_Farjudian_Mohammadian_Pattinson:2nd_Order_Euler:2020:Conf})
and those obtained in the current work ({\ie},
Corollary~\ref{cor:speed_convergence} and
Theorem~\ref{thm:speed_of_convergence_field_approx}) make it possible
to perform a static convergence analysis. These bounds, however, are
usually very conservative, and result in unnecessarily small partition
norms for obtaining a required accuracy. In practice, it is more
efficient to just run the Euler operators for various depths, in
increasing order, until the required accuracy is obtained.

Take the \ac{IVP}~\eqref{eq:expSin_ivp_non_autonomous} as an
example. A static analysis based on the
bound~\eqref{eq:conv_speed_ENTCS_equidistant_Q} would indicate that,
in order to obtain $\epsilon > 0$ accuracy, one must take the
partition norm small enough to satisfy:
\begin{equation*}
  |Q| \leq \frac{\epsilon}{15 (\me^{60} - 1)} \leq
  \frac{\epsilon}{1.72 \times 10^{27}}.
\end{equation*}
As an example, according to this estimate, to obtain an accuracy of
$\epsilon \leq 10^{-4}$, we must partition $[0,0.1]$ to at least the
depth of $100$, resulting in $2^{100}$ subintervals, which is
prohibitively large. Instead, we start from a small depth of $2$, and
increase the depth until the accuracy is reached. With this approach,
the accuracy reaches the target of $10^{-4}$ at depth $16$, in under
$30$ seconds on a personal laptop, with an
{Intel\textsuperscript{\textregistered}} {Core\texttrademark} i7-8550U
CPU at 1.80GHz and 16GB RAM, under Ubuntu 20.04 LTS environment.

\begin{remark}
  \label{remark:comparison_with_VNODE_etc}
  Although we have reported the timing of our experiments and the
  hardware that has been used, the numbers should be considered only
  in the context of comparison between the algorithms that we have
  implemented. The Euler operator $E^2$ is second-order, and as such,
  over \acp{IVP} with highly smooth vector fields, it may not exhibit
  the efficiency of the high-order methods such as
  COSY~\parencite{Berz_Makino:VerifiedIO:1998},
  VNODE~\parencite{Nedialkov_et_al:High_Order_Interval_ODE:2001} or
  ValEncIA-IVP~\parencite{RauhValencia-SCAN.2006.47}. On the other
  hand, our method is superior in its generality as it can provide
  validated solution of \acp{IVP} with non-differentiable vector
  fields as well, while the high-order methods available in the
  literature invariably require the vector field to be continuously
  differentiable.
\end{remark}


\section{Concluding Remarks}
\label{sec:concluding_remarks}

The analyses presented in the current article are applicable to any
operator that follows the general schema presented in
Section~\ref{sec:domain_temporal_discretization}. Given the
superiority of the Euler operator $E^2$ over the Runge-Kutta variant
$E^{\mathrm{R}}$, one immediate candidate is higher-order extensions
of $E^2$. In fact, the foundations are in place for such an
extension. These include the domain model
${\mathcal{W}}_{\mathcal{D}}$ of the current article, the domain of
Lipschitz functions developed
in~\parencite{EdalatLieutierPattinson:2013-MultiVar-Journal}, and an
extension of the Taylor's theorem, as presented
in~\cite[Corollary~2.14]{Edalat_Farjudian_Mohammadian_Pattinson:2nd_Order_Euler:2020:Conf}. To
obtain a full formulation for the $m$-th order variant $E^m$ of $E^2$,
the (non-scalar) differential equation $y'=f(y)$ must be
differentiated $(m-1)$ times, and the right hand side written solely
in terms of the vector field $f$ and its derivatives. This may be done
using, {\eg}, the well-known Fa{\`a} di Bruno
formula~\parencite{Constantine_Savits:Faa_di_Bruno:1996}. Although of
practical importance, this extension will require careful handling of
lengthy formulas, and may not require any significant theoretical
novelty.

A similar extension may also be considered for higher-order
Runge-Kutta methods. Again, the foundation is available, not just
based on the domain model of the current paper and that
of~\parencite{EdalatLieutierPattinson:2013-MultiVar-Journal}, but also
based on the catalogue of validated Runge-Kutta methods available
in~\parencite{Marciniak:Selected_Interval_Methods:2009}.

We have proved that the operator~$E^2$ has a second-order convergence
when the vector field of the \ac{IVP} is continuously
differentiable. The definition of the operator, however, only requires
the vector field to be Lipschitz continuous. One direction for further
investigation is convergence analysis of the operator~$E^2$ under the
assumption that the vector field is Lipschitz continuous but not
differentiable. As we have demonstrated, there is strong evidence that
the rate of convergence is not jeopardized by non-differentiability of
the Lipschitz vector field. This is potentially significant because,
to the best of our knowledge, all the higher-order validated methods
in the literature require the vector field to be at least (once)
continuously differentiable.

We presented the construction of a continuous domain for function
spaces $[X \to D]$, for any topological space $X$ and bounded-complete
continuous domain $D$. When $X$ is core-compact, the \ac{DCPO}
$[X \to D]$ is continuous. Thus, in practice, our general construction
may be more relevant when $X$ is not core-compact. We exemplified this
claim via the concrete case of differential equation solving with
temporal discretization, where the relevant topology is the upper
limit topology, which is not core-compact. It will be interesting to
see if the construction is useful for other applications as well,
{\eg}, stochastic processes with right-continuous jumps.

Probabilistic power domains have been used previously for analysis of
stochastic processes. Specifically,
in~\parencite{Edalat:Domain_stochastic_processes:LiCS:1995},
domain-theoretic models of finite-state discrete stochastic processes
have been introduced, while
in~\parencite{Bilokon_Edalat:Domain_Brownian:2017} a more general
approach has been taken for continuous time, continuous space
stochastic processes. The construction of the current article can be
modified for stochastic processes with right-continuous jumps. To be
more precise, while the upper limit topology had to be used for
temporal discretization in differential equation solving, the suitable
topology for analysis of processes with right-continuous jumps is the
lower limit topology, which is not core-compact either. Furthermore,
while the rational $P$-functions satisfying~\eqref{eq:rational_P_fun}
are left continuous with right limits, in the study of stochastic
processes, one must work with the so-called c{\`a}dl{\`a}g functions,
which are right-continuous with left
limits~\cite[Chapter~3]{Billingsley:Converg_Prob_Meas:1999}.


\printacronyms

\printbibliography

\end{document}